\documentclass[a4paper,12pt]{article} 
\usepackage{amssymb}
\usepackage{amscd}
\usepackage{amsmath} 
\usepackage[dvipdfmx]{graphicx}
\usepackage{latexsym} 
\usepackage{enumerate}
\usepackage[usenames]{color}
\usepackage{lscape}

\usepackage{theorem}
\theoremstyle{plain}
\theorembodyfont{\itshape}
\newtheorem{theorem}{Theorem}[section]
\newtheorem{proposition}[theorem]{Proposition}
\newtheorem{lemma}[theorem]{Lemma}

\theorembodyfont{\rmfamily}
\newtheorem{definition}[theorem]{Definition}
\newtheorem{remark}[theorem]{Remark}

\newtheorem{example}[theorem]{Example}

\def\ra{\mathrm{a}}
\def\rb{\mathrm{b}}
\def\rc{\mathrm{c}}
\def\rs{\mathrm{s}}
\def\rt{\mathrm{t}}

\def\rRe{\mathrm{Re}} 
\def\A{\mathbb{A}}
\def\C{\mathbb{C}}
 
\def\N{\mathbb{N}}

\def\Q{\mathbb{Q}}
\def\R{\mathbb{R}}
\def\Z{\mathbb{Z}}
\def\cS{\mathcal{S}}
\def\ve{\varepsilon}
\def\vG{\varGamma}
\def\vD{\varDelta}
\def\ba{\mbox{\boldmath $a$}}
\def\sba{\mbox{\boldmath ${\scriptstyle a}$}}

\def\bc{\mbox{\boldmath $c$}}

\def\be{\mbox{\boldmath $e$}}
\def\bj{\mbox{\boldmath $j$}}
\def\bk{\mbox{\boldmath $k$}}
\def\sbk{\mbox{\boldmath ${\scriptstyle k}$}}
\def\bl{\mbox{\boldmath $l$}}
\def\bp{\mbox{\boldmath $p$}}

\def\bt{\mbox{\boldmath $t$}}
\def\sbt{\mbox{\boldmath ${\ss t}$}}

\def\bv{\mbox{\boldmath $v$}}

\def\bal{\mbox{\boldmath $\alpha$}}
\def\bbeta{\mbox{\boldmath $\beta$}}
\def\blambda{\mbox{\boldmath $\lambda$}}
\def\sblambda{\mbox{\boldmath ${\scriptstyle \lambda}$}}
\def\bmu{\mbox{\boldmath $\mu$}}
\def\sbmu{\mbox{\boldmath ${\scriptstyle \mu}$}}
\def\bsigma{\mbox{\boldmath $\sigma$}}
\def\sbsigma{\mbox{\boldmath ${\scriptstyle \sigma}$}}
\def\btau{\mbox{\boldmath $\tau$}}
\def\sbtau{\mbox{\boldmath ${\scriptstyle \tau}$}} 
 
\def\0{\mbox{\boldmath $0$}}
\def\1{\mbox{\boldmath $1$}}
\def\2{\mbox{\boldmath $2$}}
\def\3{\mbox{\boldmath $3$}}
\def\cfL{\mbox{\Large \boldmath $\mathrm{K}$}}

\def\Phim{\Phi_{\scriptstyle \mathrm{max}}}
\def\Mm{\mathrm{M{\scriptstyle ax}}}

\def\ds{\displaystyle}
\def\ss{\scriptstyle} 
\def\ts{\textstyle}
\title{\bf Contiguous Relations, Laplace's Methods \\ and  
Continued Fractions for $\mbox{\boldmath ${}_3F_2(1)$}$\thanks{MSC (2010): 
Primary 33C20; Secondary 30B70, 30E10. Keywords: generalized hypergeometric series; 
continued fraction; contiguous relation; Laplace's method.}}  
\author{Akihito Ebisu\thanks{Faculty of Information and Computer Science, 
Chiba Institute of Technology 2-1-1, Shibazono, Narashino, Chiba, 275-0023, Japan. 
{\tt akihito.ebisu@p.chibakoudai.jp}} \ and Katsunori Iwasaki\thanks{Department of 
Mathematics, Hokkaido University, Kita 10, Nishi 8, Kita-ku, Sapporo 060-0810 Japan.  
{\tt iwasaki@math.sci.hokudai.ac.jp} (Corresponding author)}}
\date{September 1, 2017} 
\begin{document}
%%%%%%%%%%%%%%%%%%%%%%%%%%%%%%%%%%%%%%%%%%%%%%%%%%%%%%%%%%%%%%%%%%%%%%%%
\maketitle
\begin{abstract}
Using contiguous relations we construct an infinite number of   
continued fraction expansions for ratios of generalized 
hypergeometric series ${}_3F_2(1)$. 
We establish exact error term estimates for their approximants 
and prove their rapid convergences.    
To do so we develop a discrete version of Laplace's method for hypergeometric 
series in addition to the use of ordinary (continuous) Laplace's method for 
Euler's hypergeometric integrals.  
\end{abstract} 
%%%%%%%%%%%%%%%%%%%%%%%%%%%%%%% sec:intro %%%%%%%%%%%%%%%%%%%%%%%%%%%%%%%
\section{Introduction} \label{sec:intro}
%%%%%%%%%%%%%%%%%%%%%%%%%%%%%%%%%%%%%%%%%%%%%%%%%%%%%%%%%%%%%%%%%%%%%%%%%
In 1813 Gauss \cite{Gauss} introduced a general continued fraction that represents  
the ratio of two ${}_2F_1$ hypergeometric functions. 
It is interesting because it contains a variety of continued fraction expansions of 
several important elementary functions and some of more transcendental ones.   
In 1901 Van Vleck \cite{VanVleck} established a general result on its convergence.   
Gauss's continued fraction is derived from a three-term contiguous relation for 
${}_2F_1$. 
In 1956, using other contiguous relations, Frank \cite{Frank} constructed some 
more (eight or so) continued fractions of a similar sort and discussed their 
convergence.  
In 2005 Borwein, Choi and Pigulla \cite{BCP} obtained an explicit bound for 
the error term in certain special cases of the Gauss continued fraction.   
In 2011, based on Gauss's continued fraction and other means, Colman, Cuyt 
and Van Deun \cite{CCD} developed an efficient algorithm for the validated 
high-precision computation of certain ${}_2F_1$ functions.    
%%%%%
\par
%%%%% 
The generalized hypergeometric series of 
unit argument ${}_3F_2(1)$ also admits three-term contiguous relations, 
among which the  basic twelve relations were found by Wilson \cite{Wilson};  
see also Bailey \cite{Bailey2}.   
Thus it is feasible and interesting to discuss or utilize allied continued fractions 
for ${}_3F_2(1)$. 
For instance, Zhang \cite{Zhang} used contiguous relations for ${}_3F_2(1)$ 
to give new proofs of three of Ramanujan's elegant continued fractions for products 
and quotients of gamma functions, namely, entries 34, 36 and 39 in 
Ramanujan's second notebook \cite[Chapter 12]{Ramanujan}, or in its 
corrected version by Berndt, Lamphere and Wilson \cite{BLW}.    
In a similar vein, Denis and Singh \cite{DS}  dealt with entries 25 and 33 
of the same notebook.   
%%%%%
\par
%%%%% 
To give a further motivation for ${}_3F_2(1)$ continued fractions, we look 
at the special case in which one of the numerator parameters, say $a_0$, 
is equal to one:   
%%%%%%%%%%%%%%%%%%%%%%%%%%%% eqn:3F2-i %%%%%%%%%%%%%%%%%%%%%%%%%%%%%%%%%
\begin{equation} \label{eqn:3F2-i}
{}_3F_2\! 
\begin{pmatrix} 
1, & a_1, & a_2 \\ 
     & b_1, & b_2    
\end{pmatrix} := \sum_{j=0}^{\infty} 
\dfrac{(a_1; \, j) \, (a_2; \, j)}{(b_1; \, j) \, (b_2; \, j)},   
\qquad (a; \, j) := \dfrac{\vG( a+j )}{ \vG( a )},    
\end{equation}
%%%%%%%%%%%%%%%%%%%%%%%%%%%%%%%%%%%%%%%%%%%%%%%%%%%%%%%%%%%%%%%%%%%%%%%%%
where $\vG( a )$ is Euler's gamma function.  
%%%%%%%%%%%%%%%%%%%%%%%%%%%%%%% tab:3F2-i %%%%%%%%%%%%%%%%%%%%%%%%%%%%%%
\begin{table}[t]
\begin{align*}
{}_3F_2\! 
\begin{pmatrix} 
1, & 1, & 1 \\ 
   & 2, & 2    
\end{pmatrix} &= \zeta(2) = \dfrac{\pi^2}{6} \qquad \mbox{(Euler 1735)}, \\[2mm]
{}_3F_2\! 
\begin{pmatrix} 
1, & \frac{1}{4}, & \frac{1}{2} \\[2mm] 
   & \frac{5}{4}, & \frac{3}{2}    
\end{pmatrix} &= \dfrac{\pi + 2 \log 2}{4}, \\[2mm]
{}_3F_2\! 
\begin{pmatrix} 
1, & \frac{1}{4}, & \frac{1}{4} \\[2mm] 
   & \frac{5}{4}, & \frac{5}{4}    
\end{pmatrix} &= \dfrac{\pi^2 + 8 \, G}{16} \qquad 
\mbox{($G := \displaystyle \sum_{k=0}^{\infty} \frac{(-1)^k}{(2k+1)^2}$ is Catalan's constant)}, \\[2mm]
{}_3F_2\! 
\begin{pmatrix} 
1, & \frac{5}{6}, & \frac{4}{3} \\[2mm] 
   &  \frac{3}{2}, & \frac{11}{6}
\end{pmatrix} &= \frac{5}{2} \, \sqrt{3} \log(2+\sqrt{3}) \qquad \mbox{(Watson 1918)}, \\[2mm]
{}_3F_2\!
\begin{pmatrix}
1, & 1, & \frac{3}{4} \\[2mm]
   & \frac{7}{6}, & \frac{11}{6}
\end{pmatrix}
&=\frac{2^{\frac{1}{2}} \cdot 5}{ 3^{ \frac{5}{4} } }
\log\left( \frac{3^{\frac{5}{4}}-3^{\frac{3}{4}} + 2^{\frac{1}{2} } }{ 3^{\frac{5}{4}} 
-3^{\frac{3}{4}}-2^{\frac{1}{2} }}\right)  -\frac{2^{\frac{3}{2} } \cdot 5 }{ 3^{\frac{5}{4} } } \arccos
\left(\frac{3^{\frac{5}{4}}+3^{\frac{3}{4}}}{2\sqrt{5+ 3^{\frac{3}{2} } }}\right) \\[1mm]
& \hspace{3cm} \mbox{(Asakura, Otsubo and Terasoma 2016 \cite{AOT})},\\[2mm]
{}_3F_2\! 
\begin{pmatrix} 
1, & 1, & \frac{2}{3} \\[2mm] 
   & \frac{7}{6}, & \frac{11}{6}    
\end{pmatrix} 
&= \frac{5}{ 2^{\frac{5}{3}} \cdot 3^{\frac{1}{2}} } \log \left( 
\frac{ 2^{\frac{2}{3}} + 2^{\frac{1}{3}} -1 + 3^{\frac{1}{2}} }{ 2^{\frac{2}{3}} + 2^{\frac{1}{3}} -1 - 3^{\frac{1}{2}} } \right) 
- \frac{5}{ 2^{\frac{2}{3}} \cdot 3 } \arctan \left( \frac{3}{3 \cdot 2^{\frac{2}{3}} + 2^{\frac{1}{2}} + 3} \right) \\[1mm]
& \hspace{3cm} \mbox{(Yabu 2017 \cite{Yabu})}, \\[2mm]
{}_3F_2\! 
\begin{pmatrix} 
1, & a, & b \\ 
   & 1, & c    
\end{pmatrix} &= {}_2F_1\! 
\begin{pmatrix} 
a, & b \\ 
   & c    
\end{pmatrix} = \dfrac{\vG( c ) \, \vG( c-a-b ) }{\vG( c-a ) \, \vG( c-b )} 
\qquad \mbox{(Gauss 1812)}, \\[2mm]
{}_3F_2\! 
\begin{pmatrix} 
1, & a, & b \\ 
   & 2, & c    
\end{pmatrix} &= \dfrac{c-1}{(a-1)(b-1)} 
\left\{ \dfrac{\vG( c-1 ) \, \vG( c+1-a-b ) }{\vG( c-a ) \, \vG( c-b )} - 1\right\}.      
\end{align*}
\caption{Some special evaluations of the series ${}_3F_2(1, a_1, a_2; b_1, b_2)$.} 
\label{tab:3F2-i}
\end{table} 
%%%%%%%%% 
This series is well defined and non-terminating if   
%%%%%%%%%%%%%%%%%%%%%%%%%%%%%%% eqn:wd-i %%%%%%%%%%%%%%%%%%%%%%%%%%%%%%%%%%
\begin{equation} \label{eqn:wd-i} 
a_1, \, a_2, \, b_1, \, b_2 \, \not\in \, \Z_{\le 0},      
\end{equation}
%%%%%%%%%%%%%%%%%%%%%%%%%%%%%%%%%%%%%%%%%%%%%%%%%%%%%%%%%%%%%%%%%%%%%%%%%
in which case the series is absolutely convergent if and only if 
%%%%%%%%%%%%%%%%%%%%%%%%%%%%%%% eqn:saal-i %%%%%%%%%%%%%%%%%%%%%%%%%%%
\begin{equation} \label{eqn:saal-i}
\rRe \, s > 0, \qquad s := b_1+b_2-a_1-a_2-1.  
\end{equation}
%%%%%%%%%%%%%%%%%%%%%%%%%%%%%%%%%%%%%%%%%%%%%%%%%%%%%%%%%%%%%%%%%%%%%%%%% 
\par
%%%%%%
This class of infinite sums are interesting because they contain a lot of 
special evaluations, some of which are presented in Table \ref{tab:3F2-i}.  
Therefore it is important to establish a general framework for the 
precise and efficient computations of the series \eqref{eqn:3F2-i}. 
Naturally, our approach here is based on three-term contiguous 
relations and allied continued fractions. 
As an illustration of a more general story to be developed in this 
article, we shall present a continued fraction expansion of the 
series \eqref{eqn:3F2-i} with an exact error term estimate 
for its approximants that exhibits an exponentially fast convergence 
(see Theorem \ref{thm:ete-i}).     
%%%%%%
\par
%%%%% 
To state Theorem \ref{thm:ete-i},  
%%%%%%%%%%%%%%%%%%%%%%%%%%%%%%%% tab:pnd %%%%%%%%%%%%%%%%%%%%%%%%%%%%%%%%%
\begin{table}[t]
%%%%%%%%%%%%%%%%%%%%%%%%%%%%%%%% eqn:q-i %%%%%%%%%%%%%%%%%%%%%%%%%%%%%%%%%
\begin{alignat*}{2} %\label{eqn:q-i}
q_0(n) &:= \dfrac{(3 n+b_1-1)(3 n+b_2-1)-(2 n)(2 n+ a_2)}{(2 n)(2 n + a_1 -1)} 
\qquad & (n &\ge 1),  \\[2mm]
q_1(n) &:= \dfrac{(3 n + b_1)(3 n + b_2)-(2 n+1)(2 n+ a_1)}{(2 n + a_1)(2 n + a_2)} 
\qquad & (n &\ge 0),  \\[2mm]
q_2(n) &:= 
\dfrac{(3 n+b_1+1)(3 n+ b_2+1)-(2 n+ a_1+1)(2 n+ a_2+1)}{(2 n+1)(2 n + a_2 + 1)} 
\qquad & (n &\ge 0),  \\[2mm] 
r_0(n) &:= - \dfrac{(n+b_1-a_2-1)(n+b_2-a_2-1)}{(2 n-1)(2 n+a_2-1)} 
\qquad & (n &\ge 1), \\[2mm]
r_1(n) &:= - \dfrac{(n+b_1-1)(n+b_2-1)}{(2 n)(2 n+a_1-1)} 
\qquad & (n &\ge 1),  \\[2mm]
r_2(n) &:= - \dfrac{(n+b_1-a_1)(n+b_2-a_1)}{(2 n+a_1)(2 n+a_2)} 
\qquad & (n &\ge 0).   
\end{alignat*}
\caption{Partial denominators and numerators of the continued 
fraction \eqref{eqn:cf-i}.} 
\label{tab:pdn-cf}
\end{table} 
%%%%%%%%%%%%%%%%%%%%%%%%%%%%%%%%%%%%%%%%%%%%%%%%%%%%%%%%%%%%%%%%%%%%%%%
let $\{q(n)\}_{n=0}^{\infty}$ and $\{r(n)\}_{n=0}^{\infty}$ be infinite sequences 
defined by    
%%%%%%%%%%%%%%%%%%%%%%%%%%% eqn:qr-i %%%%%%%%%%%%%%%%%%%%%%%%%%%%%%%%%%% 
\begin{equation} \label{eqn:qr-i} 
q(n) := q_i((n-i)/3), \quad r(n) := r_i((n-i)/3), \quad \mbox{for} \quad 
n \equiv i \mod 3, \quad i = 0, 1, 2,     
\end{equation}
%%%%%%%%%%%%%%%%%%%%%%%%%%%%%%%%%%%%%%%%%%%%%%%%%%%%%%%%%%%%%%%%%%%%%%
where $q_i(n)$ and $r_i(n)$ are given by formulas in Table \ref{tab:pdn-cf} 
and $q_0(0) := 1$, $r_0(0) := 1$, $r_1(0) = -1$.  
The modulo $3$ structure in \eqref{eqn:qr-i} is the reflection of a 
$\Z_3$-symmetry in the relevant contiguous relations 
(see \S\ref{subsec:cr}).  
Under condition \eqref{eqn:wd-i}, all the $q(n)$ and $r(n)$ have  
non-vanishing denominators, while all the $r(n)$ have non-vanishing 
numerators if and only if the parameters satisfy     
%%%%%%%%%%%%%%%%%%%%%%%%%%% eqn:ba-i %%%%%%%%%%%%%%%%%%%%%%%%%%%%%%%
\begin{equation} \label{eqn:ba-i}
b_i - a_j \, \not\in \, \Z_{\le 0}, \qquad i, j = 1, 2.      
\end{equation}
%%%%%%%%%%%%%%%%%%%%%%%%%%%%%%%%%%%%%%%%%%%%%%%%%%%%%%%%%%%%%%%%%%%%%
Thus the (formal) infinite continued fraction  
%%%%%%%%%%%%%%%%%%%%%%%%%%%%%%% eqn:cf-i %%%%%%%%%%%%%%%%%%%%%%%%%%%%
\begin{equation} \label{eqn:cf-i}
\overset{\infty}{\underset{j=0}{\cfL}} \,\, \frac{r(j)}{q(j)} 
:= \frac{r(0)}{q(0)} \,\, \mathop{}_{+} \,\, \frac{r(1)}{q(1)} 
\,\, \mathop{}_{+} \,\, \frac{r(2)}{q(2)} \, \, \mathop{}_{+} \,\, 
\mathop{}_{\textstyle \cdots\cdots} 
\end{equation}
%%%%%%%%%%%%%%%%%%%%%%%%%%%%%%%%%%%%%%%%%%%%%%%%%%%%%%%%%%%%%%%%%%%%%%%% 
makes sense, provided that the conditions \eqref{eqn:wd-i} and \eqref{eqn:ba-i} 
are satisfied.  
%%%%%%%%%%%%%%%%%%%%%%%%%%%%%%% thm:ete-i %%%%%%%%%%%%%%%%%%%%%%%%%%%%%%%%
\begin{theorem} \label{thm:ete-i}
If conditions \eqref{eqn:wd-i}, \eqref{eqn:saal-i} and 
\eqref{eqn:ba-i} are fulfilled then continued fraction 
\eqref{eqn:cf-i} converges to series \eqref{eqn:3F2-i} 
exponentially fast and there exists an exact error term estimate 
for its approximants:    
%%%%%%%%%%%%%%%%%
\begin{equation*} 
{}_3F_2\! 
\begin{pmatrix} 
1, & a_1, & a_2 \\ 
     & b_1, & b_2    
\end{pmatrix} - 
\overset{n}{\underset{j=0}{\cfL}} \,\, \frac{r(j)}{q(j)} 
= C\! 
\begin{pmatrix} 
a_1, & a_2 \\ 
b_1, & b_2 
\end{pmatrix} \cdot 
\dfrac{(3 n)^{ \frac{1}{2} - s}}{2^{2 n + a_1 + a_2}} \cdot 
\left\{1 + O( n^{- \frac{1}{2} } ) \right\},    
\end{equation*}
%%%%%%%%%%%%%%%%%%%%%%%%%%%%%%%%%%%%%%%%%%%%%%%%%%%%%%%%%%%%%%%%%%%
as $n \to +\infty$, where the constant $C(a_1, a_2; b_1, b_2)$ is given by   
%%%%%%%%%%%%%%%%%%%%%%%%%%%%%%%
\begin{equation*} 
C\! 
\begin{pmatrix} 
a_1, & a_2 \\ 
b_1, & b_2 
\end{pmatrix} := 
\dfrac{ \pi^{\frac{3}{2}} \, \vG(b_1) \, \vG(b_2) \, \vG^2(s) }{ \vG(a_1) \, 
\vG(a_2) \, \vG( b_1-a_1 ) \, \vG( b_1 -a_2 ) \, \vG( b_2 -a_1 ) \, \vG( b_2 -a_2 ) }. 
\end{equation*}
%%%%%%%%%%%%%%%%%%%%%%%%%%%%%%%%%%%%%%%%%%%%%%%%%%%%%%%%%%%%%%%%%%%%%%%%
\end{theorem}
%%%%%%%%%%%%%%%%%%%%%%%%%%%%%%%%%%%%%%%%%%%%%%%%%%%%%%%%%%%%%%%%%%%%%%%%
\par
%%%%%% 
Theorem \ref{thm:ete-i} is only a corollary to a specific example of  
infinitely many continued fractions with exact error 
estimates we shall establish in 
Theorems \ref{thm:cf-straight} and \ref{thm:cf-cyclic} (see Example \ref{ex:1}).  
To generate infinitely many continued fractions, we naturally need infinitely many 
contiguous relations, so we then need a general theory, beyond the scopes of 
Bailey \cite{Bailey2} and Wilson \cite{Wilson}, that presides over all contiguous 
relations for ${}_3F_2(1)$.  
Our previous paper \cite{EI} develops such a theory and the present article relies 
substantially on the main results of that paper.          
%%%%%%%%%%%%%%%%%%%%%%%%%%%%%%% sec:contig %%%%%%%%%%%%%%%%%%%%%%%%%%%%%%
\section{Contiguous and Recurrence Relations} \label{sec:contig}
%%%%%%%%%%%%%%%%%%%%%%%%%%%%%%%%%%%%%%%%%%%%%%%%%%%%%%%%%%%%%%%%%%%%%%%%
The hypergeometric series of unit argument ${}_3F_2(1)$ with full five 
parameters is defined by 
%%%%%%%%%%%%%%%%%%%%%%%%%%%%%%% 
\begin{equation*} 
{}_3F_2\! 
\begin{pmatrix} 
a_0, & a_1, & a_2 \\ 
     & b_1, & b_2    
\end{pmatrix} := \sum_{j=0}^{\infty} 
\dfrac{(a_0; \, j) \, (a_1; \, j) \, (a_2; \, j)}{(1; \, j) \, 
(b_1; \, j) \, (b_2; \, j)}.    
\end{equation*}
%%%%%%%%%%%%%%%%%%%%%%%%%%%%%
With the notation $\ba = (a_0, a_1, a_2; a_3, a_4) = (a_0, a_1, a_2; b_1, b_2)$ 
this series is often denoted by ${}_3F_2(\ba)$. 
It is well defined and non-terminating as a formal sum if $\ba$ satisfies  
%%%%%%%%%%%%%%%%%%%%%%%%%%%%%%% eqn:WD %%%%%%%%%%%%%%%%%%%%%%%%%%%%%%%%%
\begin{equation} \label{eqn:WD} 
a_0, \, a_1, \, a_2, \, b_1, \, b_2 \, \not\in \, \Z_{\le 0},        
\end{equation}
%%%%%%%%%%%%%%%%%%%%%%%%%%%%%%%%%%%%%%%%%%%%%%%%%%%%%%%%%%%%%%%%%%%%%%%%
in which case ${}_3F_2(\ba)$ is absolutely convergent if and only if 
%%%%%%%%%%%%%%%%%%%%%%%%%%%%%%% eqn:3saal2 %%%%%%%%%%%%%%%%%%%%%%%%%%%%%%
\begin{equation} \label{eqn:3saal2}
\rRe \, s(\ba)  > 0, \qquad s(\ba) := b_1+b_2-a_0-a_1-a_2,   
\end{equation}
%%%%%%%%%%%%%%%%%%%%%%%%%%%%%%%%%%%%%%%%%%%%%%%%%%%%%%%%%%%%%%%%%%%%%%%%
where $s(\ba)$ is called the {\sl Saalsch\"{u}tzian index} for ${}_3F_2(\ba)$.  
We say that $\ba$ is {\sl balanced} if $s(\ba) = 0$.  
%%%%%%
\par
%%%%%% 
In order to discuss contiguous relations, however, we find it more convenient 
in many respects to replace ${}_3F_2(\ba)$ by the renormalized 
hypergeometric series defined by 
%%%%%%%%%%%%%%%%%%%%%%%%%%%%%%%
\begin{equation*} 
{}_3f_2(\ba) := \sum_{j=0}^{\infty} 
\dfrac{\vG(a_0+j) \vG(a_1+j) \vG(a_2+j) }{\vG(1+j) \vG(b_1+j) \vG(b_2+j)}.     
\end{equation*}
%%%%%%%%%%%%%%%%%%%%%%%%%%%%%%
This latter series is well defined and non-terminating as a formal sum, 
whenever 
%%%%%%%%%%%%%%%%%%%%%%%%%%%%%%%
\begin{equation*}
a_0, \, a_1, \, a_2 \, \not\in \, \Z_{\le 0} \qquad 
\mbox{(compare this with condition \eqref{eqn:WD})}, 
\end{equation*}
%%%%%%%%%%%%%%%%%%%%%%%%%%%%%
in which case series ${}_3f_2(\ba)$ is absolutely convergent 
if and only if \eqref{eqn:3saal2} is satisfied. 
Note that   
%%%%%%%%%%%%%%%%%%%%%%%%%%%%%%% eqn:Fvsf %%%%%%%%%%%%%%%%%%%%%%%%%%%%%%%%%
\begin{equation} \label{eqn:Fvsf}
{}_3f_2(\ba) 
= \dfrac{\vG( a_0 ) \, \vG( a_1 ) \, \vG( a_2 )}{\vG( b_1 ) \, \vG( b_2 )} 
\, {}_3F_2(\ba),  
\end{equation} 
%%%%%%%%%%%%%%%%%%%%%%%%%%%%%%%%%%%%%%%%%%%%%%%%%%%%%%%%%%%%%%%%%%%%%%%%
as long as both sides of equation \eqref{eqn:Fvsf} make sense.  
%%%%%%%%%%%%%%%%%%%%%%%%%%%%% subsec:cr %%%%%%%%%%%%%%%%%%%%%%%%%%%%%%%%%%
\subsection{Contiguous Relations} \label{subsec:cr}
%%%%%%%%%%%%%%%%%%%%%%%%%%%%%%%%%%%%%%%%%%%%%%%%%%%%%%%%%%%%%%%%%%%%%%%%
It follows from \cite[Theorem 1.1]{EI} that for any distinct integer vectors 
$\bk$, $\bl \in \Z^5$ different from $\0$ there exist {\sl unique} rational functions 
$u(\ba)$, $v(\ba) \in \Q(\ba)$ such that 
%%%%%%%%%%%%%%%%%%%%%%%%%%%%%%%% eqn:3tcr %%%%%%%%%%%%%%%%%%%%%%%%%%%%%%%%
\begin{equation} \label{eqn:3tcr}
{}_3f_2(\ba) = u(\ba) \cdot {}_3f_2(\ba + \bk) + v(\ba) \cdot 
{}_3f_2(\ba + \bl).      
\end{equation} 
%%%%%%%%%%%%%%%%%%%%%%%%%%%%%%%%%%%%%%%%%%%%%%%%%%%%%%%%%%%%%%%%%%%%%%%%
An identity of the form \eqref{eqn:3tcr} is called a {\sl contiguous relation} for 
${}_3f_2(1)$.        
An algorithm to calculate $u(\ba)$ and $v(\ba)$ explicitly is given in 
\cite[Recipe 5.4]{EI}. 
According to it, one calculates the connection matrix $A(\ba; \bk)$ as in 
\cite[formula (30)]{EI} and define $r(\ba; \bk) \in \Q(\ba)$ to be its $(1, 2)$-entry 
as in \cite[formula (33)]{EI}. 
One also calculates $r(\ba; \bl)$ as well as $r(\ba; \bl-\bk)$ in similar manners. 
If $\bk$ and $\bl$ are distinct then $r(\ba; \bl-\bk)$ is nonzero in 
$\Q(\ba)$ and the coefficients in \eqref{eqn:3tcr} are represented as 
%%%%%%%%%%%%%%%%%%%%%%%%%%%%%% eqn:uv %%%%%%%%%%%%%%%%%%%%%%%%%%%%%%%%%%
\begin{equation} \label{eqn:uv}
u(\ba) = \dfrac{r(\ba; \bl)}{ \det A(\ba;\bk) \cdot r(\ba+\bk; \bl-\bk)}, 
\qquad 
v(\ba) = - \dfrac{r(\ba; \bk)}{ \det A(\ba; \bk) \cdot r(\ba+\bk; \bl-\bk)},  
\end{equation}
%%%%%%%%%%%%%%%%%%%%%%%%%%%%%%%%%%%%%%%%%%%%%%%%%%%%%%%%%%%%%%%%%%%%%%
as in \cite[Proposition 5.3]{EI}, where according to \cite[formula (32)]{EI} one has    
%%%%%%%%%%%%%%%%%%%%%%%%%%%%% eqn:eps(a;k) %%%%%%%%%%%%%%%%%%%%%%%%%%%%%
\begin{equation} \label{eqn:eps(a;k)} 
\det A(\ba; \bk) = \dfrac{(-1)^{k_0+k_1+k_2} 
(s(\ba)-1; \, s(\bk)) \prod_{i=0}^2 (a_i; \, k_i) }{ \prod_{i=0}^2 \prod_{j=1}^2 
(b_j-a_i; \, l_j-k_i) }. 
\end{equation}
%%%%%%%%%%%%%%%%%%%%%%%%%%%%%%%%%%%%%%%%%%%%%%%%%%%%%%%%%%%%%%%%%%%%%%
\par
%%%%%%%
In order to formulate our main results in \S\ref{subsec:mrcf}, we need one more 
fact about the structure of $r(\ba; \bk)$ which is not discussed in \cite{EI}.    
Given a vector $\bk = (k_0,k_1,k_2; l_1, l_2) \in \Z^5$, let    
%%%%%%
\[
\langle \ba; \, \bk \rangle_{\pm} 
:= \prod_{i=0}^2 \prod_{j=1}^2 (b_j-a_i; \, (l_j-k_i)_{\pm}), 
\qquad   
|\!| \bk |\!|_+ := \sum_{i=0}^2 \sum_{j=1}^2 (l_j-k_i)_+,   
\]
%%%%%%
where $m_{\pm} := \max\{ \pm m, 0\}$.  
Note that $\prod_{i=0}^2 \prod_{j=1}^2 (b_j-a_i; \, l_j-k_i) = 
\langle \ba; \bk \rangle_+/\langle \ba + \bk; \bk \rangle_-$.   
%%%%%%%%%%%%%%%%%%%%%%% lem:factorial %%%%%%%%%%%%%%%%%%%%%%%%%%%%%%%%%
\begin{lemma} \label{lem:factorial}
For any nonzero vector $\bk \in \Z_{\ge 0}^5$ with $s(\bk) = 0$ there 
exists a nonzero polynomial $\rho(\ba; \bk) \in \Q[\ba]$ such that the rational 
function $r(\ba; \bk)$ can be written  
%%%%%%%%%%%%%%%%%%%%%%%%% eqn:rho %%%%%%%%%%%%%%%%%%%%%%%%%%%%%%%%%%%
\begin{equation} \label{eqn:rho}
r(\ba; \bk) = - \frac{ \{ s(\ba) -1 \} \rho(\ba; \bk) }{\langle \ba; \, \bk \rangle_+}, 
\qquad  \deg \rho(\ba; \bk)  \le |\!| \bk |\!|_+ - 2.   
\end{equation}
%%%%%%%%%%%%%%%%%%%%%%%%%%%%%%%%%%%%%%%%%%%%%%%%%%%%%%%%%%%%%%%%%%%%   
\end{lemma}
%%%%%%%%%%%%%%%%%%%%%%%%% begin proof %%%%%%%%%%%%%%%%%%%%%%%%%%%%%%%%
{\it Proof}.      
A nonzero polynomial $p(\ba) \in \Q[\ba]$ is said to be a {\sl denominator} 
of a rational function $r(\ba) \in \Q(\ba)$ if the product $p(\ba) \, r(\ba)$ 
becomes a polynomial.  
A denominator of the least degree, which is unique up to constant multiples, 
is referred to as the {\sl reduced denominator}. 
Any denominator is divisible by the reduced denominator in $\Q[\ba]$. 
A denominator of a matrix with entries in $\Q(\ba)$ is, by definition, 
a common denominator of those entries.   
%%%%%%%
\par
%%%%%% 
For $i = 0, 1, 2$,  $\mu = 1, 2$, let $\be_{\mu}^i := 
(\delta_{0i}, \delta_{1i}, \delta_{2i}; \delta_{1\mu}, \delta_{2\mu})$,   
where $\delta_{*\star}$ is Kronecker's delta.   
A vector of this form is said to be basic. 
A product of contiguous matrices in \cite[Table 2]{EI} yields   
%%%%%%
\[
A(\ba; \be_{\mu}^i) = \frac{1}{ (b_{\mu}-a_j)(b_{\mu}-a_k) }
\begin{pmatrix}
a_i(b_{\mu}-a_j-a_k) & s(\ba)-1 \\[1mm]
a_i a_j a_k & (a_i+1) b_{\mu} + a_j a_k - b_1 b_2 
\end{pmatrix},  
\]
%%%%%%%
where $\{i, j, k\} = \{0, 1, 2\}$. 
Any $\bk = (k_0, k_1, k_2; l_1, l_2) \in \Z_{\ge0}^5$ with $s(\bk) = 0$ 
admits a decomposition $\bk = \bv_l + \cdots + \bv_1$  
with each $\bv_i$ basic, so $A(\ba; \bk)$ can be computed by  
the chain rule    
%%%%%%%%%%%%%%%%%%%%%%%%% eqn:chain %%%%%%%%%%%%%%%%%%%%%%%%%%%%%%%
\begin{equation} \label{eqn:chain}
A(\ba; \bk) = A(\ba+\bv_{l-1}+\cdots+\bv_1; \bv_l) \cdots 
A(\ba+\bv_1; \bv_2) A(\ba; \bv_1). 
\end{equation}
%%%%%%%%%%%%%%%%%%%%%%%%%%%%%%%%%%%%%%%%%%%%%%%%%%%%%%%%%%%%%%%%%%
Thus $A(\ba; \bk)$ has a denominator each irreducible factor of which is 
of the form $b_{\mu} - a_i + \mbox{an integer}$. 
A factor of this form is said to be of type $b_{\mu}- a_i$ and the 
product of all factors of this type is referred to as the 
$b_{\mu}-a_i$ component of the denominator.     
%%%%%%%%%%%%%%%%%%%%%
\par\smallskip
%%%%%%%%%%%%%%%%%%%%%
{\bf Claim}. For each $i = 0,1,2$ and $\mu = 1,2$ the matrix 
$A(\ba; \bk)$ admits a denominator whose $b_{\mu}- a_i$ 
component is exactly the factorial function $(b_{\mu}-a_i; (l_{\mu}-k_i)_+)$. 
%%%%%%%%%%%%%%%%%%%%%
\par\smallskip
%%%%%%%%%%%%%%%%%%%%%
To show the claim we may assume $i = 0$ and $\mu = 1$ without 
loss of generality.  
%%%%%
\begin{enumerate}
\item If $m_0 := k_0-l_1 \ge 0$, then take the decomposition 
$\bk = l_1 \be_1^0 + m_0 \be_2^0 + k_1 \be_2^1 + k_2 \be_2^2$. 
\item If $m_1 := l_1-k_0 > 0$, then take the decomposition 
$\bk = k_{12} \be_2^1 + k_{22} \be_2^2 + k_0 \be_1^0 + k_{11} \be_1^1 
+ k_{21} \be_1^2$, 
where $k_{ij}$ are nonnegative integers such that $k_1 = k_{11}+k_{12}$, 
$k_2 = k_{21}+k_{22}$, $m_1 = k_{11}+k_{21}$ and $l_2 = k_{12} + k_{22}$;   
such $k_{ij}$ exist thanks to $\bk \in \Z_{\ge0}^5$ and $s(\bk) = 0$. 
\end{enumerate}
%%%%%
We use the fact that $A(\ba; m \be_{\mu}^i)$ has a denominator 
$(b_{\mu}-a_j; m)(b_{\mu}-a_k; m)$, where $\{i, j, k\} = \{0, 1, 2\}$, 
which follows by induction on $m \in \Z_{\ge0}$.  
In case (1) the decomposition of $\bk$ and the chain rule 
\eqref{eqn:chain} imply that $A(\ba; \bk)$ has a denominator 
without $b_1-a_0$ component. 
In case (2) the decomposition of $\bk$ leads to the product 
$A(\ba; \bk) = A_2(\ba; \bk) A_1(\ba; \bk)$ with 
%%%%%
\[
A_2(\ba; \bk) := A(\ba+k_{11} \be_1^1 + k_{21} \be_1^2; \, 
k_{12} \be_2^1 + k_{22} \be_2^2 + k_1 \be_1^0), \quad 
A_1(\ba; \bk) := A(\ba; k_{11} \be_1^1 + k_{21} \be_1^2). 
\]
%%%%%
Observe that $A_1(\ba; \bk)$ has a denominator whose $b_1-a_0$ 
component is $(b_1-a_0; k_{11}+k_{21}) = (b_1-a_0; m_1)$, 
while $A_2(\ba; \bk)$ has a denominator without 
$b_1-a_0$ component.  
So $A(\ba; \bk)$ has a denominator whose $b_1-a_0$ 
component is $(b_1-a_0; m_1)$. 
The claim is thus verified.  
%%%%%
\par
%%%%%
For each entry of $A(\ba; \bk)$ the Claim implies that for $i = 0,1,2$ and  
$\mu = 1, 2$ the $b_{\mu}-a_i$ component of its reduced denominator 
must divide the factorial  $(b_{\mu}-a_i; \, (l_{\mu}-k_i)_+)$, so the reduced denominator itself must divide the product 
$\langle \ba; \bk \rangle_+ = \prod_{i=0}^2 \prod_{\mu=1}^2 
(b_{\mu}-a_i; \, (l_{\mu}-k_i)_+)$. 
Thus one can take $\langle \ba; \bk \rangle_+$ as a denominator of 
$A(\ba; \bk)$. 
The {\sl index} of a rational function is the degree of its numerator 
minus that of its denominator. 
An induction on the length $l$ of product \eqref{eqn:chain} shows that 
the index $\le i-j$ for the $(i, j)$-entry of $A(\ba; \bk)$.  
Another induction shows that the $(1, 2)$-entry is divisible by $s(\ba)-1$. 
All these facts lead to expression \eqref{eqn:rho} for $r(\ba; \bk)$.  
\hfill $\Box$ 
%%%%%%%%%%%%%%%%%%%%%%%% end proof %%%%%%%%%%%%%%%%%%%%%%%%%%%%%%%%%%%%%
%%%%%%%%%%%%%%%%%%%%%%%%%%%%%%%% subsec:types %%%%%%%%%%%%%%%%%%%%%%%%%%
\subsection{Symmetry and Dichotomy} \label{subsec:types}
%%%%%%%%%%%%%%%%%%%%%%%%%%%%%%%%%%%%%%%%%%%%%%%%%%%%%%%%%%%%%%%%%%%%%%
Let $G = S_3 \times S_2$ be the group acting on $\ba = (a_0,a_1,a_2; b_1, b_2)$ by 
permuting $(a_0, a_1, a_2)$ and $(b_1, b_2)$ separately.  
It is obvious that ${}_3f_2(\ba)$ is invariant under this action, so that any element 
$\tau \in G$ transforms the contiguous relation \eqref{eqn:3tcr} into a second one      
%%%%%%%%%%%%%%%%%%%%%%%%%%%%%%%% eqn:3tcr2 %%%%%%%%%%%%%%%%%%%%%%%%%%%%%
\begin{equation} \label{eqn:3tcr2}
{}_3f_2(\ba) = {}^{\tau}\!u(\ba) \cdot {}_3f_2(\ba + \tau(\bk)) 
+ {}^{\tau}\!v(\ba) \cdot {}_3f_2(\ba + \tau(\bl)),   
\end{equation} 
%%%%%%%%%%%%%%%%%%%%%%%%%%%%%%%%%%%%%%%%%%%%%%%%%%%%%%%%%%%%%%%%%%%%%%%%
where ${}^{\tau}\!\varphi(\ba) := \varphi(\tau^{-1}(\ba))$ is 
the induced action of $\tau$ on a function $\varphi(\ba)$.  
%%%%%%
\par
%%%%%%
Take an element $\sigma \in G$ such that $\sigma^3$ is identity and set   
%%%%%%%%%%%%%%%%%%%%%%%%%%%%%%% eqn:lp %%%%%%%%%%%%%%%%%%%%%%%%%%%%%%%%%%%
\begin{equation} \label{eqn:lp}
\bl := \bk + \sigma(\bk), \qquad 
\bp := \bk + \sigma(\bl) = \bk + \sigma(\bk) + \sigma^2(\bk). 
\end{equation}
%%%%%%%%%%%%%%%%%%%%%%%%%%%%%%%%%%%%%%%%%%%%%%%%%%%%%%%%%%%%%%%%%%%%%%%%
Formula \eqref{eqn:3tcr2} with $\tau = \sigma$ followed by a shift  
$\ba \mapsto \ba + \bk$ yields        
%%%%%%%%%%%%%%%%%%%%%%%%%%%% eqn:prolong-k %%%%%%%%%%%%%%%%%%%%%%%%%%%%%%%
\begin{equation} \label{eqn:prolong-k}
{}_3f_2(\ba+\bk) = 
{}^{\sigma}\! u(\ba+\bk) \cdot {}_3f_2(\ba + \bl) 
+ {}^{\sigma}\! v(\ba + \bk) \cdot {}_3f_2(\ba + \bp),    
\end{equation}
%%%%%%%%%%%%%%%%%%%%%%%%%%%%%%%%%%%%%%%%%%%%%%%%%%%%%%%%%%%%%%%%%%%%%
and similarly formula \eqref{eqn:3tcr2} with $\tau = \sigma^2$ followed by 
another shift $\ba \mapsto \ba + \bl$ gives       
%%%%%%%%%%%%%%%%%%%%%%%%%%%% eqn:prolong-l %%%%%%%%%%%%%%%%%%%%%%%%%%%%%
\begin{equation} \label{eqn:prolong-l}
{}_3f_2(\ba+\bl) = 
{}^{\sigma^2}\! u(\ba+\bl) \cdot {}_3f_2(\ba + \bp) 
+ {}^{\sigma^2}\! v(\ba + \bl) \cdot {}_3f_2(\ba + \bp + \bk).  
\end{equation} 
%%%%%%%%%%%%%%%%%%%%%%%%%%%%%%%%%%%%%%%%%%%%%%%%%%%%%%%%%%%%%%%%%%%%
\par
%%%%%%
If $\bk$ is nonzero, nonnegative $\bk \in \Z_{\ge0}^5$ and balanced $s(\bk) = 0$, 
then so are $\bl-\bk = \sigma(\bk)$ and $\bl$ by definition \eqref{eqn:lp}, hence 
Lemma \ref{lem:factorial} applies not only to $\bk$ but also to $\sigma(\bk)$ 
and $\bl$.    
Putting formulas \eqref{eqn:eps(a;k)} and \eqref{eqn:rho} for these vectors 
into formula \eqref{eqn:uv} we have  
%%%%%%%%%%%%%%%%%%%%%%%%%%%% eqn:uv2 %%%%%%%%%%%%%%%%%%%%%%%%%%%%%%%% 
\begin{subequations} \label{eqn:uv2}
\begin{align}
u(\ba) &= \dfrac{(-1)^{k_0+k_1+k_2} \cdot \rho(\ba; \bl) \cdot 
\langle \ba; \bk\rangle_+ \cdot 
\langle \ba+\bk; \sigma(\bk) \rangle_+ }{ \rho(\ba+\bk; \sigma(\bk)) \cdot 
\langle \ba; \bl \rangle_+ \cdot \langle \ba+\bk; \bk \rangle_-  
\prod_{i=0}^2 (a_i; \, k_i)  },  
\label{eqn:uv2u} \\[2mm]
v(\ba) &= - \dfrac{(-1)^{k_0+k_1+k_2} \cdot \rho(\ba; \bk) \cdot  
\langle \ba+\bk; \sigma(\bk) \rangle_+ }{ \rho(\ba+\bk; \sigma(\bk)) \cdot 
\langle \ba+\bk; \bk \rangle_- \prod_{i=0}^2 (a_i; \, k_i)  },  
\label{eqn:uv2v} 
\end{align}
\end{subequations}
%%%%%%%%%%%%%%%%%%%% def:admissible %%%%%%%%%%%%%%%%%%%%%%%%%%%%%%%%%%%
\begin{definition} \label{def:admissible} 
For any nonzero vector $\bk \in \Z_{\ge0}^5$ with $s(\bk) = 0$ 
we consider two cases.         
\begin{enumerate} 
\item The case is said to be {\sl of  straight type} when $\sigma$ is identity, 
$\bl = 2 \bk$ and $\bp = 3 \bk$.     
\item The case is said to be {\sl of twisted type} when $\sigma$ is a cyclic 
permutation of  the upper parameters $(a_0, a_1, a_2)$ that acts    
on the lower parameters $(b_1, b_2)$ trivially,        
%%%%%%%%%%%%%%%%%%%%%%%%%% eqn:adm-kp %%%%%%%%%%%%%%%%%%%%%%%%%%%%%%%%%
\begin{equation}  \label{eqn:adm-kp} 
\bk = 
\begin{pmatrix} 
k_0, & k_1, & k_2 \\[1mm] 
     & l_1, & l_2 
\end{pmatrix}, 
\qquad 
\bp = 
\begin{pmatrix} 
p, & p, & p \\[1mm]
 & 3 l_1, & 3 l_2 
\end{pmatrix}, 
\end{equation}
%%%%%%%%%%%%%%%%%%%%%%%%%%%%%%%%%%%%%%%%%%%%%%%%%%%%%%%%%%%%%%%%%%%%%
with $p := k_0+k_1+k_2 = l_1 + l_2$, and 
if $\sigma(a_0,a_1,a_2; b_1, b_2) = (a_{\lambda}, a_{\mu}, a_{\nu}; b_1, b_2)$, then  
%%%%%%%%%%%%%%%%%%%%%%%%%%% eqn:cyclic %%%%%%%%%%%%%%%%%%%%%%%%%%%%%%%%% 
\begin{equation} \label{eqn:cyclic}   
\bl = 
\begin{pmatrix} 
k_0 + k_{\lambda}, & k_1 + k_{\mu},  & k_2 + k_{\nu} \\[2mm] 
           & 2 l_1,      &  2 l_2  
\end{pmatrix},  
\end{equation}  
%%%%%%%%%%%%%%%%%%%%%%%%%%%%%%%%%%%%%%%%%%%%%%%%%%%%%%%%%%%%%%%%%%%
where the index triple $(\lambda, \mu, \nu)$ is either $(2, 0, 1)$ or $(1, 2, 0)$.  
\end{enumerate}
\end{definition}
%%%%%%%%%%%%%%%%%%%%%%%%%%%%%%%%%%%%%%%%%%%%%%%%%%%%%%%%%%%%%%%%%%%
%%%%%%%%%%%%%%%%%%%%%%%%%% subsec:rr %%%%%%%%%%%%%%%%%%%%%%%%%%%%%%
\subsection{Recurrence Relations} \label{subsec:rr}  
%%%%%%%%%%%%%%%%%%%%%%%%%%%%%%%%%%%%%%%%%%%%%%%%%%%%%%%%%%%%%%%%%
In the situation of Definition \ref{def:admissible}, the shifts $\ba \mapsto 
\ba + n \bp$, $n \in \Z_{\ge 0}$, in the contiguous relation \eqref{eqn:3tcr} 
and its companions \eqref{eqn:prolong-k} and \eqref{eqn:prolong-l} induce 
a system of recurrence relations     
%%%%%%%%%%%%%%%%%%%%%%%%% eqn:3trr %%%%%%%%%%%%%%%%%%%%%%%%%%%%%%%%%%%%%%%
\begin{subequations} \label{eqn:3trr}
\begin{align}
f_0(n) &= q_0(n) \cdot f_1(n) + r_1(n) \cdot f_2(n),  \label{eqn:3trr0} \\ 
f_1(n) &= q_1(n) \cdot f_2(n) + r_2(n) \cdot f_0(n+1), \label{eqn:3trr1} \\
f_2(n) &= q_2(n) \cdot f_0(n+1) + r_0(n+1) \cdot f_1(n+1),  \label{eqn:3trr2}
\end{align}
\end{subequations}
%%%%%%%%%%%%%%%%%%%%%%%%%%%%%%%%%%%%%%%%%%%%%%%%%%%%%%%%%%%%%%%%%%%%%%
for $n \in \Z_{\ge 0}$, where the sequences $f_i(n)$, $q_i(n)$ and $r_i(n)$ are 
defined by  
%%%%%
\begin{alignat*}{3} 
f_0(n) &:= {}_3f_2(\ba + n \bp), \quad & 
q_0(n) &:= u(\ba + n \bp), \quad & 
r_1(n) &:= v(\ba + n \bp), \\
f_1(n) &:= {}_3f_2(\ba + n \bp + \bk), \quad & 
q_1(n) &:= {}^{\sigma}\!u(\ba + n \bp + \bk), \quad & 
r_2(n) &:= {}^{\sigma}\!v(\ba + n \bp + \bk), \\
f_2(n) &:= {}_3f_2(\ba + n \bp + \bl), \quad & 
q_2(n) &:= {}^{\sigma^2}\!u(\ba + n \bp + \bl), \quad & 
r_0(n) &:= {}^{\sigma^2}\!v(\ba + (n-1) \bp + \bl).             
\end{alignat*} 
%%%%%
In view of the modulo $3$ structure in \eqref{eqn:3trr} it is convenient to set      
%%%%%%%%%%%%%%%%%%%%%%% eqn:fqr %%%%%%%%%%%%%%%%%%%%%%%%%%%%%%%%%%%%%%
\begin{subequations} \label{eqn:fqr}
\begin{align}
f(n) &:= f_i((n-i)/3),  \label{eqn:f(n)} \\ 
q(n) &:= q_i((n-i)/3), \qquad \mbox{for} \quad n \equiv i \mod 3, 
\quad i = 0,1,2.  \label{eqn:q(n)} \\ 
r(n) &:= r_i((n-i)/3).  \label{eqn:r(n)}
\end{align}  
\end{subequations} 
%%%%%%%%%%%%%%%%%%%%%%%%%%%%%%%%%%%%%%%%%%%%%%%%%%%%%%%%%%%%%%%%%%%%%%
Then the system \eqref{eqn:3trr} is unified into a single three-term 
recurrence relation  
%%%%%%%%%%%%%%%%%%%%%%%%%%%%% eqn:u3trr %%%%%%%%%%%%%%%%%%%%%%%%%%%%%%%
\begin{equation} \label{eqn:u3trr}
f(n) = q(n) \cdot f(n+1) + r(n+1) \cdot f(n+2), \qquad n \in \Z_{\ge 0}.    
\end{equation}  
%%%%%%%%%%%%%%%%%%%%%%%%%%%%%%%%%%%%%%%%%%%%%%%%%%%%%%%%%%%%%%%%%%%%%%%%
\par
%%%%%%
If $\bk$ is nonnegative, $\bk \in \Z_{\ge 0}^5$, then 
so are $\bl$ and $\bp$ by formula \eqref{eqn:lp}, 
hence all $f(n)$, $n \in \Z_{\ge 0}$, are well defined under single 
assumption \eqref{eqn:WD}. 
If moreover $\bk$ is {\sl balanced}, $s(\bk) = 0$, then so are $\bl$ and 
$\bp$ again by formula \eqref{eqn:lp}, hence all $f(n)$, $n \in \Z_{\ge 0}$ 
have the same Saalsch\"utzian index.    
Thus all these series are convergent under the single 
assumption \eqref{eqn:3saal2}. 
In what follows we refer to $\bk$ as the {\sl seed} vector while $\bp$ as the 
{\sl shift} vector. 
We remark that $\bk$ is primary in the sense that $\bl$ and $\bp$ 
are derived from $\bk$ by the rule \eqref{eqn:lp}, but $\bp$ is likewise  
important because it is $\bp$ rather than $\bk$ that is directly responsible 
for the asymptotic behavior of the sequence $f(n)$.   
%%%%%%%%%%%%%%%%%%%%%%%%%%%%% subsec:simult %%%%%%%%%%%%%%%%%%%%%%%%%%%%%%
\subsection{Simultaneousness} \label{subsec:simult}
%%%%%%%%%%%%%%%%%%%%%%%%%%%%%%%%%%%%%%%%%%%%%%%%%%%%%%%%%%%%%%%%%%%%%%%%
In place of the series ${}_3f_2(\ba)$ we consider another series 
%%%%%%%%%%%%%%%%%%%%%%%%%%%%%%%%% eqn:3g2 %%%%%%%%%%%%%%%%%%%%%%%%%%%%%%%
\begin{equation} \label{eqn:3g2}
{}_3g_2(\ba) = 
{}_3g_2\begin{pmatrix}
a_0, & a_1, & a_2 \\
     & b_1, & b_2 
\end{pmatrix}
:= 
{}_3f_2 \! 
\begin{pmatrix} 
a_0, & a_0-b_1+1, & a_0-b_2+1 \\
     & a_0-a_1+1, & a_0-a_2+1   
\end{pmatrix}. 
\end{equation}
%%%%%%%%%%%%%%%%%%%%%%%%%%%%%%%%%%%%%%%%%%%%%%%%%%%%%%%%%%%%%%%%%%%%%%%% 
Let $\bk$, $\bl$ and $\bp$ be vectors as in \eqref{eqn:lp} such that $s(\bk) = 0$ 
and hence $s(\bl) = s(\bp) = 0$. 
By assertion (3) of \cite[Theorem 1.1]{EI} the contiguous relation 
\eqref{eqn:3tcr} for ${}_3f_2(\ba)$ is {\sl simultaneously} satisfied by 
${}_3h_2(\ba) := \exp( \pi \sqrt{-1} \, s(\ba)) \, {}_3g_2(\ba)$, but    
the factor $\exp( \pi \sqrt{-1} \, s(\ba))$ is irrelevant by 
$s(\bk) = s(\bl) = 0$, thus \eqref{eqn:3tcr} is satisfied by ${}_3g_2(\ba)$ itself.   
Let $g_i(n)$ and $g(n)$ be defined from ${}_3g_2(\ba)$ in the same manner 
as $f_i(n)$ and $f(n)$ are defined from ${}_3f_2(\ba)$ in \S\ref{subsec:rr}, that is, let       
%%%%%%%%%%%%%%%%%%%%%%%%%% eqn:gi(n)-g(n) %%%%%%%%%%%%%%%%%%%%%%%%%%%%%%%%%
\begin{subequations} \label{eqn:gi(n)-g(n)} 
\begin{gather}
g_0(n) := {}_3g_2(\ba + n \bp), \quad 
g_1(n) := {}_3g_2(\ba + n \bp + \bk), \quad 
g_2(n) := {}_3g_2(\ba + n \bp + \bl),  \label{eqn:gi(n)} \\[2mm]
g(n) := g_i((n-i)/3) \qquad \mbox{for} \quad n \equiv i \mod 3, \quad i = 0,1,2. 
\label{eqn:gg(n)} 
\end{gather}   
\end{subequations}
%%%%%%%%%%%%%%%%%%%%%%%%%%%%%%%%%%%%%%%%%%%%%%%%%%%%%%%%%%%%%%%%%%%%%%%
Then the sequences $f(n)$ in \eqref{eqn:f(n)} and $g(n)$ in \eqref{eqn:gg(n)} 
solve the {\sl same} recurrence relation \eqref{eqn:u3trr}. 
With this observation we are now ready to consider continued fractions.    
%%%%%%%%%%%%%%%%%%%%%%%%%%%%%%% sec:cf %%%%%%%%%%%%%%%%%%%%%%%%%%%%%%%%%%
\section{Continued Fractions} \label{sec:cf}
%%%%%%%%%%%%%%%%%%%%%%%%%%%%%%%%%%%%%%%%%%%%%%%%%%%%%%%%%%%%%%%%%%%%%%%
First we present a general principle to establish an exact error estimate for 
the approximants to a continued fraction.   
Next we announce the final goal of this article, Theorems 
\ref{thm:cf-straight} and \ref{thm:cf-cyclic}, which will be achieved by 
the principle after a rather long journey of asymptotic analysis. 
% that extends from \S \ref{sec:cspm} to \S \ref{sec:casorati}.        
%%%%%%%%%%%%%%%%%%%%%%%%%%%%%%% subsec:gee %%%%%%%%%%%%%%%%%%%%%%%%%%%%%%
\subsection{A General Error Estimate} \label{subsec:gee}
%%%%%%%%%%%%%%%%%%%%%%%%%%%%%%%%%%%%%%%%%%%%%%%%%%%%%%%%%%%%%%%%%%%%%%%
Let $\{ q(n) \}_{n=0}^{\infty}$ and  $\{ r(n) \}_{n=1}^{\infty}$ be 
sequences of complex numbers such that $r(n)$ is nonzero  
for every $n \in \N := \Z_{\ge 1}$.  
We consider a sequence of finite continued fractions 
%%%%%%%%%%%%%%%%%%%%%%%%%%%%%%% eqn:cf %%%%%%%%%%%%%%%%%%%%%%%%%%%%%%%%%%
\begin{equation} \label{eqn:cf}
q(0) + \overset{n}{\underset{j=1}{\cfL}} \,\, \frac{r(j)}{q(j)} 
:= q(0) \, + \, \frac{r(1)}{q(1)} \, \mathop{}_{+} \,\,  
\mathop{}_{\cdots\cdots} \,\, \mathop{}_{+} \, \frac{r(n)}{q(n)}, \qquad 
n \in \Z_{\ge 0}.  
\end{equation}
%%%%%%%%%%%%%%%%%%%%%%%%%%%%%%%%%%%%%%%%%%%%%%%%%%%%%%%%%%%%%%%%%%%%%%%%
The convergence of \eqref{eqn:cf} can be described in terms of 
the three-term recurrence relation 
%%%%%%%%%%%%%%%%%%%%%%%%%% eqn:3trr-x %%%%%%%%%%%%%%%%%%%%%%%%%%%%%%%%%%%%
\begin{equation} \label{eqn:3trr-x}
x(n) = q(n) \cdot x(n+1) + r(n+1) \cdot x(n+2), \qquad n \in \Z_{\ge 0}.  
\end{equation}
%%%%%%%%%%%%%%%%%%%%%%%%%%%%%%%%%%%%%%%%%%%%%%%%%%%%%%%%%%%%%%%%%%%%%%%%
A nontrivial solution $X(n)$ to equation \eqref{eqn:3trr-x} is said to be 
{\sl recessive} if $X(n)/Y(n) \to 0$ as $n \to +\infty$ for any 
solution $Y(n)$ not proportional to $X(n)$. 
Recessive solution, if it exists, is unique up to nonzero constant multiples.  
Any non-recessive solution is said to be {\sl dominant}.   
%%%%%%%%%%%%%%%%%%%%%%%%%%%%%%% thm:pincherle %%%%%%%%%%%%%%%%%%%%%%%%%%%%%
\begin{theorem}[Pincherle \cite{Pincherle}] \label{thm:pincherle}
Sequence \eqref{eqn:cf} is convergent if and only if the recurrence 
equation \eqref{eqn:3trr-x} has a recessive solution $X(n)$, in which case 
\eqref{eqn:cf} converges to the ratio $X(0)/X(1)$.  
\end{theorem} 
%%%%%%%%%%%%%%%%%%%%%%%%%%%%%%%%%%%%%%%%%%%%%%%%%%%%%%%%%%%%%%%%%%%%%%%%
\par
%%%%%%%% 
Let us make this theorem more quantitative. 
For any nontrivial solution $x(n)$ to equation \eqref{eqn:3trr-x} and any 
positive integer $m \in \N$ one has  
%%%%%%%%
\[
\frac{x(0)}{x(1)} = q(0) + 
\overset{m-1}{\underset{j=1}{\cfL}} \,\, \frac{r(j)}{q(j)} 
\, \mathop{}_{+} \, \dfrac{r(m)}{q(m) + \frac{r(m+1)}{ \frac{x(m+1)}{x(m+2)}}}. 
\]
%%%%%%%%
Thus if $x(n; m)$ is a nontrivial solution to \eqref{eqn:3trr-x} that vanishes at $n = m+2$, 
then 
%%%%%%%%
\[
\frac{x(0; m)}{x(1; m)} = q(0) + \overset{m}{\underset{j=1}{\cfL}} \,\, \frac{r(j)}{q(j)}, 
\quad \mbox{or equivalently}, \quad 
\frac{x(1; m)}{x(0; m)} = \overset{m}{\underset{j=0}{\cfL}} \,\, \frac{r(j)}{q(j)}, 
\quad r(0) := 1. 
\] 
%%%%%%%%%
One can express the solution $x(n;m)$ in the form   
%%%%%%%%
\[
x(n; m) = X(n) - R(m) \cdot Y(n), \qquad 
R(m) := \frac{X(m+2)}{Y(m+2)}, \qquad m, n \in \Z_{\ge 0},
\] 
%%%%%%%%%
where $X(n)$ and $Y(n)$ are recessive and dominant 
solutions to \eqref{eqn:3trr-x} respectively, so that $R(m) \to 0$ as 
$m \to +\infty$.  
Hence if $X(0)$ is nonzero then so is $x(0; m)$ for every $m \gg 0$ and  
%%%%%%%%%
\[
\frac{X(1)}{X(0)} - \frac{x(1; m)}{x(0; m)} = 
\frac{X(1)}{X(0)} - \frac{X(1) - R(m) \cdot Y(1)}{X(0) - R(m) \cdot Y(0)}
= \frac{\omega(0) \cdot R(m) }{X(0)^2 \, \{1 - R(m) \cdot Y(0)/X(0) \} },  
\]
%%%%%%%%%
where $\omega(n) := X(n) \cdot Y(n+1) - X(n+1) \cdot Y(n)$ is the Casoratian 
of $X(n)$ and $Y(n)$, thus 
%%%%%%%%%%%%%%%%%%%%%%%%%%%%%%% eqn:error %%%%%%%%%%%%%%%%%%%%%%%%%%%%%%%
\begin{equation} \label{eqn:error} 
\frac{X(1)}{X(0)} -  \overset{n}{\underset{j=0}{\cfL}} \,\, \frac{r(j)}{q(j)} 
= \frac{\omega(0) \cdot R(n) }{X(0)^2} 
\left\{1+ O\left (\frac{R(n) \cdot Y(0)}{X(0)} \right) \right\} \quad 
\mbox{as} \quad n \to + \infty. 
\end{equation}
%%%%%%%%%%%%%%%%%%%%%%%%%%%%%%%%%%%%%%%%%%%%%%%%%%%%%%%%%%%%%%%%%%%%%%%
\par
%%%%%% 
In order to apply this general estimate to continued fractions for ${}_3f_2(1)$, 
we want to set up the situation in which the sequences $f(n)$ in 
\eqref{eqn:f(n)}  and $g(n)$ in \eqref{eqn:gg(n)} are recessive and dominant 
solutions, respectively,  to the recurrence relation \eqref{eqn:u3trr}.   
We present in \S\ref{sec:cspm} a sufficient condition for $f(n)$ to be recessive,  
while we impose in \S \ref{sec:ds} a further constraint that insures the 
dominance of $g(n)$.  
In fact, upon assuming those conditions, we deduce asymptotic representations 
for $f(n)$ and $g(n)$ showing that they are actually recessive and dominant 
respectively.  
The asymptotic analysis there is used not only to prove such a qualitative 
assertion but also to get a precise asymptotic behavior for the ratio 
$R(n) = f(n+2)/g(n+2)$.  
We have also to evaluate the initial term $\omega(0)$ for the Casoratian  
of $f(n)$ and $g(n)$; this final task is done in \S \ref{sec:casorati}.   
%%%%%%%%%%%%%%%%%%%%%%%%% subsec:mrcf %%%%%%%%%%%%%%%%%%%%%%%%%%%%%%%%%% 
\subsection{Main Results on Continued Fractions}  \label{subsec:mrcf}
%%%%%%%%%%%%%%%%%%%%%%%%%%%%%%%%%%%%%%%%%%%%%%%%%%%%%%%%%%%%%%%%%%%%%%
Let $\{ q(n)\}_{n=0}^{\infty}$ and $\{r(n)\}_{n=1}^{\infty}$ be sequences 
\eqref{eqn:q(n)} and \eqref{eqn:r(n)} derived from $u(\ba)$ and $v(\ba)$ 
as in formula \eqref{eqn:uv2}. 
Consider the continued fraction $\mathbf{K}_{j=0}^{\infty} \, r(j)/q(j)$, where 
$r(0) := 1$ by convention.  
It is said to be {\sl well defined} if $q(j)$ and $r(j)$ take finite values with 
$r(j)$ nonzero for every $j \ge 0$.   
%%%%%%
\par
%%%%%%
Let $\cS(\R)$ be the set of all real vectors 
$\bp = (p_0, p_1, p_2; q_1, q_2) \in \R^5$ such that 
%%%%%%%%%%%%%%%%%%%%%%%%%%%%%% eqn:bp %%%%%%%%%%%%%%%%%%%%%%%%%%%%%%%%%
\begin{equation} \label{eqn:bp}
s(\bp) = 0; \qquad p_1, \, p_2 \le p_0 < q_1 \le q_2 < p_1 + p_2. 
\end{equation}  
%%%%%%%%%%%%%%%%%%%%%%%%%%%%%%%%%%%%%%%%%%%%%%%%%%%%%%%%%%%%%%%%%%%%%%
Note that \eqref{eqn:bp} in particular implies $p_1, p_2 > 0$ and that $\cS(\R)$ 
is a $4$-dimensional polyhedral convex cone defined by a linear equation and a 
set of linear inequalities. 
It is the space to which the shift vector $\bp$ in \eqref{eqn:lp} should 
belong; or rather as an integer vector it should lie on   
%%%%%%%%%%%%%%%%%%%%%%%%%%%%% eqn:cS(Z) %%%%%%%%%%%%%%%%%%%%%%%%%%%%%%
\begin{equation} \label{eqn:cS(Z)} 
\cS(\Z) := \cS(\R) \cap \Z^5. 
\end{equation}
%%%%%%%%%%%%%%%%%%%%%%%%%%%%%%%%%%%%%%%%%%%%%%%%%%%%%%%%%%%%%%%%%%%%
The following functions of $\bp \in \cS(\R)$ play important 
roles in several places of this article:     
%%%%%%%%%%%%%%%%%%%% eqn:A(p) %%%%% eqn:v(p) %%%%%%%%%%%%%%%%%%%%%%%%%%
\begin{align}  
D(\bp) &:= \dfrac{(-1)^{q_1+q_2} p_0^{p_0} p_1^{p_1} p_2^{p_2}}{\prod_{i=0}^2 
\prod_{j=1}^2 (q_j - p_i)^{q_j - p_i}},   \label{eqn:A(p)}  \\[2mm] 
\vD(\bp) &:= e_1^2 e_2^2 + 18 \, e_1 e_2 e_3  
- 2 \, e_2^3 - 8 \, e_1^3 e_3 - 27 \, e_3^2,     
\label{eqn:v(p)} 
\end{align}   
%%%%%%%%%%%%%%%%%%%%%%%%%%%%%%%%%%%%%%%%%%%%%%%%%%%%%%%%%%%%%%%%%%%%%
where $e_1 := p_0+p_1+p_2 = q_1+q_2$,  
$e_2 := p_0 p_1 + p_1 p_2 + p_2 p_0 + q_1 q_2$ and $e_3 := p_0 p_1 p_2$. 
We remark that $\vD(\bp)$ is the discriminant (up to a positive 
constant multiple) of the cubic equation 
%%%%%%
\[
(x-p_0)(x-p_1)(x-p_2) + x (x-q_1) (x-q_2) = 0,  
\]   
%%%%%%
which plays an important role in \S \ref{subsec:2nd}.  
Moreover, for $\bk = (k_0, k_1, k_2; l_1, l_2) \in \Z^5$ we put  
%%%%%%%%%%%%%%%%%%%%%%%%%%% eqn:gamma(a) %%%%%%%%%%%%%%%%%%%%%%%%%%%%
\begin{equation} \label{eqn:gamma(a)} 
\gamma(\ba; \bk) 
:= \dfrac{\vG(a_0) \vG(a_1) \vG(a_2) \vG^2(s(\ba)) }{ \prod_{i=0}^2 
\prod_{j=1}^2 \vG(b_j-a_i+(l_j-k_i)_+)}.     
\end{equation}
%%%%%%%%%%%%%%%%%%%%%%%%%%%%%%%%%%%%%%%%%%%%%%%%%%%%%%%%%%%%%%%%%%%%
\par
%%%%%% 
We are now able to state the main results of this article; they are stated 
in terms of the seed vector $\bk$, but a large part of their proofs will be 
given in terms of the shift vector $\bp$. 
For continued fractions of straight type in Definition \ref{def:admissible} 
we have the following theorem.     
%%%%%%%%%%%%%%%%%%%%%%%%%%%%% thm:cf-straight %%%%%%%%%%%%%%%%%%%%%%%%%
\begin{theorem}[Straight Case] \label{thm:cf-straight}  
If $\bk = (k_0, k_1, k_2; l_1, l_2) \in \cS(\Z)$ satisfies either 
%%%%%%%%%%%%%%%%%%%%%%%%%%%% eqn:cond-straight %%%%%%%%%%%%%%%%%%%%%%%
\begin{equation} \label{eqn:cond-straight}
(\mathrm{a}) \quad  \vD(\bk) \le 0 \qquad \mbox{or} \qquad 
(\mathrm{b}) \quad 2 l_1^2 - 2(k_1 + k_2) l_1 + k_1 k_2 \ge 0, 
\end{equation}
%%%%%%%%%%%%%%%%%%%%%%%%%%%%%%%%%%%%%%%%%%%%%%%%%%%%%%%%%%%%%%%%%%%%%%
then  $|D(\bk)| > 1$ and there exists an error estimate of continued fraction 
expansion  
%%%%%%%%%%%%%%%%%%%%%%%%%%%%%% eqn:cf-straight %%%%%%%%%%%%%%%%%%%%%%%%%%
\begin{equation} \label{eqn:cf-straight}
\frac{{}_3f_2(\ba+\bk)}{{}_3f_2(\ba)} - 
\overset{n}{\underset{j=0}{\cfL}} \,\, \frac{r(j)}{q(j)} 
= \dfrac{ c_{\rs}(\ba; \bk) }{{}_3f_2(\ba)^2}  
\cdot D(\bk)^{-n} \cdot n^{-s(\sba) + \frac{1}{2} } \cdot 
\left\{ 1 + O (n^{-\frac{1}{2} } ) \right\},  
\end{equation} 
%%%%%%%%%%%%%%%%%%%%%%%%%%%%%%%%%%%%%%%%%%%%%%%%%%%%%%%%%%%%%%%%%%%%%%
as $n \to + \infty$, provided that $\rRe \, s(\ba)$ is positive, ${}_3f_2(\ba)$ is 
nonzero and the continued fraction $\mathbf{K}_{j=0}^{\infty} \, r(j)/q(j)$ is well 
defined,  where $D(\bk)$ is defined in \eqref{eqn:A(p)} with $\bp$ 
replaced by $\bk$, while 
%%%%%%%%%%%%%%%%%%%%%%%%%%%%
\begin{equation*}   
c_{\rs}(\ba; \bk) := \rho(\ba; \bk) \cdot e_{\rs}(\ba; \bk) \cdot \gamma(\ba; \bk), 
\end{equation*} 
%%%%%%%%%%%%%%%%%%%%%%%%%%%
with $\rho(\ba; \bk) \in \Q[\ba]$ being the polynomial in \eqref{eqn:rho}, 
explicitly computable from $\bk$, 
%%%%%%%%%%%%%%%%%%%%%%%%%%%% eqn:es(a;k) %%%%%%%%%%%%%%%%%%%%%%%%%%%%%%%
\begin{equation} \label{eqn:es(a;k)} 
e_{\rs}(\ba; \bk) 
:= (2\pi)^{ \frac{3}{2} } \, 
\dfrac{ \prod_{i=0}^2 \prod_{j=1}^2 
(l_j-k_i)^{2(l_j-k_i)+b_j-a_i- \frac{1}{2} } }{ s_2(\bk)^{2 s(\sba)-1} 
\prod_{i=0}^2 k_i^{2 k_i + a_i - \frac{1}{2} } },  
\end{equation}
%%%%%%%%%%%%%%%%%%%%%%%%%%%%%%%%%%%%%%%%%%%%%%%%%%%%%%%%%%%%%%%%%%%%%%
with $s_2(\bk) := k_0 k_1 + k_1 k_2 + k_2 k_0 - l_1 l_2$ and $\gamma(\ba; \bk)$ 
defined by formula \eqref{eqn:gamma(a)}.  
\end{theorem}
%%%%%%%%%%%%%%%%%%%%%%%%%%%%%%%%%%%%%%%%%%%%%%%%%%%%%%%%%%%%%%%%%%%%%%
\par
%%%%%%% 
A numerical inspection shows that about 43 \% of the vectors in 
$\cS(\Z)$ satisfy condition \eqref{eqn:cond-straight} (see Remark \ref{rem:ds}). 
In the straight case with $\bk \in \cS(\Z)$ formulas \eqref{eqn:uv2} become simpler: 
%%%%%%%%%%%%%%%%%%%%%%%%%%%% eqn:uv3 %%%%%%%%%%%%%%%%%%%%%%%%%%%%%%%%%%
\begin{equation} \label{eqn:uv3}
u(\ba) = \dfrac{(-1)^{l_1+l_2} \, \rho(\ba; 2 \bk)}{\rho(\ba+\bk; \bk) \, 
\prod_{i=0}^2 (a_i; \, k_i)}, 
\qquad 
v(\ba) = -\dfrac{(-1)^{l_1+l_2} \, \rho(\ba; \bk) \cdot 
\langle \ba+\bk; \bk \rangle_+}{\rho(\ba+\bk; \bk) \, 
\prod_{i=0}^2 (a_i; \, k_i)}.  
\end{equation} 
%%%%%%%%%%%%%%%%%%%%%%%%%%%%%%%%%%%%%%%%%%%%%%%%%%%%%%%%%%%%%%%%%%%%%%
\par
%%%%%% 
We turn our attention to continued fractions of twisted type in 
Definition \ref{def:admissible}.  
%%%%%%%%%%%%%%%%%%%%%%%%%%%% thm:cf-cyclic %%%%%%%%%%%%%%%%%%%%%%%%%%%%%%
\begin{theorem}[Twisted Case] \label{thm:cf-cyclic} 
If $\bk = (k_0, k_1, k_2; l_1, l_2) \in \Z^5_{\ge0}$ satisfies the condition   
%%%%%%%%%%%%%%%%%%%%%%%%%%%% eqn:cond-cyclic %%%%%%%%%%%%%%%%%%%%%%%%%%%
\begin{equation} \label{eqn:cond-cyclic}
k_0+k_1+k_2 = l_1 + l_2, \qquad
l_1 \le l_2 \le \tau \, l_1,  \quad 
 \tau := (1+ \sqrt{3})/2 = 1.36602540\cdots, 
\end{equation} 
%%%%%%%%%%%%%%%%%%%%%%%%%%%%%%%%%%%%%%%%%%%%%%%%%%%%%%%%%%%%%%%%%%%%% 
then there exists an error estimate of continued fraction expansion  
%%%%%%%%%%%%%%%%%%%%%%%%%%% eqn:cf-cyclic %%%%%%%%%%%%%%%%%%%%%%%%%%%%%%
\begin{equation} \label{eqn:cf-cyclic}
\frac{{}_3f_2(\ba+\bk)}{{}_3f_2(\ba)} - 
\overset{n}{\underset{j=0}{\cfL}} \,\, \frac{r(j)}{q(j)} 
= \dfrac{c_{\rt}(\ba; \bk)}{{}_3f_2(\ba)^2}  
\cdot E(l_1, l_2)^{-n} \cdot n^{-s(\sba) + \frac{1}{2} } \cdot 
\left\{ 1 + O (n^{-\frac{1}{2} } ) \right\}, 
\end{equation}
%%%%%%%%%%%%%%%%%%%%%%%%%%%%%%%%%%%%%%%%%%%%%%%%%%%%%%%%%%%%%%%%%%%%%
as $n \to +\infty$, provided that $\rRe \, s(\ba)$ is positive, ${}_3f_2(\ba)$ is 
nonzero and the continued fraction $\mathbf{K}_{j=0}^{\infty} \, r(j)/q(j)$ 
is well defined,  where $E(l_1, l_2)$ and $c_{\rt}(\ba; \bk)$ are given by 
%%%%%%%%%%%%%%%%%%%%%%%% eqn:A(l1,l2) %%%%%%%%%%%%%%%%%%%%%%%%%%%%%%%%%% 
\begin{align}
E(l_1, l_2) 
&:= \dfrac{(-l_1 - l_2)^{l_1 + l_2}}{(2 l_1 - l_2)^{2 l_1 - l_2} (2 l_2 - l_1)^{2 l_2 - l_1}}, 
\qquad |E(l_1, l_2)| > 1,  \label{eqn:A(l1,l2)} \\[2mm] 
c_{\rt}(\ba; \bk) &:= \rho(\ba; \bk) \cdot e_{\rt}(\ba; \bk) \cdot 
\gamma(\ba; \bk),  \nonumber
\end{align} 
%%%%%%%%%%%%%%%%%%%%%%%%%%%%%%%%%%%%%%%%%%%%%%%%%%%%%%%%%%%%%%%%%%%%%
with $\rho(\ba; \bk) \in \Q[\ba]$ being the polynomial in 
\eqref{eqn:rho}, explicitly computable from $\bk$,  
%%%%%%%%%%%%%%%%%%%%%%%%%%%%% eqn:ec(a;k) %%%%%%%%%%%%%%%%%%%%%%%%%%%%%%
\begin{equation}  \label{eqn:ec(a;k)} 
e_{\rt}(\ba; \bk) 
:= (2 \pi)^{ \frac{3}{2} } \, 
\dfrac{ (2l_1-l_2)^{2(2 l_1-l_2) + 2 b_1-b_2+s(\sba)- \frac{3}{2} } \cdot 
 (2l_2-l_1)^{2(2 l_2-l_1) + 2 b_2-b_1 + s(\sba)- \frac{3}{2} } }{ 3^{s(\sba) -\frac{1}{2} } 
\cdot 
 (l_1+l_2)^{2(l_1+l_2) +a_0+a_1+a_2- \frac{3}{2} } \cdot 
(l_1^2 -l_1 l_2 + l_2^2)^{2 s(\sba)-1}},  
\end{equation}
%%%%%%%%%%%%%%%%%%%%%%%%%%%%%%%%%%%%%%%%%%%%%%%%%%%%%%%%%%%%%%%%%%%%%%%
and $\gamma(\ba; \bk)$ being defined by formula \eqref{eqn:gamma(a)}. 
%%%%%%%%%%%%%%%%%%%%%%%%%%%%%%%%%%%%%%%%%%%%%%%%%%%%%%%%%%%%%%%%%%%%%%%
\end{theorem}
%%%%%%%%%%%%%%%%%%%%%%%%%%%%%%%%%%%%%%%%%%%%%%%%%%%%%%%%%%%%%%%%%%%%%%%
\par
%%%%%%
The proofs of Theorems \ref{thm:cf-straight} and \ref{thm:cf-cyclic} will be 
completed at the end of \S \ref{sec:casorati}. 
%%%%%%%%%%%%%%%%%%%%%%%%%%%%%%% sec:cspm %%%%%%%%%%%%%%%%%%%%%%%%%%%%%%%
\section{Continuous Laplace Method} \label{sec:cspm}
%%%%%%%%%%%%%%%%%%%%%%%%%%%%%%%%%%%%%%%%%%%%%%%%%%%%%%%%%%%%%%%%%%%%%%%%
We shall find a class of directions $\bp = (p_0, p_1, p_2; q_1, q_2) \in \R^5$ in which 
the sequence 
%%%%%%%%%%%%%%%%%%%%%%%%%%%% eqn:tentative %%%%%%%%%%%%%%%%%%%%%%%%%%%%%%%
\begin{equation} \label{eqn:tentative} 
f(n) = {}_3f_2(\ba+n \bp) = {}_3f_2\begin{pmatrix}
a_0 + p_0 n, & a_1 + p_1 n, & a_2 + p_2 n \\
             & b_1 + q_1 n, & b_2 + q_2 n 
\end{pmatrix}, 
\qquad n \in \Z_{\ge 0}, 
\end{equation} 
%%%%%%%%%%%%%%%%%%%%%%%%%%%%%%%%%%%%%%%%%%%%%%%%%%%%%%%%%%%%%%%%%%%%%%%%
behaves like $n^{\alpha}$ as $n \to +\infty$ for some $\alpha \in \R$, where 
we assume $s(\bp) = 0$ so that the Saalsch\"utian indices for $f(n)$ are 
independent of $n$, always equal to $s(\ba)$. 
We remark that the current $f(n)$ corresponds to the sequence $f_0(n)$ in 
\S \ref{subsec:rr}, not to $f(n)$ in formula \eqref{eqn:f(n)}.     
%%%%%%%
\par
%%%%%%%
In terms of the series ${}_3f_2(\ba)$, Thomae's transformation 
\cite[Corollary 3.3.6]{AAR} reads 
%%%%%%%
\[
{}_3f_2\begin{pmatrix}
a_0, & a_1, & a_2 \\
     & b_1, & b_2 
\end{pmatrix}
= \dfrac{\vG(a_1) \vG(a_2)}{\vG(b_1-a_0) \vG(b_2-a_0)} \, 
{}_3f_2\begin{pmatrix}
s(\ba), & b_1-a_0, & b_2-a_0 \\
     & s(\ba) + a_1, & s(\ba) + a_2 
\end{pmatrix}. 
\]
%%%%%%%
To investigate the asymptotic behavior of $f(n)$, take 
Thomae's transformation of \eqref{eqn:tentative} to have    
%%%%%%%%%%%%%%%%%%%%%%%%%%%% eqn:tentative2 %%%%%%%%%%%%%%%%%%%%%%%%%%% 
\begin{subequations} \label{eqn:tentative2}
\begin{align}
f(n) &= \psi_1(n) \cdot f_1(n), \label{eqn:tentative21} \\[1mm]  
\psi_1(n) &:= \dfrac{\vG( a_1 + p_1 n ) 
\vG( a_2 + p_2 n )}{ \vG( b_1 - a_0 + (q_1 - p_0) n ) 
\vG( b_2 - a_0 + (q_2 - p_0) n )}, 
\label{eqn:tentative22} \\[2mm]
f_1(n) &= 
{}_3f_2\begin{pmatrix}
s(\ba), & b_1 - a_0 + (q_1 - p_0) n, & b_2 - a_0 + (q_2 - p_0) n \\
             & s(\ba) + a_1 + p_1 n, & s(\ba) + a_2 + p_2 n 
\end{pmatrix}, \label{eqn:tentative23}
\end{align} 
\end{subequations}
%%%%%%%%%%%%%%%%%%%%%%%%%%%%%%%%%%%%%%%%%%%%%%%%%%%%%%%%%%%%%%%%%%%%%%%%
and then apply ordinary Laplace's method to the Euler integral 
representation for \eqref{eqn:tentative23}. 
Since this analysis is not limited to ${}_3f_2(1)$, 
we shall deal with more general ${}_{p+1}f_p(1)$ series.         
%%%%%%%%%%%%%%%%%%%%%%%%%%%% subsec:eir %%%%%%%%%%%%%%%%%%%%%%%%%%%%%%%%
\subsection{Euler Integral Representations} \label{subsec:eir}
%%%%%%%%%%%%%%%%%%%%%%%%%%%%%%%%%%%%%%%%%%%%%%%%%%%%%%%%%%%%%%%%%%%%%%%%
The renormalized generalized hypergeometric series ${}_{p+1}f_p(z)$ 
is defined by 
%%%%%%%%%%%%%%%%%%%%%%%%%%%%%%% eqn:hgp-z %%%%%%%%%%%%%%%%%%%%%%%%%%%%%%%
\begin{equation} \label{eqn:hgp-z}
{}_{p+1}f_p\left(
\begin{matrix}
a_0, & a_1, & \dots, & a_p \\
     & b_1, & \dots, & b_p 
\end{matrix}
; z \right) := 
\sum_{k=0}^{\infty} \dfrac{\vG(a_0+k) \vG(a_1+k) \cdots 
\vG(a_p+k)}{\vG(1+k) \vG(b_1+k) \cdots \vG(b_p+k)} \, z^k, 
\end{equation}
%%%%%%%%%%%%%%%%%%%%%%%%%%%%%%%%%%%%%%%%%%%%%%%%%%%%%%%%%%%%%%%%%%%%%%%%
where $\ba = (a_0, \dots, a_p; b_1, \dots, b_p) \in \C^{p+1} 
\times \C^p$ are parameters such that none of $a_0, \dots, a_p$ 
is a negative integer or zero. 
Then \eqref{eqn:hgp-z} is absolutely convergent on the open unit 
disk $|z| < 1$. 
%%%%%%
\par
%%%%%%
It is well known that if the parameters $\ba$ satisfy the condition  
%%%%%%%%%%%%%%%%%%%%%%%% eqn:bba %%%%%%%%%%%%%%%%%%%%%%%%%%%%%%%%%%%%%%%%
\begin{equation} \label{eqn:bba} 
\rRe \, b_i > \rRe \, a_i > 0 \qquad (i = 1, \dots, p),   
\end{equation}
%%%%%%%%%%%%%%%%%%%%%%%%%%%%%%%%%%%%%%%%%%%%%%%%%%%%%%%%%%%%%%%%%%%%%%%%
then the improper integral of Euler type  
%%%%%%%%%%%%%%%%%%%%%%%%
\begin{equation*} 
E_p(\ba ; z) 
:= \int_{I^p} \phi_p( \bt ; \ba ; z) \, d \bt, \qquad    
\phi_p( \bt ; \ba ; z) 
:=\dfrac{ \prod_{i=1}^p 
t_i^{a_i-1} (1-t_i)^{b_i-a_i-1}}{ (1- z \, t_1 \cdots t_p)^{a_0} } 
\end{equation*}
%%%%%%%%%%%%%%%%%%%%%%
is absolutely convergent and the series \eqref{eqn:hgp-z} admits an 
integral representation 
%%%%%%%%%%%%%%%%%%%%%%%% eqn:eif-z %%%%%%%%%%%%%%%%%%%%%%%%%%%%%%%%%%%%%%
\begin{equation} \label{eqn:eif-z}
{}_{p+1}f_p( \ba ; z)
= \dfrac{ \vG( a_0 ) \cdot 
E_p( \ba ; z )}{\prod_{i=1}^p \vG( b_i-a_i )} 
\qquad \mbox{on the open unit disk} \quad |z| < 1.    
\end{equation}
%%%%%%%%%%%%%%%%%%%%%%%%%%%%%%%%%%%%%%%%%%%%%%%%%%%%%%%%%%%%%%%%%%%%%%%%
where $I = (0, \, 1)$ is the open unit interval,  
$\bt = (t_1, \dots, t_p) \in I^p$ and $d \bt = d t_1 \cdots d t_p$. 
%%%%%%
\par
%%%%%%
We are more interested in ${}_{p+1}f_p(1)$, that is, in the series 
\eqref{eqn:hgp-z} at unit argument $z = 1$: 
%%%%%%%%%%%%%%%%%%%%%%%%%%%%%%% eqn:hgp %%%%%%%%%%%%%%%%%%%%%%%%%%%%%%%%
\begin{equation} \label{eqn:hgp}
{}_{p+1}f_p( \ba ) = {}_{p+1}f_p( \ba; 1 ) 
:= \sum_{k=0}^{\infty} \dfrac{\vG(a_0+k) \vG(a_1+k) \cdots 
\vG(a_p+k)}{\vG(1+k) \vG(b_1+k) \cdots \vG(b_p+k)}.  
\end{equation}
%%%%%%%%%%%%%%%%%%%%%%%%%%%%%%%%%%%%%%%%%%%%%%%%%%%%%%%%%%%%%%%%%%%%%%%
It is well known that series \eqref{eqn:hgp} is absolutely convergent 
if and only if    
%%%%%%%%%%%%%%%%%%%%%%%%%%%%%%% eqn:saal %%%%%%%%%%%%%%%%%%%%%%%%%%%%%%%%
\begin{equation} \label{eqn:saal}
\rRe \, s(\ba) > 0, \qquad 
s( \ba ) := b_1 + \cdots + b_p - a_0 - a_1 - \cdots - a_p,    
\end{equation}
%%%%%%%%%%%%%%%%%%%%%%%%%%%%%%%%%%%%%%%%%%%%%%%%%%%%%%%%%%%%%%%%%%%%%%%%
in which case we have ${}_{p+1}f_p( \ba; z ) \to {}_{p+1}f_p( \ba )$ as 
$z \to 1$ within the open unit disk $|z| < 1$. 
%%%%%%%%%%%%%%%%%%%%%%%% lem:eu %%%%%%%%%%%%%%%%%%%%%%%%%%%%%%%%%%%%%%%%%
\begin{lemma} \label{lem:eu} 
If conditions \eqref{eqn:bba} and \eqref{eqn:saal} are satisfied, 
then the integral  
%%%%%%%%%%%%%%%%%%%%%%%% eqn:eu %%%%%%%%%%%%%%%%%%%%%%%%%%%%%%%%%%%%%%
\begin{equation} \label{eqn:eu}
E_p( \ba ) 
:= \int_{I^p} \phi_p( \bt ; \ba ) \, d \bt, \qquad   
\phi_p( \bt ; \ba ) 
:=  \dfrac{ \prod_{i=1}^p t_i^{a_i-1} 
(1-t_i)^{b_i-a_i-1}}{(1- t_1 \cdots t_p)^{a_0}} 
\end{equation}
%%%%%%%%%%%%%%%%%%%%%%%%%%%%%%%%%%%%%%%%%%%%%%%%%%%%%%%%%%%%%%%%%%%%%%%%
is absolutely convergent and the series \eqref{eqn:hgp} admits an 
integral representation 
%%%%%%%%%%%%%%%%%%%%%%% eqn:eif %%%%%%%%%%%%%%%%%%%%%%%%%%%%%%%%%%%%%%
\begin{equation} \label{eqn:eif}
{}_{p+1}f_p( \ba ) = \dfrac{\vG( a_0 ) \cdot 
E_p( \ba ) }{\prod_{i=1}^p \vG( b_i-a_i )}.    
\end{equation}
%%%%%%%%%%%%%%%%%%%%%%%%%%%%%%%%%%%%%%%%%%%%%%%%%%%%%%%%%%%%%%%%%%%%%%%% 
\end{lemma}
%%%%%%%%%%%%%%%%%%%%%%%%%%%%%%%%%%%%%%%%%%%%%%%%%%%%%%%%%%%%%%%%%%%%%%%%
%%%%%%%%%%%%%%%%%%%%%%%% begin proof %%%%%%%%%%%%%%%%%%%%%%%%%%%%%%%%%%%%
{\it Proof}. 
If $r$ denotes the distance of $\bt$ from $\1 := (1, \dots,1)$ then one has  
%%%%%%%%%%%%%%%%%%%%%%%% eqn:polar %%%%%%%%%%%%%%%%%%%%%%%%%%%%%%%%%%%%%%
\begin{equation} \label{eqn:polar}
\phi_p(\bt; \ba) = O(r^{s(\sba)-p}) \qquad \mbox{as} \quad I^p \ni \bt \to \1,   
\end{equation}
%%%%%%%%%%%%%%%%%%%%%%%%%%%%%%%%%%%%%%%%%%%%%%%%%%%%%%%%%%%%%%%%%%%%%%    
The absolute convergence of integral \eqref{eqn:eu} off a neighborhood 
$U$ of $\1$ is due to condition \eqref{eqn:bba}, while that on $U$  
follows from condition \eqref{eqn:saal} and estimate \eqref{eqn:polar}. 
In view of 
%%%%%%%%%%
\[
\lim_{I \ni z \to 1} \phi_p( \bt ; \ba ; z)  = \phi_p( \bt; \ba ),   
\qquad 
|\phi_p( \bt ; \ba ; z)| \le 
\begin{cases}
\phi_p( \bt; \rRe \, \ba; 0) & (\rRe \, a_0 \le 0, \, z \in I), \\ 
\phi_p( \bt; \rRe \, \ba)    & (\rRe \, a_0 > 0, \, z \in I), 
\end{cases}
\]
%%%%%%%%%
formula \eqref{eqn:eif} is derived from formula \eqref{eqn:eif-z} by 
Lebesgue's convergence theorem.  \hfill $\Box$ \par\medskip
%%%%%%%%%%%%%%%%%%%%%%%% end proof %%%%%%%%%%%%%%%%%%%%%%%%%%%%%%%%%%%%%
The series \eqref{eqn:hgp} is symmetric in $a_0, a_1, \dots, a_p$, 
but the integral representation \eqref{eqn:eif} is symmetric only in 
$a_1, \dots, a_p$.          
This fact is efficiently used in the next subsection.    
%%%%%%%%%%%%%%%%%%%%%%%% subsec:aa-ei %%%%%%%%%%%%%%%%%%%%%%%%%%%%%%%%%% 
\subsection{Asymptotic Analysis of Euler Integrals} \label{subsec:aa-ei}
%%%%%%%%%%%%%%%%%%%%%%%%%%%%%%%%%%%%%%%%%%%%%%%%%%%%%%%%%%%%%%%%%%%%%%%%
Observing that the $0$-th numerator parameter of the sequence $f_1(n)$ in 
\eqref{eqn:tentative23} is independent of $n$, we consider a sequence of the form 
%%%%%%%
\[
f_1(n) := {}_{p+1}f_p\begin{pmatrix}
a_0, & a_1 + k_1 n, & \dots, & a_p + k_p n \\[1mm]
     & b_1 + l_1 n, & \dots, & b_p + l_p n 
\end{pmatrix}, \qquad n \in \Z_{\ge 0}.     
\]
%%%%%%% 
The associated Euler integrals have an almost product structure 
which allows a particularly simple treatment in applying 
Laplace's approximation method.  
%%%%%%%%%%%%%%%%%%%%%%%% prop:cspm %%%%%%%%%%%%%%%%%%%%%%%%%%%%%%%%%%%%%%
\begin{proposition} \label{prop:cspm} 
If $\bk = (0, k_1, \dots, k_p; l_1, \dots, l_p) \in \R^{2 p+1}$ is a 
real vector such that 
%%%%%%%%%%%%%%%%%%%%%%%% eqn:liki %%%%%%%%%%%%%%%%%%%%%%%%%%%%%%%%%%%%%%%%
\begin{equation} \label{eqn:liki}
l_i > k_i > 0 = k_0, \qquad i = 1, \dots, p,  
\end{equation}
%%%%%%%%%%%%%%%%%%%%%%%%%%%%%%%%%%%%%%%%%%%%%%%%%%%%%%%%%%%%%%%%%%%%%%%%
then $E_p(\ba + n \, \bk)$ admits an asymptotic representation as 
$n \to +\infty$,  
%%%%%%%%%%%%%%%%%%%%%%%% eqn:cspm %%%%%%%%%%%%%%%%%%%%%%%%%%%%%%%%%%%%%%
\begin{equation} \label{eqn:cspm}
E_p\begin{pmatrix}
a_0, & a_1 + k_1 n, & \dots, & a_p + k_p n \\[1mm]
     & b_1 + l_1 n, & \dots, & b_p + l_p n 
\end{pmatrix} 
= C \cdot \Phim^n \cdot n^{- \frac{p}{2} } \cdot \left\{1+ O(1/n) \right\},  
\end{equation}
%%%%%%%%%%%%%%%%%%%%%%%%%%%%%%%%%%%%%%%%%%%%%%%%%%%%%%%%%%%%%%%%%%%%%%%%
uniform for $\ba =(a_0, \dots, a_p; b_1, \dots, b_p)$ in any compact 
subset of $(\C \setminus \Z_{\le 0}) \times \C^p \times \C^p$, where 
%%%%%%%%%%%%%%%%%%%%%%%% eqn:Phim-C %%%%%%%%%%%%%%%%%%%%%%%%%%%%%%%%%%%%%
\begin{subequations} \label{eqn:Phim-C}
\begin{align} 
\Phim &:= \prod_{i=1}^p 
\frac{k_i^{k_i} (l_i - k_i)^{l_i -k_i}}{l_i^{l_i}},  \label{eqn:Phim-C1} \\
C &:= (2 \pi)^{ \frac{p}{2} } 
\left(1- \frac{k_1 \cdots k_p}{l_1 \cdots l_p} \right)^{-a_0} 
\prod_{i=1}^p 
\frac{k_i^{a_i- \frac{1}{2} } 
(l_i - k_i)^{b_i-a_i- \frac{1}{2} }}{l_i^{b_i - \frac{1}{2} }}.  \label{eqn:Phim-C2}
\end{align}
\end{subequations}
%%%%%%%%%%%%%%%%%%%%%%%%%%%%%%%%%%%%%%%%%%%%%%%%%%%%%%%%%%%%%%%%%%%%%%%%
\end{proposition}
%%%%%%%%%%%%%%%%%%%%%%%%%%%%%%%%%%%%%%%%%%%%%%%%%%%%%%%%%%%%%%%%%%%%%%%%
%%%%%%%%%%%%%%%%%%%%%%%%%%% begin proof %%%%%%%%%%%%%%%%%%%%%%%%%%%%%%%%
{\it Proof}. 
The proof is an application of the standard Laplace method to the 
integral \eqref{eqn:eif}, so only an outline of it is presented below.  
Replacing $\ba$ with $\ba + n \, \bk$ in 
definition \eqref{eqn:eu}, we have  
%%%%%%%
\begin{equation*}
E_p(\ba + n \, \bk) 
= \int_{I^p} \Phi(\bt)^n \cdot u(\bt) \, d \bt  
= \int_{I^p} e^{- n \, \phi(\sbt)} \cdot u(\bt) \, d \bt,  
\end{equation*}
%%%%%%%
where $\Phi(\bt)$, $\phi(\bt)$ and $u(\bt)$ are defined by 
%%%%%%%
\[
\Phi(\bt) := \prod_{i=1}^p t_i^{k_{i}}(1-t_i)^{l_i-k_i}, \qquad 
\phi(\bt) := - \log \Phi(\bt), \qquad u(\bt) := \phi_p(\bt; \ba). 
\]
%%%%%%%
Observe that $\phi(\bt)$ attains a unique minimum at 
$\bt_0 := (k_1/l_1, \dots, k_p/l_p)$ in the interval $I^p$, 
since  
%%%%%%%
\[
\frac{\partial \phi}{\partial t_i} = 
- \frac{k_i}{t_i} + \frac{l_i-k_i}{1-t_i} = \frac{l_i t_i - k_i}{t_i(1-t_i)}, 
\qquad 
\frac{\partial^2 \phi}{\partial t_i^2} 
= \frac{k_i}{t_i^2} + \frac{l_i - k_i}{(1-t_i)^2} > 0, 
\qquad 
\frac{\partial^2 \phi}{\partial t_i \partial t_j} = 0 \quad (i \neq j).  
\]
%%%%%%%
The standard formula for Laplace's approximation then leads to 
%%%%%%%
\[
\begin{split}
\int_{I^p} e^{- n \, \phi(\sbt)} \cdot u(\bt) \, d \bt 
&= \dfrac{u(\bt_0)}{\sqrt{\mathrm{Hess}(\phi; \bt_0)}} 
\left(\frac{2 \pi}{n} \right)^{ \frac{p}{2} } \exp(-n \, \phi(\bt_0)) \,  
\left\{ 1 + O (1/n) \right\} \\[2mm] 
&= C \cdot \Phim^n \cdot n^{- \frac{p}{2} } \, \left\{1+ O(1/n) \right\}   
\qquad \mbox{as} \quad n \to \infty, 
\end{split} 
\]
%%%%%%%
where $\mathrm{Hess}(\phi; \bt_0)$ is the Hessian of $\phi$ at $\bt_0$ 
while $\Phim$ and $C$ are given by formulas 
\eqref{eqn:Phim-C}.  \hfill $\Box$ 
%%%%%%%%%%%%%%%%%%%%%%%%%%% end proof %%%%%%%%%%%%%%%%%%%%%%%%%%%%%%%%
%%%%%%%%%%%%%%%%%%%%%%%%%%% subsec:rs %%%%%%%%%%%%%%%%%%%%%%%%%%%%%%%%
\subsection{Recessive Sequences} \label{subsec:rs} 
%%%%%%%%%%%%%%%%%%%%%%%%%%%%%%%%%%%%%%%%%%%%%%%%%%%%%%%%%%%%%%%%%%%%%
We return to the special case of ${}_3f_2(1)$ series and prove the following.   
%%%%%%%%%%%%%%%%%%%%%%%%%% thm:rs %%%%%%%%%%%%%%%%%%%%%%%%%%%%%%%%%%%%%
\begin{theorem} \label{thm:rs} 
If $\bp = (p_0,p_1,p_2;q_1,q_2) \in \R^5$ is balanced, $s(\bp) = 0$, and 
%%%%%%%%%%%%%%%%%%%%%%%%%%%%%%%% eqn:rs-a %%%%%%%%%%%%%%%%%%%%%%%%%%%%%
\begin{equation} \label{eqn:rs-a}
p_1 > q_1 - p_0 > 0, \qquad  p_2 > q_2 - p_0 > 0, 
\end{equation}
%%%%%%%%%%%%%%%%%%%%%%%%%%%%%%%%%%%%%%%%%%%%%%%%%%%%%%%%%%%%%%%%%%%%%
then the sequence $f(n) = {}_3f_2(\ba + n \bp)$ in \eqref{eqn:tentative} admits 
an asymptotic representation  
%%%%%%%%%%%%%%%%%%%%%%%%% eqn:f(n)-asym %%%%%%%%%%%%%%%%%%%%%%%%%%%%%%%
\begin{equation} \label{eqn:f(n)-asym}
{}_3f_2(\ba + n \bp) = 
\vG(s(\ba)) \cdot s_2(\bp)^{-s(\sba)} \cdot n^{-2 s(\sba)} \cdot \{ 1+O(1/n) \},  
\quad \mbox{as} \quad n \to + \infty, 
\end{equation}
%%%%%%%%%%%%%%%%%%%%%%%%%%%%%%%%%%%%%%%%%%%%%%%%%%%%%%%%%%%%%%%%%%%%
uniform in any compact subset of $\rRe \, s(\ba) > 0$, 
where $s_2(\bp) := p_0 p_1 + p_1p_2 + p_2 p_0 -q_1 q_2$.  
\end{theorem} 
%%%%%%%%%%%%%%%%%%%%%%%%%%%%%%%%%%%%%%%%%%%%%%%%%%%%%%%%%%%%%%%%%%%%
%%%%%%%%%%%%%%%%%%%%%%%%% begin proof %%%%%%%%%%%%%%%%%%%%%%%%%%%%%%%%%
{\it Proof}. 
By formulas \eqref{eqn:tentative2} and \eqref{eqn:eif} the sequence 
\eqref{eqn:tentative} can be written $f(n) = \psi_2(n) \, e_2(n)$ with  
%%%%%%%%%%%%%%%%%%%%%%%%%%% eqn:psi2e2 %%%%%%%%%%%%%%%%%%%%%%%%%%%%%%%%
\begin{subequations} \label{eqn:psi2e2}
\begin{align}
\psi_2(n) &:=  
\dfrac{\vG( s(\ba) ) \, \vG( a_1 + p_1 n ) \,  
\vG( a_2 + p_2 n )}{ \prod_{j=1}^2 \prod_{i=0, j} \vG( b_j-a_i+(q_j-p_i) n ) }, 
\label{eqn:psi2}
\\[2mm]
e_2(n) &:=  
E_2\begin{pmatrix}
s(\ba), & b_1 - a_0 + (q_1 - p_0) n, & b_2 - a_0 + (q_2 - p_0) n \\[1mm] 
   & s(\ba) + a_1 + p_1 n,              & s(\ba) + a_2 + p_2 n 
\end{pmatrix}.  \label{eqn:e2}
\end{align} 
\end{subequations}
%%%%%%%%%%%%%%%%%%%%%%%%%%%%%%%%%%%%%%%%%%%%%%%%%%%%%%%%%%%%%%%%%%%
\par
%%%%%
Conditions $s(\bp) = 0$ and \eqref{eqn:rs-a} imply that $p_1$, $p_2 > 0$ and 
$q_j - p_i > 0$ for every $j = 1$, $2$ and $i = 0$, $j$, so Stirling's formula 
applied to \eqref{eqn:psi2} yields an asymptotic representation  
%%%%%%%%%%%%%%%%%%%%%%%%% eqn:psi2-asym %%%%%%%%%%%%%%%%%%%%%%%%%%%%%
\begin{equation}  \label{eqn:psi2-asym}
\psi_2(n) = B \cdot A^n \cdot n^{1-2 s(\sba)} \, \{ 1 + O(1/n) \},  
\end{equation}
%%%%%%%%%%%%%%%%%%%%%%%%%%%%%%%%%%%%%%%%%%%%%%%%%%%%%%%%%%%%%%%%%% 
as $n \to + \infty$, where $A$ and $B$ are given by 
%%%%%%
\[ 
A := \frac{p_1^{p_1} \, p_2^{p_2}}{ \prod_{j=1}^2 \prod_{i=0,j} (q_j-p_i)^{q_j-p_i} }, 
\qquad 
B := \frac{\vG(s(\ba)) \cdot p_1^{a_1- \frac{1}{2} } \, 
p_2^{a_2- \frac{1}{2} }}{ 2 \pi \, \prod_{j=1}^2 \prod_{i=0, j}    
(q_j-p_i)^{b_j-a_i- \frac{1}{2} } }. 
\]
%%%%%%%%%%%%%%
\par
%%%%%%%%%%%%%%
When $p = 2$, $k_1 = q_1-p_0$, $k_2 = q_2-p_0$, $l_1 = p_1$, $l_2 = p_2$, 
condition \eqref{eqn:liki} becomes \eqref{eqn:rs-a}, so 
Proposition \ref{prop:cspm} applies to the sequence \eqref{eqn:e2}. 
In this situation we have $\Phim = A^{-1}$ in \eqref{eqn:Phim-C1} 
and $C = B^{-1} \cdot \vG(s(\ba)) \cdot 
\{p_1p_2-(q_1-p_0)(q_2-p_0)\}^{-s(\sba)}$ in \eqref{eqn:Phim-C2}, 
where we have $p_1 p_2 -(q_1-p_0)(q_2-p_0) = s_2(\bp)$ 
from $s(\bp) = 0$.   
Thus formula \eqref{eqn:cspm} reads  
%%%%%%%%%%%%%%%%%%%%%%%%% eqn:e2-asym %%%%%%%%%%%%%%%%%%%%%%%%%%%%
\begin{equation} \label{eqn:e2-asym}
e_2(n) = B^{-1} \cdot \vG(s(\ba)) \cdot s_2(\bp)^{-s(\sba)} 
\cdot A^{-n} \cdot n^{-1} \, \{1 + O(1/n) \} \quad \mbox{as} \quad 
n \to + \infty.    
\end{equation}
%%%%%%%%%%%%%%%%%%%%%%%%%%%%%%%%%%%%%%%%%%%%%%%%%%%%%%%%%%%%%%%%% 
Combining formulas \eqref{eqn:psi2-asym} and \eqref{eqn:e2-asym} we have 
the asymptotic representation \eqref{eqn:f(n)-asym}.  \hfill $\Box$   
%%%%%%%%%%%%%%%%%%%%%%%% sec:dspm %%%%%%%%%%%%%%%%%%%%%%%%%%%%%%%%
\section{Discrete Laplace Method} \label{sec:dspm}
%%%%%%%%%%%%%%%%%%%%%%%%%%%%%%%%%%%%%%%%%%%%%%%%%%%%%%%%%%%%%%%%
When a solution to a recurrence equation is given in terms of 
hypergeometric series, we want to know its asymptotic behavior and 
thereby to check whether it is actually a dominant solution.  
To this end, regarding the series as a ``discrete" integral, 
we develop a {\sl discrete Laplace method} as an analogue to 
the usual (continuous) Laplace method for ordinary integrals.  
While Theorems \ref{thm:cf-straight} and \ref{thm:cf-cyclic} on continued 
fractions are the final goal of this article, the main result of this section, 
Theorem \ref{thm:dspm}, and the method leading to it are the 
methodological core of the article.       
%%%%%%%%%%%%%%%%%%%%%%%%%%%%%%%%%%%%%%%%%%%%%%%%%%%%%%%%%%%%%%%%%%%%%%%%
\subsection{Formulation} \label{subsec:miar}
%%%%%%%%%%%%%%%%%%%%%%%%%%%%%%%%%%%%%%%%%%%%%%%%%%%%%%%%%%%%%%%%%%%%%%%%
Let $\bsigma = (\sigma_i) \in \R^I$, $\blambda = (\lambda_i) \in \R^I$, 
$\btau = (\tau_j) \in \R^J$, $\bmu = (\mu_j) \in \R^J$ be {\sl real} numbers 
indexed by finite sets $I$ and $J$. 
Suppose that the pairs $(\bsigma, \btau)$ and $(\blambda, \bmu)$ are {\sl balanced} 
to the effect that   
%%%%%%%%%%%%%%%%%%%%%%%%%%%%%%%%% eqn:ap-bq2 %%%%%%%%%%%%%%%%%%%%%%%%%%%%
\begin{equation} \label{eqn:ap-bq2}
\sum_{i \in I} \sigma_i = \sum_{j \in J} \tau_j, \qquad 
\sum_{i \in I} \lambda_i = \sum_{j \in J} \mu_j.    
\end{equation}
%%%%%%%%%%%%%%%%%%%%%%%%%%%%%%%%%%%%%%%%%%%%%%%%%%%%%%%%%%%%%%%%%%%%%%%%  
Let $\bal(n) = (\alpha_i(n)) \in \C^I$ and $\bbeta(n) = (\beta_j(n))  \in \C^J$ be 
sequences in $n \in \N$ of complex numbers indexed by $i \in I$ and $j \in J$. 
Suppose that they are {\sl bounded}, that is, for some constant $R > 0$,   
%%%%%%%%%%%%%%%%%%%%%%%%%%%%%%%%% eqn:bdd %%%%%%%%%%%%%%%%%%%%%%%%%%%%%%%
\begin{equation} \label{eqn:bdd}
|\alpha_i(n) | \le R \quad (i \in I); \qquad |\beta_j(n) | \le R \quad (j \in J), 
\qquad {}^{\forall}n \in \N.  
\end{equation} 
%%%%%%%%%%%%%%%%%%%%%%%%%%%%%%%%%%%%%%%%%%%%%%%%%%%%%%%%%%%%%%%%%%%%%%%%
In practical applications $\bal(n)$ and $\bbeta(n)$ will typically be independent 
of $n$, however allowing such a moderate dependence upon $n$ as in 
\eqref{eqn:bdd} is quite helpful in developing the theory.  
%%%%%%
\par
%%%%%%
Given $0 \le r_0 < r_1 \le + \infty$, we consider the sum of gamma products  
%%%%%%%%%%%%%%%%%%%%%%%%%%% eqn:g(n) %%%%%%%%%%%%%%%%%%%%%%%%%%%%%%%%%%%%%
\begin{equation} \label{eqn:g(n)}
g(n) := \sum_{k = \lceil r_0 n \rceil}^{\lceil r_1 n \rceil-1} G(k; n), \quad  
G(k; n) := \dfrac{\prod_{i \in I}
 \vG( \sigma_i \, k + \lambda_i \, n + \alpha_i(n))}{\prod_{j \in J} 
\vG( \tau_j \, k + \mu_j \, n + \beta_j(n))}, \quad n \in \N,       
\end{equation}
%%%%%%%%%%%%%%%%%%%%%%%%%%%%%%%%%%%%%%%%%%%%%%%%%%%%%%%%%%%%%%%%%%%%%%%%  
where $\lceil x \rceil := \min \{ m \in \Z \,:\, x \le m \}$ denotes the ceiling function.   
We remark that the {\sl reflection} of discrete variable 
$k \mapsto \lceil r_0 n \rceil + \lceil r_1 n \rceil -1 -k$ in \eqref{eqn:g(n)} induces an involution 
%%%%%%%%%%%%%%%%%%%%%%%%%%%%%%%%% eqn:refl %%%%%%%%%%%%%%%%%%%%%%%%%%%%%%%%
\begin{subequations} \label{eqn:refl}
\begin{alignat}{3}
\sigma_i' &= - \sigma_i, \qquad  & \lambda_i' &= \lambda_i + \sigma_i (r_0 + r_1),  \qquad & 
\alpha_i'(n) &= \alpha_i(n) - \sigma_i \, r(n),  \label{eqn:refl-s} \\
   \tau_j' &= - \tau_j,     \qquad  & \mu_j'      &= \mu_j + \tau_j (r_0 + r_1), \qquad & 
\beta_j'(n) &= \beta_j(n)  - \tau_j \, r(n),   \label{eqn:refl-t}
\end{alignat}
\end{subequations}
%%%%%%%%%%%%%%%%%%%%%%%%%%%%%%%%%%%%%%%%%%%%%%%%%%%%%%%%%%%%%%%%%%%%%%%%
where $r(n) := (r_0+r_1) n +1 - \lceil r_0 n \rceil - \lceil r_1 n \rceil $ and the resulting     
data are indicated with a prime, while the reflection leaves $r_0$ and $r_1$ unchanged.  
Since $-1 < r(n) \le 1$, if $\alpha_i(n)$ and $\beta_j(n)$ are bounded then so are 
$\alpha_i'(n)$ and $\beta_j'(n)$.  
This reflectional symmetry is helpful in some occasions. 
Moreover, for any integer $s \le r_0$ the {\sl shift} $k \mapsto k + s n$ in 
\eqref{eqn:g(n)} results in the translations   
%%%%%%%%%%%%%%%%%%%%%%%%%%%%%%%% eqn:shift %%%%%%%%%%%%%%%%%%%%%%%%%%%%%%%%%%%%%
\begin{equation} \label{eqn:shift}
r_0 \mapsto r_0 - s, \quad r_1 \mapsto r_1 - s; \quad 
\sigma_i \mapsto \sigma_i, \quad \lambda_i \mapsto  \lambda_i + \sigma_i s; \quad 
\tau_j \mapsto \tau_j, \quad \mu_j \mapsto \mu_j + \tau_j s.     
\end{equation}    
%%%%%%%%%%%%%%%%%%%%%%%%%%%%%%%%%%%%%%%%%%%%%%%%%%%%%%%%%%%%%%%%%%%%%%%%
Taking $s = \lfloor r_0 \rfloor$ we may assume $0 \le r_0 < 1$, where 
$\lfloor x \rfloor := \max \{ m \in \Z \,:\, m \le x \}$ is the floor function.    
This normalization is also sometimes convenient.      
%%%%%%
\par
%%%%%% 
It is insightful to rewrite the gamma product $G(k; n)$ as    
%%%%%%%%%%%%%%%%%%%%%%%% eqn:H %%%%%%%%%%%%%%%%%%%%%%%%%%%%%%%%%%%%%%%%%%
\begin{equation} \label{eqn:H}
G(k; n) = H\left( k/n ; n \right), 
\qquad 
H(x; n) := \dfrac{\prod_{i \in I} \vG( l_i(x) \, n + \alpha_i(n))}{\prod_{j \in J} 
\vG( m_j(x) \, n + \beta_j(n))},   
\end{equation}
%%%%%%%%%%%%%%%%%%%%%%%%%%%%%%%%%%%%%%%%%%%%%%%%%%%%%%%%%%%%%%%%%%%%%%%% 
where $l_i(x)$ and $m_j(x)$ are affine functions defined by  
%%%%%%%%%%%%%%%%%%%%%%%%%%%%%%%%%
\begin{equation*} 
l_i(x) := \sigma_i x + \lambda_i \quad (i \in I), \qquad 
m_j(x) := \tau_j x + \mu_j \qquad (j \in J). 
\end{equation*}
%%%%%%%%%%%%%%%%%%%%%%%%%%%%%%%%
We remark that condition \eqref{eqn:ap-bq2} is equivalent to the balancedness of 
affine functions  
%%%%%%%%%%%%%%%%%%%%%%%%%%%%%%%%% eqn:balanced %%%%%%%%%%%%%%%%%%%%%%%%%%%
\begin{equation} \label{eqn:balanced}
\sum_{i \in I} l_i(x) = \sum_{j \in J} m_j(x),  
\qquad {}^{\forall} x \in \mathbb{R}. 
\end{equation}
%%%%%%%%%%%%%%%%%%%%%%%%%%%%%%%%%%%%%%%%%%%%%%%%%%%%%%%%%%%%%%%%%%%%%%%%
\par
%%%%%%%    
The sum $g(n)$ is said to be {\sl admissible} if   
%%%%%%%%%%%%%%%%%%%%%%%%%%%%%%%%% eqn:ap-bq1 %%%%%%%%%%%%%%%%%%%%%%%%%%%%
\begin{subequations} \label{eqn:ap-bq1}
\begin{alignat}{4} 
\sigma_i &\neq 0; \qquad & l_i( r_0 ) &\ge 0, \qquad & 
l_i( r_1 ) &\ge 0 \qquad & (i &\in I), 
\label{eqn:ap-bq11} \\ 
\tau_j &\neq 0; \qquad & m_j( r_0 ) &\ge 0, \qquad & 
m_j( r_1 ) &\ge 0 \qquad & (j &\in J),    
\label{eqn:ap-bq12} 
\end{alignat}
\end{subequations}
%%%%%%%%%%%%%%%%%%%%%%%%%%%%%%%%%%%%%%%%%%%%%%%%%%%%%%%%%%%%%%%%%%%%%%%%
where if $r_1 = +\infty$ then by $l_i( r_1 ) \ge 0$ and $m_j( r_1 ) \ge 0$ we mean 
$\sigma_i > 0$ and $\tau_j > 0$. 
Condition \eqref{eqn:ap-bq1} says that $l_i(x)$ and $m_j(x)$ are non-constant affine functions 
taking nonnegative values at both ends of the interval $[r_0, \, r_1]$, so they must be 
positive in its interior, that is,   
%%%%%%%%%%%%%%%%%%%%%%%%%%%%%%%%% eqn:positive %%%%%%%%%%%%%%%%%%%%%%%%%%%%
\begin{equation} \label{eqn:positive}
l_i(x) > 0 \quad (i \in I); \quad m_j(x) > 0 \quad (j \in J), 
\qquad r_0 < {}^{\forall} x < r_1.  
\end{equation}  
%%%%%%%%%%%%%%%%%%%%%%%%%%%%%%%%%%%%%%%%%%%%%%%%%%%%%%%%%%%%%%%%%%%%%%%%
\par
%%%%%%
To work near the endpoints of the interval we introduce four index subsets       
%%%%%%%%%%%%%%%%%%%%%%%%%%%% eqn:IJ01 %%%%%%%%%%%%%%%%%%%%%%%%%%%%%%%%%%%
\begin{subequations} \label{eqn:IJ01}
\begin{alignat}{2}
I_0 &:= \{ i \in I \,:\, l_i( r_0 ) = 0 \}, \qquad & I_1 &:= \{ i \in I \,:\, l_i( r_1 ) = 0 \},  \label{eqn:I01} \\
J_0 &:= \{ j \in J \,:\, m_j( r_0 ) = 0 \}, \qquad & J_1 &:= \{ j \in J \,:\, m_j( r_1 ) = 0 \}.  \label{eqn:J01} 
\end{alignat}
\end{subequations}
%%%%%%%%%%%%%%%%%%%%%%%%%%%%%%%%%%%%%%%%%%%%%%%%%%%%%%%%%%%%%%%%%%%%%%%%
Then there exists a positive constant $c > 0$ such that 
%%%%%%%%%%%%%%%%%%%%%%%%%%%%%%%%% eqn:positive2 %%%%%%%%%%%%%%%%%%%%%%%%%%%
\begin{equation} \label{eqn:positive2}
l_i(x) \ge c \quad (i \in I \setminus (I_0 \cup I_1)); \quad 
m_j(x) \ge c \quad (j \in J \setminus (J_0 \cup J_1)), 
\quad r_0 \le {}^{\forall} x  \le r_1.  
\end{equation}  
%%%%%%%%%%%%%%%%%%%%%%%%%%%%%%%%%%%%%%%%%%%%%%%%%%%%%%%%%%%%%%%%%%%%%%%%
This ``uniformly away from zero'' property will be important in applying a version 
of Stirling's formula which is given later in \eqref{eqn:stirling}, especially when 
$I_0 \cup I_1 \cup J_0 \cup J_1 = \emptyset$ (regular case).      
%%%%%%%%%%%%%%%%%%%%%%%%%%% lem:I0I1 %%%%%%%%%%%%%%%%%%%%%%%%%%%%%%%%%%%%%
\begin{lemma} \label{lem:I0I1}  We have $\sigma_i > 0$ for $i \in I_0$ while 
$\sigma_i < 0$ for $i \in I_1$, in particular  $I_0 \cap I_1 = \emptyset$. 
Similarly we have $\tau_j > 0$ for $j \in J_0$ while $\tau_j < 0$ for $j \in J_1$, 
in particular  $J_0 \cap J_1 = \emptyset$. 
If  
%%%%%%%%%%%%%%%%%%%%%%%%%%%%%%% eqn:punc0 %%%%%%%%%%%%%%%%%%%%%%%%%%%%%%%
\begin{equation} \label{eqn:punc0}
\alpha^{(\nu)}_i(n) := \alpha_i(n) + \sigma_i (\lceil r_{\nu} n \rceil - r_{\nu} n ) 
\not \in \Z_{\le 0}  + |\sigma_i| \, \Z_{\le -\nu}, 
\quad  i \in I_{\nu}, \,\,  \nu = 0, 1, \,\, n \in \N,     
\end{equation}
%%%%%%%%%%%%%%%%%%%%%%%%%%%%%%%%%%%%%%%%%%%%%%%%%%%%%%%%%%%%%%%%%%%%%%%%
then the sum $g(n)$ is well defined, that is, every summand 
$G(k; n) = H(k/n; n)$ in \eqref{eqn:g(n)} takes a finite value for any 
$n \ge (R+1)/c$ with $R$ and $c$ given in \eqref{eqn:bdd} and 
\eqref{eqn:positive2}.      
%%%%%%%%%%%%%%%%%%%%%%%%%%%%%%%%%%%%%%%%%%%%%%%%%%%%%%%%%%%%%%%%%%%%%%%%
\end{lemma} 
%%%%%%%%%%%%%%%%%%%%%%%%%% begin proof %%%%%%%%%%%%%%%%%%%%%%%%%%%%%%%%%%%
{\it Proof}. 
By condition \eqref{eqn:ap-bq11}, if $i \in I_0$ then $0 \le l_i( r_1 ) = l_i( r_1 ) - 
l_i( r_0 ) = (r_1-r_0) \sigma_i$ with $r_1-r_0 > 0$ and $\sigma_i \neq 0$, 
which forces $\sigma_i > 0$, while if $i \in I_1$ then  $0 \le l_i( r_0 ) = l_i( r_0 ) - 
l_i( r_1 ) = (r_0-r_1) \sigma_i$ with $r_0-r_1 < 0$ and $\sigma_i \neq 0$, 
which forces $\sigma_i < 0$.  
A similar argument using \eqref{eqn:ap-bq12} leads to the assertions for 
$J_0$ and $J_1$.   
The sum $g(n)$ fails to make sense only when the argument of an upper gamma 
factor of a summand $G(k; n)$ takes a negative integer value or zero, that is,  
%%%%%
\[
\sigma_i k + \lambda_i n + \alpha_i(n)  = l_i(k/n) \, n + \alpha_i(n)  \in \Z_{\le 0}, 
\qquad {}^{\exists} i \in I, \quad 
\lceil r_0 n \rceil \le {}^{\exists} k \le \lceil r_1 n \rceil -1.  
\]
%%%%% 
This cannot occur for $i \in I \setminus (I_0 \cup I_1)$ and $n \ge (R+1)/c$,  
since \eqref{eqn:bdd} and \eqref{eqn:positive2} imply that 
$l_i(k/n) \, n + \rRe \, \alpha_i(n) \ge c n -R \ge 1$ for any  
$k \in \Z$ such that $r_0 \le k/n \le r_1$. 
Observe that 
%%%%%
\[
\sigma_i k + \lambda_i n + \alpha_i(n) = \sigma_i \, l + l_i(r_0) n + 
\alpha^{(0)}_i(n) = \sigma_i \, l + \alpha^{(0)}_i(n), \qquad i \in I_0,  
\]
%%%%%
where $l := k - \lceil r_0 n \rceil$ ranges over 
$0, 1, \dots, \lceil r_1 n \rceil - \lceil r_0 n \rceil -1$. 
This cannot be a negative integer or zero, if condition \eqref{eqn:punc0} is 
satisfied for $\nu = 0$.      
A similar argument can be made for $\nu = 1$, since condition \eqref{eqn:punc0} 
for $\nu =1$ is obtained from that for $\nu=0$ by applying reflectional 
symmetry \eqref{eqn:refl}.   
Thus if \eqref{eqn:punc0} is satisfied then $g(n)$ is well defined for 
$n \ge (R+1)/c$. 
\hfill $\Box$ \par\medskip
%%%%%%%%%%%%%%%%%%%%%%%%%%%%% end proof %%%%%%%%%%%%%%%%%%%%%%%%%%%%%%%%%%
To carry out analysis it is convenient to quantify condition \eqref{eqn:punc0} 
by writing   
%%%%%%%%%%%%%%%%%%%%%%%%%%%% eqn:punc %%%%%%%%%%%%%%%%%%%%%%%%%%%%%%%%%%%
\begin{equation} \label{eqn:punc}
\delta_{\nu}(n) := \min \Big\{ \, 1, \, 
\prod_{i \in I_{\nu} } \mathrm{dist} (\alpha_i^{(\nu)} (n), \, 
\Z_{\le 0} + |\sigma_i| \Z_{\le - \nu}) \, \Big\} > 0, \quad \nu = 0, 1, \,\, 
n \in \N, 
\end{equation}
%%%%%%%%%%%%%%%%%%%%%%%%%%%%%%%%%%%%%%%%%%%%%%%%%%%%%%%%%%%%%%%%%%%%%%%
where $\mathrm{dist}(z, Z)$ stands for the distance between a point $z$ and 
a set $Z$ in $\C$, and cut off by $1$ is simply to make 
$\delta_{\nu}(n) \le 1$ as it really works only when $0 < \delta_{\nu}(n) \ll 1$.   
Condition \eqref{eqn:punc} or \eqref{eqn:punc0} is referred to as the 
{\sl genericness} for the data $\bal(n)$. 
%%%%%%%%%%%%%%%%%%%%%%%%%% subsec:results %%%%%%%%%%%%%%%%%%%%%%%%%%%%%%%%
\subsection{Main Results on Discrete Laplace Method} \label{subsec:results}
%%%%%%%%%%%%%%%%%%%%%%%%%%%%%%%%%%%%%%%%%%%%%%%%%%%%%%%%%%%%%%%%%%%%%%%
To state the main result of this section we introduce the following quantities: 
%%%%%%%%%%%%%%%%%%%%%%%%%% eqn:Phi %%%% eqn:u %%%%%%%%%%%%%%%%%%%%%%%%%%%
\begin{subequations} \label{eqn:G-asymp2}
\begin{align}
\Phi(x) &:= \prod_{i \in I} l_i(x)^{l_i(x)} \prod_{j \in J} m_j(x)^{-m_j(x)},  
\label{eqn:Phi} \\
u(x; n) &:= (2 \pi)^{\frac{|I|-|J|}{2}} 
\prod_{i \in I} l_i(x)^{\alpha_i(n) - \frac{1}{2} } 
\prod_{j \in J} m_j(x)^{ \frac{1}{2} - \beta_j(n)},  \label{eqn:u} 
\end{align}
\end{subequations}
%%%%%%%%%%%%%%%%%%%%%%%%%%%%%%%%%%%%%%%%%%%%%%%%%%%%%%%%%%%%%%%%%%%%%%%
where $|I|$ and $|J|$ are the cardinalities of $I$ and $J$.  
We refer to $\Phi(x)$ as the {\sl multiplicative phase function} for the 
sum $g(n)$ in \eqref{eqn:g(n)}. 
%%%%%%
\par
%%%%%%
Thanks to positivity \eqref{eqn:positive} the function $\Phi(x)$ is smooth  
and positive on $(r_0, \, r_1)$. 
If we employ the convention $0^0 = 1$, which is natural in view of    
the limit $x^x = 1$ as $x \to +0$, then $\Phi(x)$ is continuous and positive 
at $x = r_0$ as well as at $x = r_1$ when $r_1 < +\infty$, even if some of the 
$l_i(x)$'s or $m_j(x)$'s vanish at one or both endpoints. 
When $r_1 = + \infty$, some calculations using balancedness condition 
\eqref{eqn:ap-bq2} shows that   
%%%%%%
\[
\Phi(x) = \left( \bsigma^{\sblambda} / \btau^{\sbmu} \right) \cdot 
\left( \bsigma^{\sbsigma} / \btau^{\sbtau} \right)^x \cdot 
\left \{ 1 + O(1/x) \right\} \quad \mbox{as} \quad x \to + \infty, 
\]
%%%%%%   
where $\bsigma^{\sbsigma} := \prod_{i \in I} \sigma_i^{\sigma_i}$, 
$\bsigma^{\sblambda} := \prod_{i \in I} \sigma_i^{\lambda_i}$ and 
so on; note that all of $\sigma_i$ and $\tau_j$ are positive due to the 
admissibility condition \eqref{eqn:ap-bq1} for the $r_1 = + \infty$ case. 
Thus it is natural to define  
%%%%%%%%%%%%%%%%%%%%%%%%% eqn:Phi(inf) %%%%%%%%%%%%%%%%%%%%%%%%%%%%%%%%%%
\begin{equation} \label{eqn:Phi(inf)}
\Phi(+ \infty) := 
\begin{cases}
0 \quad & (\mbox{if} \,\,\, \bsigma^{\sbsigma} < \btau^{\sbtau} ), \\[1mm]
\bsigma^{\sblambda} / \btau^{\sbmu} \quad & 
(\mbox{if} \,\,\,  \bsigma^{\sbsigma} = \btau^{\sbtau} ),  \\[1mm]
+\infty \quad & (\mbox{if} \,\,\, \bsigma^{\sbsigma} > \btau^{\sbtau} ).  
\end{cases}
\end{equation}
%%%%%%%%%%%%%%%%%%%%%%%%%%%%%%%%%%%%%%%%%%%%%%%%%%%%%%%%%%%%%%%%%%%%%%
With this understanding we assume the {\sl continuity at infinity}:    
%%%%%%%%%%%%%%%%%%%%%%%%%%%%% eqn:infty %%%%%%%%%%%%%%%%%%%%%%%%%%%%%%
\begin{equation} \label{eqn:infty}
\bsigma^{\sbsigma} \le \btau^{\sbtau} \qquad (\mbox{when} \,\,\, r_1 = + \infty).   
\end{equation}
%%%%%%%%%%%%%%%%%%%%%%%%%%%%%%%%%%%%%%%%%%%%%%%%%%%%%%%%%%%%%%%%%%%
Then $\Phi(x)$ is continuous on $[r_0, \, r_1]$ even when $r_1 = + \infty$ and 
it makes sense to define 
%%%%%
\[
\Phim := \max_{r_0 \le x \le r_1} \Phi(x),  
\]
%%%%%     
as a positive finite number. 
Therefore the function    
%%%%%%%%%%%%%%%%%%%%%%%%%%%%%%%%%% eqn:phi %%%%%%%%%%%%%%%%%%%%%%%%%
\begin{equation} \label{eqn:phi}
\phi(x) := - \log \Phi(x)  
\end{equation}
%%%%%%%%%%%%%%%%%%%%%%%%%%%%%%%%%%%%%%%%%%%%%%%%%%%%%%%%%%%%%%%%%%
is a real-valued, continuous function on $[r_0, \, r_1)$, smooth in $(r_0, \, r_1)$; 
if $r_1 < + \infty$ then it is also continuous at $x = r_1$; otherwise, $\phi(x)$ 
is either continuous at $x = + \infty$ or tends to $+ \infty$ as $x \to + \infty$.   
We refer to $\phi(x)$ as the {\sl additive phase function} for the sum $g(n)$ in 
\eqref{eqn:g(n)}.  
%%%%%%%
\par
%%%%%%%
When $r_1 = + \infty$ we have to think of the (absolute)  
{\sl convergence} of infinite series \eqref{eqn:g(n)}.  
If the strict inequality $\bsigma^{\sbsigma} < \btau^{\sbtau}$ holds in 
\eqref{eqn:infty} then it certainly converges. 
Otherwise, in order to guarantee its convergence, suppose that 
there is a constant $\sigma > 0$ such that for any $n \in \N$,    
%%%%%%%%%%%%%%%%%%%%%%%%%%% eqn:conv %%%%%%%%%%%%%%%%%%%%%%%%%%%%%%
\begin{equation} \label{eqn:conv}
\rRe \, \gamma(n) \le -1-\sigma  \qquad 
(\mbox{if $\bsigma^{\sbsigma} = \btau^{\sbtau}$}), 
\end{equation}
%%%%%%%%%%%%%%%%%%%%%%%%%%%%%%%%%%%%%%%%%%%%%%%%%%%%%%%%%%%%%%%%%
where
%%%%%%%%%%%%%%%%%%%%%%%%%%% eqn:ga(n) %%%%%%%%%%%%%%%%%%%%%%%%%%%%%
\begin{equation} \label{eqn:ga(n)} 
\gamma(n) := \sum_{i \in I} \alpha_i(n) - \sum_{j \in J} 
\beta_j(n) + \frac{|J|-|I|}{2}.    
\end{equation}
%%%%%%%%%%%%%%%%%%%%%%%%%%%%%%%%%%%%%%%%%%%%%%%%%%%%%%%%%%%%%%%%%%
\par
%%%%%%
Thanks to positivity \eqref{eqn:positive} the function $u(x; n)$ is also smooth 
and nowhere vanishing on $(r_0, \, r_1)$, but it may be singular at one or both 
ends of the interval when some of the $l_i(x)$'s or $m_j(x)$'s vanish there.   
To deal with this situation we say that $g(n)$ is {\sl left-regular} if 
$I_0 \cup J_0 = \emptyset$;  {\sl right-regular} if $I_1 \cup J_1 = \emptyset$; 
and {\sl regular} if $I_0 \cup J_0 \cup I_1 \cup J_1 = \emptyset$. 
If $g(n)$ is left-regular resp. right-regular with $r_1 < + \infty$, then 
$u(x; n)$ is continuous at $x = r_0$ resp. $x = r_1$. 
When $r_1 < + \infty$ the reflectional symmetry \eqref{eqn:refl} exchanges 
left and right regularities to each other.       
We remark that if $r_1 = + \infty$ then right-regularity automatically follows 
from admissibility. 
%%%%%
\par
%%%%%
The maximum of $\Phi(x)$ or equivalently the minimum of $\phi(x)$ plays 
a leading role in our analysis, so it is important to think of the first 
and second derivatives of $\phi(x)$. 
Differentiations of \eqref{eqn:phi} with balancedness condition 
\eqref{eqn:ap-bq2} took into account yield      
%%%%%%%%%%%%%%%%%%%%%%%%%%%%%%%% eqn:dash %%%%%%%%%%%%%%%%%%%%%%%%%%%%%%%
\begin{subequations} \label{eqn:dash}
\begin{align}
\phi'(x)  &= \log \prod_{j \in J} m_j(x)^{\tau_j} \prod_{i \in I} l_i(x)^{-\sigma_i}, 
\label{eqn:dash1} \\
\phi''(x) &= \sum_{j \in J} \dfrac{\tau_j^2}{m_j(x)} - \sum_{i \in I} 
\dfrac{\sigma_i^2}{l_i(x)}. \label{eqn:dash2} 
\end{align}
\end{subequations}   
%%%%%%%%%%%%%%%%%%%%%%%%%%%%%%%%%%%%%%%%%%%%%%%%%%%%%%%%%%%%%%%%%%%%%%%%
\par
%%%%%
Denote by $\Mm$ the set of all maximum points of $\Phi(x)$ on $[r_0, \, r_1]$. 
Suppose that $\Phi(x)$ attains its maximum $\Phim$ only within $(r_0, \, r_1)$, 
that is, $r_0$, $r_1 \not \in \Mm$.  
Moreover suppose that  every maximum point is {\sl nondegenerate} to the effect that  
%%%%%%%%%%%%%%%%%%%%%%%%%%%%%%%% eqn:max %%%%%%%%%%%%%%%%%%%%%%%%%%%%%%%%
\begin{equation} \label{eqn:max} 
\Mm \Subset (r_0, \, r_1), \qquad 
\phi''( x_0 ) > 0 \quad \mbox{at any} \quad x_0 \in \Mm,    
\end{equation}
%%%%%%%%%%%%%%%%%%%%%%%%%%%%%%%%%%%%%%%%%%%%%%%%%%%%%%%%%%%%%%%%%%%%%%%%
which is referred to as {\sl properness} of the maximum.         
By formula \eqref{eqn:dash1} any $x \in \Mm$ is a root of   
%%%%%%%%%%%%%%%%%%%%%%%%%%%%%%%% eqn:max-eq %%%%%%%%%%%%%%%%%%%%%%%%%%%%%
\begin{equation} \label{eqn:max-eq}
\chi( x ) :=  \prod_{j \in J} m_j( x )^{\tau_j} - \prod_{i \in I} l_i( x )^{\sigma_i} = 0,   
\qquad x \in (r_0, \, r_1),   
\end{equation}
%%%%%%%%%%%%%%%%%%%%%%%%%%%%%%%%%%%%%%%%%%%%%%%%%%%%%%%%%%%%%%%%%%%%%%% 
which is called the {\sl characteristic equation} for $g(n)$, while $\chi( x )$ 
is referred to as the {\sl characteristic function} for $g(n)$.  
It is easy to see that equation \eqref{eqn:max-eq} has only a finite number of roots, 
unless $\chi(x) \equiv 0$, so $\Mm$ must be a finite set.      
Note that $\phi'(x)$ and $\chi(x)$ have the same sign.   
%%%%%%%%%%
\par
%%%%%%%%%% 
Equation \eqref{eqn:max-eq} can be used to determine the set $\Mm$ explicitly. 
In applications to hypergeometric series, one usually puts $\sigma_i$, $\tau_j = 
\pm 1$ 
and $\lambda_i$, $\mu_j \in \Z$, thus \eqref{eqn:max-eq} is equivalent to an 
algebraic equation with integer coefficients and hence any $x \in \Mm$ 
must be an algebraic number. 
In this case with $r_1 = +\infty$, since $\bsigma^{\sbsigma} = \btau^{\sbtau} = 
\bsigma^{\sblambda} = \btau^{\sbmu} = 1$, the continuity at infinity \eqref{eqn:infty} 
is trivially satisfied with $\Phi(+ \infty) = 1$ in \eqref{eqn:Phi(inf)}, thus condition 
$\Mm \Subset (r_0, \, +\infty)$ in \eqref{eqn:max} includes $\Phim > 1$.  
%%%%%%%%%%%%%%%%%%%%%%%%%%%%%%%% thm:dspm %%%%%%%%%%%%%%%%%%%%%%%%%%%%%%
\begin{theorem} \label{thm:dspm}  
If balancedness \eqref{eqn:ap-bq2}, boundedness \eqref{eqn:bdd}, 
admissibility \eqref{eqn:ap-bq1}, genericness \eqref{eqn:punc} and properness 
\eqref{eqn:max} are all satisfied, with continuity at infinity \eqref{eqn:infty} 
and convergence \eqref{eqn:conv} being added when $r_1 = + \infty$, then the 
sum $g(n)$ in \eqref{eqn:g(n)} admits an asymptotic representation    
%%%%%%%%%%%%%%%%%%%%%%%%%%%%%%%% eqn:dspm %%%%%%%%%%%%%%%%%%%%%%%%%%%%%%
\begin{equation} \label{eqn:dspm}
g(n) = n^{\gamma(n) + \frac{1}{2} } \cdot 
\Phim^{\, n} \cdot \left\{\, C(n)  + \Omega(n)  \right\},     
\end{equation} 
%%%%%%%%%%%%%%%%%%%%%%%%%%%%%%%%%%%%%%%%%%%%%%%%%%%%%%%%%%%%%%%%%%%%%%%
where $\gamma(n)$ is defined in formula \eqref{eqn:ga(n)} while $\Phim$ and 
$C(n)$ are defined by 
%%%%%%%%%%%%%%%%%%%%%%%%%%%%%%% eqn:PhimC %%%%%%%%%%%%%%%%%%%%%%%%%%%%%%
\begin{equation} \label{eqn:PhimC}
\Phim := \Phi(x_0) \quad \mbox{for any} \quad x_0 \in \Mm; \qquad 
C(n) := \sqrt{2 \pi } \sum_{x_0 \in \Mm} \dfrac{u( x_0; n)}{ \sqrt{\phi''( x_0 )}}  
\end{equation} 
%%%%%%%%%%%%%%%%%%%%%%%%%%%%%%%%%%%%%%%%%%%%%%%%%%%%%%%%%%%%%%%%%%%%%%%%
in terms of the notations in \eqref{eqn:G-asymp2}, whereas the error term 
$\Omega(n)$ is estimated as  
%%%%%%%%%%%%%%%%%%%%%%%% eqn:uniform %%%%%%%%%%%%%%%%%%%%%%%%%%%%%%%%%%%
\begin{equation} \label{eqn:uniform}
|\Omega(n)| \le K \{ n^{-\frac{1}{2}} + 
\lambda^{-n} (\delta_0(n)^{-1} + \delta_1(n)^{-1}) \},  
\qquad {}^{\forall} n \ge N,  
\end{equation}
%%%%%%%%%%%%%%%%%%%%%%%%%%%%%%%%%%%%%%%%%%%%%%%%%%%%%%%%%%%%%%%%%%%%%%
for some constants $K > 0$, $\lambda > 1$ and $N \in \N$, where $\delta_0(n)$ 
and $\delta_1(n)$ are defined in \eqref{eqn:punc}.  
This estimate is valid uniformly for all $\bal(n)$ and $\bbeta(n)$ satisfying 
conditions \eqref{eqn:bdd} and \eqref{eqn:punc} along with 
\eqref{eqn:conv} when $r_1 = +\infty$, in which case $I_1 = \emptyset$ 
and so $\delta_1(n) = 1$.    
\end{theorem} 
%%%%%%%%%%%%%%%%%%%%%%%%%%%%%%%%%%%%%%%%%%%%%%%%%%%%%%%%%%%%%%%%%%%%%%%
\par
%%%%%
Things are simpler when $\Mm$ consists of a single point $x_0 \in (r_0, \, r_1)$,  
in which case the main idea for proving Theorem \ref{thm:dspm} is to divide 
the sum \eqref{eqn:g(n)} into five components: 
%%%%%%%%%%%%%%%%%%%%%%%%
\begin{equation*}  
g(n) = g_0(n) + h_0(n) + h(n) + h_1(n) + g_1(n),      
\end{equation*}
%%%%%%%%%%%%%%%%%%%%%%%
with each component being a partial sum of \eqref{eqn:g(n)} defined by  
%%%%%%%%%%%%%%%%%%%%%%%%%%%%%% eqn:divide %%%%%%%%%%%%%%%%%%%%%%%%%%%%%%%%
\begin{alignat}{3}
g_0(n) &:= \mbox{sum of $G(k; n)$ over} \,\, & \lceil r_0 n \rceil & \le k \le \lceil (r_0 + \ve) n \rceil-1, \quad 
& & \mbox{(left end)} \nonumber \\ 
h_0(n) &:= \mbox{sum of $G(k; n)$ over} \,\, & \lceil (r_0 + \ve) n \rceil & \le k \le \lceil (x_0 -\ve) n \rceil-1, \quad 
& & \mbox{(left side)} \nonumber \\
h(n) &:= \mbox{sum of $G(k; n)$ over} \,\, & \lceil (x_0-\ve) n \rceil & \le k \le \lceil (x_0 +\ve) n \rceil-1, \quad 
& & \mbox{(top)}  \label{eqn:divide} \\
h_1(n) &:= \mbox{sum of $G(k; n)$ over} \,\, & \lceil (x_0 +\ve) n \rceil & \le k \le \lceil (r_1-\ve) n \rceil-1, \quad 
& & \mbox{(right side)} \nonumber \\
g_1(n) &:= \mbox{sum of $G(k; n)$ over} \,\, & \lceil (r_1 -\ve) n & \le k \le \lceil r_1 n \rceil-1,  \quad 
& & \mbox{(right end)} \nonumber 
\end{alignat}
%%%%%%%%%%%%%%%%%%%%%%%%%%%%%%%%%%%%%%%%%%%%%%%%%%%%%%%%%%%%%%%%%%%%%%%%
where if $r_1 = + \infty$ then the right-end component should be omitted.    
In order for the division \eqref{eqn:divide} to make sense, the number $\ve$ must satisfy     
%%%%%%%%%%%%%%%%%%%%%%%%%%%%%% eqn:ve %%%%%%%%%%%%%%%%%%%%%%%%%%%%%%%%%%%
\begin{equation} \label{eqn:ve}
0 < \ve < \ve_0 := \min \{ (x_0 -r_0)/2, \, (r_1-x_0)/2 \}. 
\end{equation}
%%%%%%%%%%%%%%%%%%%%%%%%%%%%%%%%%%%%%%%%%%%%%%%%%%%%%%%%%%%%%%%%%%%%%%%%
How to take $\ve \in (0, \, \ve_0)$ will be specified in the course of establishing 
Theorem \ref{thm:dspm}. 
%%%%%
\par
%%%%%
We want to think of $h(n)$ as the principal part of $g(n)$, while other 
four components as remainders.  
Thus estimating the top component $h(n)$ is the central issue of this section, 
but treatment of both ends $g_0(n)$ and $g_1(n)$ is also far from trivial.       
For the sake of simplicity we shall deal with the case $|\Mm| = 1$ only, 
but even when $|\Mm| \ge 2$ things are essentially the same and it 
will be clear how to modify the arguments.  
The reflectional symmetry \eqref{eqn:refl} reduces the discussion 
at the right end or right side to the discussion at the left counterpart.     
The top and side sums are regular, so we shall begin by estimating regular 
sums in \S \ref{subsec:regular}. 
%%%%%%
\par
%%%%%%   
In the sequel we shall often utilize the following version of Stirling's formula:    
For any positive number $c > 0$ and any compact subset $A \Subset \C$ we have     
%%%%%%%%%%%%%%%%%%%%%%%%%%%%%%%% eqn:stirling %%%%%%%%%%%%%%%%%%%%%%%%%%%%%%
\begin{equation} \label{eqn:stirling} 
\vG( x n + a ) = (2 \pi)^{1/2} \, x^{a- \frac{1}{2} } \, n^{a- \frac{1}{2} } \, x^{x n} \, (n/e)^{x n} \,  
\left\{ \, 1 + O( 1/n ) \, \right\} \quad \mbox{as} \quad n \to +\infty, 
\end{equation}
%%%%%%%%%%%%%%%%%%%%%%%%%%%%%%%%%%%%%%%%%%%%%%%%%%%%%%%%%%%%%%%%%%%%%%%%
where Landau's symbol $O(1/n)$ is uniform with respect to 
$(x, a) \in \R_{\ge c} \times A$.  
%%%%%%%%%%%%%%%%%%%%%%%% subsect:regular %%%%%%%%%%%%%%%%%%%%%%%%%%%%%%%%%%
\subsection{Regular Sums and Side Components} \label{subsec:regular} 
%%%%%%%%%%%%%%%%%%%%%%%%%%%%%%%%%%%%%%%%%%%%%%%%%%%%%%%%%%%%%%%%%%%%%%%%
\par
%%%%%
In this subsection we assume that $g(n)$ in \eqref{eqn:g(n)} satisfies 
balancedness \eqref{eqn:ap-bq2}, boundedness \eqref{eqn:bdd} and admissibility 
\eqref{eqn:ap-bq1}, along with continuity at infinity \eqref{eqn:infty} and 
convergence \eqref{eqn:conv} if $r_1 = + \infty$, while 
properness \eqref{eqn:max} is not assumed  and 
genericness \eqref{eqn:punc} is irrelevant to regular sums.   
%%%%%%%%%%%%%%%%%%%%%%%%%% lem:H %%%%%%%%%%%%%%%%%%%%%%%%%%%%%%%%%%%%%%%%
\begin{lemma} \label{lem:H} 
If the sum $g(n)$ in \eqref{eqn:g(n)} is regular then there exists an integer 
$N_0 \in \N$ and a constant $C_0 > 0$ such that $H(x; n)$ in formula 
\eqref{eqn:H} can be written     
%%%%%%%%%%%%%%%%%%%%%%%% eqn:H-asymp %%%%%%%%%%%%%%%%%%%%%%%%%%%%%%%%%%%%
\begin{subequations} \label{eqn:H-asymp}
\begin{align}  
H(x; n) &= u(x; n)  \cdot n^{\gamma(n)} \cdot \Phi(x)^n \cdot \{1+ e(x; n) \}, 
\label{eqn:H-asymp1} \\[1mm]
|e( x; n)| &\le C_0/n,  \qquad 
{}^{\forall} n \ge N_0, \,\, r_0 \le {}^{\forall} x \le r_1.    
\label{eqn:H-asymp2}     
\end{align}
\end{subequations}
%%%%%%%%%%%%%%%%%%%%%%%%%%%%%%%%%%%%%%%%%%%%%%%%%%%%%%%%%%%%%%%%%%%%%%%%
\end{lemma} 
%%%%%%%%%%%%%%%%%%%%%%%% begin proof %%%%%%%%%%%%%%%%%%%%%%%%%%%%%%%%%%%%
{\it Proof}. 
Since $g(n)$ is regular, that is, $I_0 \cup I_1 \cup J_0 \cup J_1 = \emptyset$, we 
have the uniform positivity \eqref{eqn:positive2} for all $i \in I$, $j \in J$ and 
$x \in [r_0, \, r_1]$.   
This together with boundedness \eqref{eqn:bdd} allows us to apply Stirling's 
formula \eqref{eqn:stirling} to all gamma factors  
$\vG(l_i(x) n + \alpha_i(n))$ and $\vG(m_j(x) n + \beta_j(n))$ of $H(x; n)$ in  \eqref{eqn:H}. 
Taking definitions \eqref{eqn:G-asymp2} and \eqref{eqn:ga(n)} into account  
we use formula \eqref{eqn:stirling} to have 
%%%%%
\[
H(x; n) = u(x; n) \cdot n^{\gamma(n)} \cdot \Phi(x)^n \cdot
\left\{ (n/e)^{\sum_{i \in I} l_i(x) - \sum_{j \in J} m_j(x)} \right \}^n \cdot 
\{1 + O(1/n) \}, 
\]
%%%%%
where the $O(1/n)$ term is uniform with respect to $x \in [r_0, \, r_1]$ as 
well as to $\bal(n)$ and $\bbeta(n)$ satisfying condition \eqref{eqn:bdd}.  
Then balancedness \eqref{eqn:balanced} yields the desired 
formula \eqref{eqn:H-asymp}. \hfill $\Box$ 
%%%%%%%%%%%%%%%%%%%%%%%% end proof %%%%%%%%%%%%%%%%%%%%%%%%%%%%%%%%%%%%%%%
%%%%%%%%%%%%%%%%%%%%%%%% prop:reg-bd %%%%%%%%%%%%%%%%%%%%%%%%%%%%%%%%%%%%%
\begin{proposition} \label{prop:reg-bd} 
If the sum $g(n)$ in formula \eqref{eqn:g(n)} is regular then it admits an estimate  
%%%%%%%%%%%%%%%%%%%%%%%%
\begin{equation*} 
|g(n)| \le C_1 \cdot n^{\rRe \, \gamma(n) + 1} \cdot \Phim^{\, n}, \qquad 
{}^{\forall} n \ge N_0,   
\end{equation*} 
%%%%%%%%%%%%%%%%%%%%%%%
for a constant $C_1 > 0$ and an integer $N_0 \in \N$ which is the same as in 
Lemma $\ref{lem:H}$.    
\end{proposition} 
%%%%%%%%%%%%%%%%%%%%%%%% begin proof %%%%%%%%%%%%%%%%%%%%%%%%%%%%%%%%%%%
{\it Proof}. 
From representation \eqref{eqn:H-asymp} we have 
%%%%%%%%%%%%%%%%%%%%%%% eqn:reg-bd-es1 %%%%%%%%%%%%%%%%%%%%%%%%%%%%%%%%%
\begin{equation} \label{eqn:reg-bd-es1}
|H(x; n)| \le (1+C_0) \cdot 
|u(x; n)|\cdot n^{\rRe \, \gamma(n)} \cdot \Phim^{\, n}, \qquad  
r_0 \le {}^{\forall} x \le r_1, \,\, {}^{\forall} n \ge N_0.  
\end{equation}
%%%%%%%%%%%%%%%%%%%%%%%%%%%%%%%%%%%%%%%%%%%%%%%%%%%%%%%%%%%%%%%%%%%%%% 
\par
%%%%%%%
First we consider the case $r_1 < + \infty$. 
Since $g(n)$ is regular and $\bal(n)$ and $\bbeta(n)$ are bounded by 
assumption \eqref{eqn:bdd}, the definition \eqref{eqn:u} implies that $u(x; n)$ 
is bounded for $(x, n) \in [r_0, \, r_1] \times \Z_{\ge N_0}$. 
Replacing the constant $C_0$ by a larger one if necessary, we have 
$|H(x; n)| \le C_0 \cdot n^{\rRe \, \gamma(n)} 
\cdot \Phim^{\, n}$ for any $x \in [r_0, \, r_1]$ and $n \ge N_0$. 
Thus by definitions \eqref{eqn:g(n)} and \eqref{eqn:H} we have for any 
$n \ge N_0$,   
%%%%%
\begin{align*}
|g(n)| & \le \sum_{k=\lceil r_0 n \rceil}^{\lceil r_1 n \rceil -1} |H(k/n; n)| 
\le C_0 \cdot n^{\rRe \, \gamma(n)} \cdot \Phim^{\, n} 
\sum_{k=\lceil r_0 n \rceil}^{\lceil r_1 n \rceil -1} 1 \\ 
& = C_0 \cdot n^{\rRe \, \gamma(n)} \cdot \Phim^{\, n} \cdot 
(\lceil r_1 n \rceil - \lceil r_0 n \rceil) 
\le C_1 \cdot n^{\rRe \, \gamma(n) + 1} \cdot \Phim^{\, n},  
\end{align*}
%%%%%%%%%%
with the constant $C_1 := C_0 (1+r_1-r_0)$.   
%%%%%
\par
%%%%%
We proceed to the case $r_1 = + \infty$ and $\bsigma^{\sbsigma} = 
\btau^{\sbtau}$ in which condition \eqref{eqn:conv} takes place.    
Since $g(n)$ is regular and $\bal(n)$ and $\bbeta(n)$ are bounded by  
\eqref{eqn:bdd}, the definition \eqref{eqn:u} implies that 
%%%%%%
\[
u(x; n) = 
(2 \pi)^{ \frac{|I|-|J|}{2} } \, 
\prod_{i \in I} \sigma_i^{\alpha_i(n) - \frac{1}{2}} 
\prod_{j \in J} \tau_j^{\frac{1}{2}-\beta_j(n)} \cdot x^{\gamma(n)} 
\cdot \{ 1 + O(1/x) \} \qquad \mbox{as} \quad x \to + \infty, 
\]
%%%%%%
uniformly for $n \in \N$. 
By condition \eqref{eqn:conv} there exists a constant $C_2 > 0$ 
such that       
%%%%%%%%%%%
\[
|u(x; n)| \le C_2 \, (2+x)^{-1-\sigma}, \qquad 
{}^{\forall} x \ge r_0, \,\, {}^{\forall} n \ge N_0.     
\]
%%%%%%%%%%
In view of definitions \eqref{eqn:g(n)} and \eqref{eqn:H} this estimate 
together with formula \eqref{eqn:reg-bd-es1} yields 
%%%%%%%
\begin{align*}
|g(n)| &\le \sum_{k=\lceil r_0 n \rceil}^{\infty} |H(k/n; n)| 
\le C_2 \, (1+C_0) \cdot n^{\rRe \, \gamma(n) +1} \cdot \Phim^{\, n} 
\sum_{k=\lceil r_0 n \rceil}^{\infty} 
\left( 2+\frac{k}{n} \right)^{-1-\sigma} \frac{1}{n} \\
& \le C_2 \, (1+C_0) \cdot n^{\rRe \, \gamma(n) +1} \cdot \Phim^{\, n} 
\int_{r_0}^{\infty} \left( 1+ x \right)^{-1-\sigma} \, d x = C_1 \cdot 
n^{\rRe \, \gamma(n) + 1} \cdot \Phim^{\, n}, 
\end{align*}
%%%%%%%
for any integer $n \ge N_0$, where $C_1 := C_2 \, (1+C_0)\, 
(1+r_0)^{-\sigma}/\sigma $.  
%%%%%%%
\par
%%%%%%%
The proof ends with the case where $r_1 = + \infty$ and 
$\bsigma^{\sbsigma} < \btau^{\sbtau}$. 
By Stirling's formula \eqref{eqn:stirling} there exists a constant 
$C_3 > 0$ such that $|H(x; n)| \le C_3 \cdot (x n)^{\rRe \, \gamma(n)} 
\cdot \Phi(x)^n$ and $\Phi(x) \le C_3 \cdot \rho^x$ with 
$\rho := \bsigma^{\sbsigma}/\btau^{\sbtau} < 1$  
for any $x \ge r_0$ and $n \ge N_0$. 
Take a number $r_2 > r_0$ so large that 
$d := C_3 \cdot \rho^{r_2/2} < \Phim$ and let $g(n) = g_1(n) + g_2(n)$ 
be the decomposition according to the division 
$[r_0, \,+\infty) = [r_0, \, r_2) \cup [r_2, \, +\infty)$.  
Then an estimate for the $r_1 < + \infty$ case applies to 
$g_1(n)$, while $|H(x; n)| \le C_3 \cdot d^n \cdot (x n)^c \cdot 
\rho^{x n/2}$ for $x \ge r_2$, where $c := \sup_{n \ge N_0} 
\rRe \, \gamma(n)$, and hence     
%%%%%%%
\[
|g_2(n)| \le \sum_{k=\lceil r_2 n \rceil}^{\infty} |H(k/n; n)| \le 
C_3 \cdot d^n \sum_{k = \lceil r_2 n \rceil}^{\infty} k^c \cdot \rho^{k/2} 
\le C_3 \cdot d^n \sum_{k = 1}^{\infty} k^c \cdot \rho^{k/2} 
= C_4 \cdot d^n   
\]
%%%%%%%   
for any $n \ge N_0$. 
It is clear from $0 < d < \Phim$ that the proposition follows. 
\hfill $\Box$ \par\medskip
%%%%%%%%%%%%%%%%%%%%%%%% end proof %%%%%%%%%%%%%%%%%%%%%%%%%%%%%%%%%%%%%
Proposition \ref{prop:reg-bd} can be used to estimate the side components $h_0(n)$ 
and $h_1(n)$ in \eqref{eqn:divide}. 
%%%%%%%%%%%%%%%%%%%%%%%% lem:side %%%%%%%%%%%%%%%%%%%%%%%%%%%%%%%%%%%%%%
\begin{lemma} \label{lem:side} 
For any $0 < \ve < \ve_0$ there exist $N_1^{\ve} \in \N$ and $C_1^{\ve} > 0$ 
such that  
%%%%%
\[
| h_0(n) | \le C_1^{\ve} \cdot n^{\rRe \, \gamma(n) + 1} \cdot (\Phi_0^{\ve})^n, 
\quad 
| h_1(n) | \le C_1^{\ve} \cdot n^{\rRe \, \gamma(n) + 1} \cdot (\Phi_1^{\ve})^n,  
\quad 
{}^{\forall} n \ge N_1^{\ve},  
\]
%%%%% 
where $\Phi_0^{\ve} := \displaystyle \max_{r_0 + \ve \le x \le x_0-\ve} \Phi(x)$ 
and $\Phi_1^{\ve} := \displaystyle \max_{x_0 + \ve \le x \le r_1-\ve} \Phi(x)$.   
\end{lemma}
%%%%%%%%%%%%%%%%%%%%%%%% begin proof %%%%%%%%%%%%%%%%%%%%%%%%%%%%%%%%%%%%
{\it Proof}. 
We have only to apply Proposition \ref{prop:reg-bd} with $r_0$ and $r_1$ 
replaced by $r_0 + \ve$ and $x_0-\ve$ to deduce the estimate for $h_0(n)$.  
In a similar manner we apply the proposition this time with $r_0$ and $r_1$ 
replaced by $x_0 + \ve$ and $r_1-\ve$ to get the estimate for $h_1(n)$. 
\hfill $\Box$ 
%%%%%%%%%%%%%%%%%%%%%%%% subsec:t-s %%%%%%%%%%%%%%%%%%%%%%%%%%%%%%%%%%%%%
\subsection{Top Component} \label{subsec:t-s} 
%%%%%%%%%%%%%%%%%%%%%%%%%%%%%%%%%%%%%%%%%%%%%%%%%%%%%%%%%%%%%%%%%%%%%%% 
We consider the top component $h(n)$ in \eqref{eqn:divide}. 
Recall the setting in \S \ref{subsec:results} that $\Mm = \{ x_0 \} \Subset 
(r_0, \, r_1)$, $\Phim = \Phi( x_0 ) = e^{- \phi( x_0)}$, $\phi'( x_0 ) = 0$ and 
$\phi''( x_0 ) > 0$. 
Since the sum $h(n)$ is regular, Lemma \ref{lem:H} implies that $H(x; n)$ 
can be written as in \eqref{eqn:H-asymp1} with 
estimate \eqref{eqn:H-asymp2} now being 
%%%%%%%%%%%%%%%%%%%%%%%%% eqn:C0(ve) %%%%%%%%%%%%%%%%%%%%%%%%%%%%%%%%%%%%
\begin{equation} \label{eqn:C0(ve)}
|e( x; n)| \le C_0(\ve)/n, \qquad 
{}^{\forall} n \ge N_0(\ve), \,\, x_0-\ve \le {}^{\forall} x \le x_0+\ve.    
\end{equation}
%%%%%%%%%%%%%%%%%%%%%%%%%%%%%%%%%%%%%%%%%%%%%%%%%%%%%%%%%%%%%%%%%%%%%%%%
\par
%%%%%    
The local study of $H(x; n)$ near $x = x_0 $ is best performed in terms of 
new variables  
%%%%%%%%%%%%%%%%%%%%%%%%%%
\begin{equation*} 
y := x- x_0 \quad \mbox{(shift)}; \qquad 
z := \sqrt{n} \, y \quad \mbox{(scale change)}. 
\end{equation*} 
%%%%%%%%%%%%%%%%%%%%%%%
Taylor expansions around $x = x_0 $ show that $\phi( x )$ and $u (x; n) $ can 
be written 
%%%%%%%%%%%%%%%%%%%%%%%% eqn:taylor %%%%%%%%%%%%%%%%%%%%%%%%%%%%%%%%%%%%
\begin{subequations}
\begin{alignat}{3}
\phi(x) &= \phi( x_0 ) + a \, y^2 + \eta( y ), \quad &  
|\eta( y )| &\le b \, |y|^3, \qquad & |{}^{\forall} y \, | &\le \ve_1, 
\label{eqn:taylor-p} \\
u(x; n) &= u( x_0; n) + v( y; n ), \quad &  
|v( y; n )| &\le c \, |y|, \qquad & |{}^{\forall} y \, | &\le \ve_1,  
\label{eqn:taylor-u}    
\end{alignat}
\end{subequations}
%%%%%%%%%%%%%%%%%%%%%%%%%%%%%%%%%%%%%%%%%%%%%%%%%%%%%%%%%%%%%%%%%%%%%%%%
with $a := \ts \frac{1}{2} \, \phi''( x_0 ) > 0$ and some positive constants $b$, $c$, 
$\ve_1 > 0$.   
It is clear that $a$ and $b$ are independent of $n$.  
We can also take $c$ and $\ve_1$ uniformly in $n$ because 
$\bal(n)$ and $\bbeta(n)$ are bounded by assumption \eqref{eqn:bdd}. 
If we put 
%%%%%%%%%%%%%%%%%%%%%%%% eqn:Habc %%%%%%%%%%%%%%%%%%%%%%%%%%%%%%%%%%%%%%%
\begin{subequations}
\begin{align}
H_{\ra}(x; n) &:= u(x_0; n) \cdot n^{\gamma(n)} \cdot \Phi( x )^n 
= u( x_0; n ) \cdot n^{\gamma(n)} \cdot \Phim^{\, n} \cdot 
e^{-n \{ a \, y^2 + \eta( y ) \}}, 
\label{eqn:Ha} \\
H_{\rb}(x; n) &:= v( y; n ) \cdot n^{\gamma(n)} \cdot \Phi( x )^n  =  
n^{\gamma(n)} \cdot \Phim^{\, n} \cdot v( y; n ) \cdot 
e^{-n \{ a \, y^2 + \eta( y ) \}}, \label{eqn:Hb} \\
H_{\rc}(x; n) &:= u(x; n) \cdot n^{\gamma(n)} \cdot \Phi(x)^n \cdot e(x; n),      
\label{eqn:Hc} 
\end{align}
\end{subequations}   
%%%%%%%%%%%%%%%%%%%%%%%%%%%%%%%%%%%%%%%%%%%%%%%%%%%%%%%%%%%%%%%%%%%%%%%%
then formula \eqref{eqn:H-asymp1} yields 
$H(x; n) = H_{\ra}(x; n) + H_{\rb}(x; n) + H_{\rc}(x; n)$, which in turn gives 
%%%%%%
\[
h(n) = h_{\ra}(n) + h_{\rb}(n) + h_{\rc}(n), \qquad 
h_{\nu}(n) := \sum_{k = l}^{m-1} H_{\nu}(k/n; n), \quad 
\nu = \ra, \rb, \rc,  
\]  
%%%%%%
where  $l := \lceil (x_0-\ve) n \rceil$ and $m := \lceil (x_0 + \ve) n \rceil$. 
%%%%%%
\par
%%%%%%
To estimate $h_{\ra}(n)$ we use some a priori estimates, which will be collected  
in \S \ref{subsec:apriori}.    
%%%%%%%%%%%%%%%%%%%%%%%%%% lem:ha %%%%%%%%%%%%%%%%%%%%%%%%%%%%%%%%%%%%%%%
\begin{lemma} \label{lem:ha}  
For any $0 < \ve < \ve_2 := \min\{ \ve_0, \, \frac{\ve_1}{2}, \, \frac{a}{4 b} \}$ 
and $n \ge N_1(\ve) := \max\{ 2/\ve, \, N_0(\ve) \}$,   
%%%%%%%%%%%%%%%%%%%%%%%%% eqn:ha %%%%%%%%%%%%%%%%%%%%%%%%%%%%%%%%%%%%%%%%
\begin{subequations} \label{eqn:ha}
\begin{gather} 
h_{\ra}(n) = \sqrt{\pi/a} \, \cdot u(x_0; n) \cdot 
n^{\gamma(n) +1/2} \cdot \Phim^{\, n} \cdot \left\{ 1 + e_{\ra}(n) \cdot n^{-1/2} \right \}, 
\label{eqn:ha1} \\[2mm]
|e_{\ra}(n)| \le M_5(a, b; \ve) := 2 M_3(a, b) + (5/a) \cdot (2 \ve)^{-3/2},   \label{eqn:ha2}
\end{gather}
\end{subequations}
%%%%%%%%%%%%%%%%%%%%%%%%%%%%%%%%%%%%%%%%%%%%%%%%%%%%%%%%%%%%%%%%%%%%%%%%
where $M_3(a, b)$ is defined in Lemma $\ref{lem:approxi}$ and currently  
$a := \frac{1}{2} \phi''(x_0) > 0$.   
\end{lemma}
%%%%%%%%%%%%%%%%%%%%%%%%% begin proof %%%%%%%%%%%%%%%%%%%%%%%%%%%%%%%%%%%%
{\it Proof}. 
Put $\psi(z; a) := e^{- a \, z^2 + \delta(z)}$ with 
$\delta(z) := - n \cdot \eta\left(n^{-1/2} z \right)$. 
Then \eqref{eqn:Ha} and \eqref{eqn:taylor-p} read  
%%%%%%%%%%%%%%%%%%%%%%%%%% eqn:Ha-e %%%%%%%%%%%%%%%%%%%%%%%%%%%%%%%%%%%%%
\begin{subequations} \label{eqn:Ha-e}
\begin{align}
H_{\ra}(x; n) &= u(x_0; n) \cdot n^{\gamma(n) + 1/2} \cdot  \Phim^{\, n} \cdot 
\psi(z; a) \cdot \frac{1}{\sqrt{n}},   \label{eqn:Ha-e1} \\
|\delta(z)| &\le \frac{b}{\sqrt{n}} \, |z|^3 \qquad 
(|{}^{\forall} z \, | \le \ve_1 \sqrt{n} \,).  \label{eqn:Ha-e2}
\end{align}  
\end{subequations}
%%%%%%%%%%%%%%%%%%%%%%%%%%%%%%%%%%%%%%%%%%%%%%%%%%%%%%%%%%%%%%%%%%%%%%%%
Consider the sequence $\vD : \xi_k := (k - x_0 n)/\sqrt{n}$ $(k = l, \dots, m)$.        
From the definitions of $l$ and $m$,    
%%%%%%%%%%%%%%%%%%%%%%%%%% eqn:range %%%%%%%%%%%%%%%%%%%%%%%%%%%%%%%%%%%%
\begin{equation} \label{eqn:range} 
- \ve \sqrt{n} \le \xi_l < -\ve \sqrt{n} + 1/\sqrt{n}, \qquad 
\ve \sqrt{n} \le \xi_m < \ve \sqrt{n} + 1/\sqrt{n},   
\end{equation}
%%%%%%%%%%%%%%%%%%%%%%%%%%%%%%%%%%%%%%%%%%%%%%%%%%%%%%%%%%%%%%%%%%%%%%%%
which together with $0 < \ve < \ve_2$ and $n \ge N_1(\ve)$ implies inclusion 
$[\xi_l, \, \xi_m] \subset [- \ve_1 \sqrt{n}, \, \ve_1 \sqrt{n}]$, so the estimate 
\eqref{eqn:Ha-e2} is available for all $z \in [\xi_l, \, \xi_m]$.  
From formula \eqref{eqn:Ha-e1} we have  
%%%%%%%%%%%%%%%%%%%%%%%%%% eqn:R-sum %%%%%%%%%%%%%%%%%%%%%%%%%%%%%%%%%%%%
\begin{equation} \label{eqn:R-sum}
h_{\ra}(n) = u(x_0; n) \cdot n^{\gamma(n) + \frac{1}{2} } \cdot \Phim^{\, n} \cdot R(\psi, \vD) 
\quad \mbox{with} \quad 
R(\psi; \vD) := \sum_{k=l}^{m-1} \psi(\xi_k; a) \frac{1}{\sqrt{n}},     
\end{equation}
%%%%%%%%%%%%%%%%%%%%%%%%%%%%%%%%%%%%%%%%%%%%%%%%%%%%%%%%%%%%%%%%%%%%%%%%
where $R(\psi; \vD)$ is the left Riemann sum of $\psi(z; a)$ 
for equipartition $\vD$ of the interval $[\xi_l, \, \xi_m]$.  
%%%%%%%
\par
%%%%%%%     
Let $\varphi(z; a) := e^{-a z^2}$.  
Since $|\xi_k - z| \le 1/\sqrt{n}$ for any $z \in [\xi_k, \, \xi_{k+1}]$, 
Lemma \ref{lem:approxi} yields   
%%%%%%%
\[
\begin{split}
& \left| R(\psi; \vD) - \int_{\xi_l}^{\xi_m} \varphi(z; a) \, d z \right|
= \left| \sum_{k=l}^{m-1} \int_{\xi_k}^{\xi_{k+1}} 
\{ \psi(\xi_k; a) - \varphi(z; a) \} \, d z \right| \\
& \quad \le \sum_{k=l}^{m-1} \int_{\xi_k}^{\xi_{k+1}} 
\left| \psi(\xi_k; a) - \varphi(z; a) \right| \, d z 
\le \frac{M_3(a, b)}{\sqrt{n}} 
\sum_{k=l}^{m-1} \int_{\xi_k}^{\xi_{k+1}}\varphi(z; a/4) \, d z \\
& \quad =  \frac{M_3(a, b)}{\sqrt{n}} 
\int_{\xi_l}^{\xi_m} \varphi(z; a/4) \, d z 
\le \dfrac{M_3(a, b)}{\sqrt{n}} 
\ds \int_{-\infty}^{\infty} \varphi(z; a/4) \, d z 
= 2 M_3(a, b) \sqrt{\dfrac{\pi}{a n}},  
\end{split}
\]
%%%%%%
where estimate \eqref{eqn:approxi} is used in the second inequality. 
By the partition of Gaussian integral 
%%%%%%
\[
\sqrt{\pi/a} = \int_{-\infty}^{\infty} \varphi(z; a) \, 
d z = \int_{-\infty}^{\xi_l} + \int_{\xi_l}^{\xi_m} + 
\int_{\xi_m}^{\infty} \varphi(z; a) \, d z,  
\]  
%%%%%% 
and bounds $\xi_l \le -\ve \sqrt{n}/2$ and $\xi_m \ge \ve \sqrt{n}$, 
which follow from \eqref{eqn:range} and $n \ge 2/\ve$, we have  
%%%%%%%%%%%%%%%%%%%%%% eqn:R-es %%%%%%%%%%%%%%%%%%%%%%%%%%%%%
\begin{align}
\left| R(\psi; \vD) - \sqrt{\pi/a} \, \right| 
& \le \int_{-\infty}^{-\ve \sqrt{n}/2} \varphi(z; a) \, d z + 
2 M_3(a, b) \sqrt{\dfrac{\pi}{a n}} + 
\int_{\ve \sqrt{n}}^{\infty} \varphi(z; a) \, d z \nonumber \\
& \le 2 M_3(a, b) \sqrt{\dfrac{\pi}{a n}} + 
\frac{5 \sqrt{\pi}}{2 a^{3/2} \ve^2 \cdot n} \le 
\sqrt{\frac{\pi}{a n}} M_5(a, b; \ve),  \label{eqn:R-es}
\end{align}
%%%%%%%%%%%%%%%%%%%%%%%%%%%%%%%%%%%%%%%%%%%%%%%%%%%%%%%%%%%%%
with $M_5(a, b; \ve) := 2 M_3(a, b) + (5/a) \cdot (2 \ve)^{-3/2}$, 
where the estimate  
%%%%%%
\[
\int_z^{\infty} \varphi(t; a) \, d t \le 
\frac{1}{2} \sqrt{\frac{\pi}{a}} \, \varphi(z; a) 
\le \frac{\sqrt{\pi}}{2 a^{3/2} z^2},  \qquad 
{}^{\forall} z \ge 0,   
\]
%%%%%%
and $\sqrt{n} \ge \sqrt{2/\ve}$ are used in the second and third 
inequalities respectively. 
Upon writing $R(\psi; \vD) = \sqrt{\pi/a} \, \{ 1 + e_2(n) \cdot n^{-1/2} \}$, 
formula \eqref{eqn:ha} follows from \eqref{eqn:R-sum} and 
\eqref{eqn:R-es}. \hfill $\Box$  
%%%%%%%%%%%%%%%%%%%%%%%%%% end proof %%%%%%%%%%%%%%%%%%%%%%%%%%%%%%%%%%%%%
%%%%%%%%%%%%%%%%%%%%%%%%%% lem:hbc %%%%%%%%%%%%%%%%%%%%%%%%%%%%%%%%%%%%%%
\begin{lemma} \label{lem:hbc}  
For any $0 < \ve < \ve_2$ and $n \ge N_1(\ve)$ we have    
%%%%%%%%%%%%%%%%%%%%%%%%% eqn:hb %%%%%%%%%%%%%%%%%%%%%%%%%%%%%%%%%%%%%%%%
\begin{equation} \label{eqn:hb} 
|h_{\rb}(n) | \le c \, M_6(a) \cdot n^{\rRe \, \gamma(n)} \cdot \Phim^{\, n}, 
\end{equation}
%%%%%%%%%%%%%%%%%%%%%%%%%%%%%%%%%%%%%%%%%%%%%%%%%%%%%%%%%%%%%%%%%%%%%%%%
where $M_6(a) := 2 M_4(a/2) \sqrt{\pi/a} + 2/a$ with $M_4(a)$ defined in Lemma 
$\ref{lem:Lipschitz2}$ and $a := \frac{1}{2} \phi''(x_0) > 0$. 
For any $0 < \ve < \ve_0$ there exists a constant $C_2(\ve) > 0$ such that  
%%%%%%%%%%%%%%%%%%%%%%%%% eqn:hc %%%%%%%%%%%%%%%%%%%%%%%%%%%%%%%%%%%%%%%%
\begin{equation} \label{eqn:hc}
| h_{\rc}(n) | \le C_2(\ve) \cdot n^{\rRe \, \gamma(n)} \cdot \Phim^{\, n}, 
\qquad {}^{\forall} n \ge N_0(\ve).   
\end{equation}
%%%%%%%%%%%%%%%%%%%%%%%%%%%%%%%%%%%%%%%%%%%%%%%%%%%%%%%%%%%%%%%%%%%%%%%%    
\end{lemma}
%%%%%%%%%%%%%%%%%%%%%%%%% begin proof %%%%%%%%%%%%%%%%%%%%%%%%%%%%%%%%%%%%
{\it Proof}. 
If $| y | \le 2 \ve_2$ $(\le \ve_1)$ then estimate 
\eqref{eqn:taylor-p} yields    
%%%%%%
\[
a \, y^2 + \eta( y ) 
\ge a \, y^2 - b \, |y|^3 = a \, y^2 
\left(1 - \ts \frac{b}{a} \, |y| \right) \\ 
\ge a \, y^2 \left( 1 - \ts \frac{2 b}{a}\, \ve_2 \right) 
\ge \ts \frac{a}{2} \, y^2,  
\]
%%%%%%
which together with estimate \eqref{eqn:taylor-u} and definition 
\eqref{eqn:Hb} gives 
%%%%%%%%%%%%%%%%%%%%%%%% eqn:Hb-e %%%%%%%%%%%%%%%%%%%%%%%%%%%%%%%%%%%%%%%
\begin{alignat}{2}
|H_{\rb}(x; n)| 
&\le c \cdot n^{\rRe \, \gamma(n)} \cdot \Phim^{\, n} \cdot 
e^{ - \frac{a}{2} n y^2} \, | y |, \qquad & 
|{}^{\forall} y \, | &\le 2 \ve_2, \nonumber \\[2mm]
&= c \cdot n^{\rRe \, \gamma(n)} \cdot \Phim^{\, n} \cdot 
\varphi_1(z; a/2) \cdot \ts \frac{1}{\sqrt{n}}, \qquad &  
|{}^{\forall} z \,| &\le 2 \ve_2 \sqrt{n},  \label{eqn:Hb-e}
\end{alignat} 
%%%%%%%%%%%%%%%%%%%%%%%%%%%%%%%%%%%%%%%%%%%%%%%%%%%%%%%%%%%%%%%%%%%%%%%%
where $\varphi_1(z; a) := |z| \, e^{- a z^2}$. 
If $0 < \ve < \ve_2$ and $n \ge N_1(\ve)$ then $[\xi_l, \, \xi_m] \subset 
[- 2 \ve_2 \sqrt{n}, \, 2 \ve_2 \sqrt{n}]$ follows from \eqref{eqn:range}, so 
estimate \eqref{eqn:Hb-e} is available for all $z \in [\xi_l, \, \xi_m]$, yielding  
%%%%%
\[
| h_{\rb}(n) | \le \sum_{k=l}^{m-1} | H_{\rb} (k/n; n)| 
\le c \cdot n^{\rRe \, \gamma(n)} \cdot \Phim^{\, n} \cdot R(\varphi_1; \vD),   
\]
%%%%%
where the Riemann sum $R(\varphi_1; \vD) := \sum_{k=l}^{m-1} 
\varphi_1(\xi_k; a/2) \cdot \frac{1}{\sqrt{n}}$ is estimated as  
%%%%%
\begin{align*}
R(\varphi_1; \vD) & \le 
\sum_{k=l}^{m-1} \int_{\xi_k}^{\xi_{k+1}} |\varphi_1(\xi_k; a/2) - \varphi_1(z; a/2) | \, d z 
+ \int_{\xi_l}^{\xi_m} \varphi_1(z; a/2) \, d z \\
& \le M_4(a/2) \sum_{k=l}^{m-1} \int_{\xi_k}^{\xi_{k+1}} | \xi_k - z | \, \varphi(z; a/4) \, d z 
+ \int_{-\infty}^{\infty} \varphi_1(z; a/2) \, d z \\
& \le \frac{M_4(a/2)}{\sqrt{n}} \int_{-\infty}^{\infty} \varphi(z; a/4) \, d z + 
\int_{-\infty}^{\infty} \varphi_1(z; a/2) \, d z = 
2 M_4(a/2) \sqrt{\frac{\pi}{a n}} + \frac{2}{a} \le M_6(a), 
\end{align*} 
%%%%%
where the second inequality is obtained by Lemma \ref{lem:Lipschitz2}. 
Now \eqref{eqn:hb} follows readily. 
%%%%%
\par
%%%%%
Since $\bal(n)$ and $\bbeta(n)$ are bounded by \eqref{eqn:bdd}, 
there exists a constant $C_1(\ve) > 0$ such that $|u(x; n)| \le C_1(\ve)$ for any 
$n \in \N$ and $x \in [x_0-\ve, \, x_0+\ve]$, which together with 
\eqref{eqn:C0(ve)} yields      
%%%%%
\[
|H_{\rc}(x; n)| \le C_1(\ve) \cdot n^{\rRe \, \gamma(n)} \cdot \Phim^{\, n} 
\cdot C_0(\ve)/n,  
\quad {}^{\forall} n \ge N_0(\ve), \,\, x_0-\ve \le {}^{\forall} x \le x_0+\ve.     
\]
%%%%%
Since $m-l = \lceil (x_0 + \ve) n\rceil - \lceil (x_0 - \ve) n\rceil 
\le (2 \ve + 1) n$, we have for any $n \ge N_0(\ve)$,    
%%%%% 
\[
| h_{\rc}(n) | \le \sum_{k=l}^{m-1} |H(k/n; n)| \le C_0(\ve) \cdot C_1(\ve) \cdot 
\frac{m-l}{n} \cdot n^{\rRe \, \gamma(n)} \cdot \Phim^{\, n} = 
C_2(\ve) \cdot n^{\rRe \, \gamma(n)} \cdot \Phim^{\, n},  
\]
%%%%%%
where $C_2(\ve) := (2 \ve + 1) \cdot C_0(\ve) \cdot C_1(\ve)$. 
This establishes estimate \eqref{eqn:hc}. \hfill $\Box$
%%%%%%%%%%%%%%%%%%%%%%%% end proof %%%%%%%%%%%%%%%%%%%%%%%%%%%%%%%%%%%%%%
%%%%%%%%%%%%%%%%%%%%%%%% prop:h(n) %%%%%%%%%%%%%%%%%%%%%%%%%%%%%%%%%%%%%%%
\begin{proposition} \label{prop:h(n)} 
For any $0 < \ve < \ve_2$, there is a constant $M(\ve) > 0$ such that 
%%%%%%%%%%%%%%%%%%%%%%%% eqn:h(n) %%%%%%%%%%%%%%%%%%%%%%%%%%%%%%%%%%%%%%%%
\begin{equation} \label{eqn:h(n)} 
h(n) = \sqrt{2 \pi} \, \frac{ u(x_0; n) }{\sqrt{ \phi''(x_0)} } \cdot 
n^{\gamma(n) + \frac{1}{2} } \cdot \Phim^{\, n} \cdot 
\left\{ 1 + \frac{e(n)}{ \sqrt{n} } \right \}, 
\qquad |e(n)| \le M(\ve), 
\end{equation}
%%%%%
for any $n \ge N_1(\ve)$, where $\ve_2$ and $N_1(\ve)$ are given in Lemma  
$\ref{lem:ha}$. 
\end{proposition}
%%%%%%%%%%%%%%%%%%%%%%%% begin proof %%%%%%%%%%%%%%%%%%%%%%%%%%%%%%%%%%%%%
{\it Proof}. This readily follows from $h(n) = h_{\ra}(n) + h_{\rb}(n) + h_{\rc}(n)$ and  
Lemmas \ref{lem:ha} and \ref{lem:hbc}. \hfill $\Box$
%%%%%%%%%%%%%%%%%%%%%%%% end proof %%%%%%%%%%%%%%%%%%%%%%%%%%%%%%%%%%%%%%
%%%%%%%%%%%%%%%%%%%%%%%% subsec:end %%%%%%%%%%%%%%%%%%%%%%%%%%%%%%%%%%%%%
\subsection{Irregular Sums and End Components} \label{subsec:end} 
%%%%%%%%%%%%%%%%%%%%%%%%%%%%%%%%%%%%%%%%%%%%%%%%%%%%%%%%%%%%%%%%%%%%%%%%
We shall estimate the left-end component $g_0(n)$ in \eqref{eqn:divide}. 
When $r_1 < +\infty$ the estimate for the right-end component $g_1(n)$ 
follows from the left-end counterpart by reflectional symmetry \eqref{eqn:refl}.   
If we make the translation $k \mapsto l := k - \lceil r_0 n \rceil$ 
for convenience, we can write   
%%%%%%%%%%%%%%%%%%%%%%%%%
\begin{alignat*}{3}
\sigma_i k + \lambda_i n + \alpha_i(n) &= \sigma_i \, l + 
\bar{\lambda}_i n + \bar{\alpha}_i(n),  \qquad & 
\bar{\lambda}_i &:= l_i(r_0) \qquad & (i & \in I),   \\ 
\tau_j k + \mu_j n + \beta_j(n) &= \tau_j \, l + 
\bar{\mu}_j n + \bar{\beta}_j(n), \qquad & 
\bar{\mu}_j &:= m_j(r_0) \qquad & (j & \in J),  
\end{alignat*} 
%%%%%%%%%%%%%%%%%%%%%%%%
where $\bar{\alpha}_i(n) :=\alpha_i(n) + \sigma_i (\lceil r_0 n \rceil - r_0 n) $ and  
$\bar{\beta}_j(n) := \beta_j(n) + \tau_j (\lceil r_0 n \rceil - r_0 n)$. 
Note that $\bar{\alpha}_i(n)$ here is the same as $\alpha_i^{(0)}(n)$ in 
formula \eqref{eqn:punc0}.  
Put $I_0^+ := I \setminus I_0$ and $J_0^+ := J \setminus J_0$, 
where the index sets $I_0$ and $J_0$ are defined in \eqref{eqn:IJ01}.   
Then $G(k; n)$ factors as 
%%%%%%%%%%%%%%%%%%%%%%%%%% eqn:G0 %%%%%%%%%%%%%%%%%%%%%%%%%%%%%%%%%%%%%%%
\begin{subequations} \label{eqn:G0}
\begin{gather} 
%%%%%%%%%%%%%%%%%%%%%%%%%% eqn:G=G0p %%%%%%%%%%%%%%%%%%%%%%%%%%%%%%%%%%%%
G(k; n) = G_0( l; n) \cdot G_0^+( l ; n),  \qquad  l := k - \lceil r_0 n \rceil, 
\label{eqn:G=G0p} \\[2mm]   
%%%%%%%%%%%%%%%%%%%%%%%%% eqn:G0p %%%%%%%%%%%%%%%%%%%%%%%%%%%%%%%%%%%%%%%  
G_0(l; n) := \dfrac{\prod_{i \in I_0}
\vG( \sigma_i \, l + \bar{\alpha}_i(n))}{\prod_{j \in J_0} 
\vG( \tau_j \, l + \bar{\beta}_j(n))}, \quad 
G_0^+(l; n) := \dfrac{\prod_{i \in I_0^+}
 \vG( \sigma_i \, l + \bar{\lambda}_i \, n + 
\bar{\alpha}_i(n))}{\prod_{j \in J_0^+} \vG( \tau_j \, l + 
\bar{\mu}_j \, n + \bar{\beta}_j(n))}.    \label{eqn:G0p}  
\end{gather}
\end{subequations}  
%%%%%%%%%%%%%%%%%%%%%%%%%%%%%%%%%%%%%%%%%%%%%%%%%%%%%%%%%%%%%%%%%%%%%%%%
From Lemma \ref{lem:I0I1} one has $\sigma_i > 0$ for $i \in I_0$ and $\tau_j > 0$ 
for $j \in J_0$, whereas condition \eqref{eqn:balanced} at $x = r_0$ implies that  
$( \bar{\lambda}_i )_{i \in I_0^+}$ and $( \bar{\mu}_j )_{j \in J_0^+}$ 
are balanced to the effect that     
%%%%%%%%%%%%%%%%%%%%%%%%%% eqn:balance-bq %%%%%%%%%%%%%%%%%%%%%%%%%%%%%%%
\begin{equation} \label{eqn:balance-bq} 
\sum_{i \in I_0^+} \bar{\lambda}_i = \sum_{j \in J_0^+} \bar{\mu}_j.    
\end{equation}
%%%%%%%%%%%%%%%%%%%%%%%%%%%%%%%%%%%%%%%%%%%%%%%%%%%%%%%%%%%%%%%%%%%%%%%% 
However, since $(\sigma_i)_{i \in I_0}$ and $(\tau_j)_{j \in J_0}$, resp. 
$(\sigma_i)_{i \in I_0^+}$ and $(\tau_j)_{j \in J_0^+}$, may not be balanced, we put  
%%%%%%%%%%%%%%%%%%%%%%%%%% eqn:rho0 %%%%%%%%%%%%%%%%%%%%%%%%%%%%%%%%%%%%%%
\begin{equation} \label{eqn:rho0}
\rho_0 := \sum_{i \in I_0} \sigma_i - \sum_{j \in J_0} \tau_j, \qquad 
\rho_0^+ := \sum_{i \in I_0^+} \sigma_i - \sum_{j \in J_0^+} \tau_j, 
\qquad \rho_0 = - \rho_0^+,  
\end{equation}
%%%%%%%%%%%%%%%%%%%%%%%%%%%%%%%%%%%%%%%%%%%%%%%%%%%%%%%%%%%%%%%%%%%%%%%%
where the relation $\rho_0 = - \rho_0^+$ follows from the first 
condition of \eqref{eqn:ap-bq2}. 
%%%%%
\par
%%%%%
We begin by giving an asymptotic behavior of $G_0(l; n)$ as $l \to \infty$ in terms of  
%%%%%%%%%%%%%%%%%%%%%%%%%% eqn:G0-asymp2 %%%%%%%%%%%%%%%%%%%%%%%%%%%%%%%%
\begin{align*} 
\Phi_0 &:= e^{-\rho_0} \prod_{i \in I_0} \sigma_i^{\sigma_i} 
\prod_{j \in J_0} \tau_j^{-\tau_j},  \\
u_0(n) &:= (2 \pi)^{\frac{|I_0|-|J_0|}{2}} 
\prod_{i \in I_0} \sigma_i^{\bar{\alpha}_i(n) - \frac{1}{2} } 
\prod_{j \in J_0} \tau_j^{ \frac{1}{2} - \bar{\beta}_j(n)},  \\  
\gamma_0(n) &:= \sum_{i \in I_0} \bar{\alpha}_i(n) - 
\sum_{j \in J_0} \bar{\beta}_j(n) + \frac{|J_0|-|I_0|}{2}. 
\end{align*} 
%%%%%%%%%%%%%%%%%%%%%%%%%%%%%%%%%%%%%%%%%%%%%%%%%%%%%%%%%%%%%%%%%%%%%%%%
Note that $\Phi_0$ is positive and $u_0(n)$ is nonzero due to the positivity of 
$\sigma_i$ and $\tau_j$ for $i \in I_0$ and $j \in J_0$. 
We use the following general fact about the gamma function.  
%%%%%%%%%%%%%%%%%%%%%%%%%% lem:dist %%%%%%%%%%%%%%%%%%%%%%%%%%%%%%%%%%%%%%
\begin{lemma} \label{lem:dist}
For any $z \in \C \setminus \Z_{\le 0}$ and any integer $m$ such that 
$m \ge 1+ |\rRe \, z|$, 
%%%%%
\[
|\vG(z)| \le \dfrac{2 |\vG(z+m)|}{\mathrm{dist}(z, \, \Z_{\le 0})}. 
\]
%%%%%
\end{lemma}
%%%%%%%%%%%%%%%%%%%%%%%%%% begin proof %%%%%%%%%%%%%%%%%%%%%%%%%%%%%%%%%%%
{\it Proof}. If $\rRe \, z > 0$ we have $\mathrm{dist}(z, \, \Z_{\le 0}) = |z|$ 
and the results follows readily.  
If $\rRe \, z \le 0$ then $\rRe (z+m) \ge 1$ and so the sequence 
$|z|, |z+1|, \cdots, |z+m-1|$ contains $\mathrm{dist}(z, \Z_{\le 0})$ as its 
minimum with the next smallest $\ge 1/2$ and all the rest $\ge 1$, 
thus $|(z; \, m)| = |z||z+1|\cdots |z+m-1| \ge \mathrm{dist}(z, \Z_{\le 0})/2$, 
hence $|\vG(z)| = |\vG(z+m)/(z; \, m)| \le 2|\vG(z+m)|/\mathrm{dist}(z, \Z_{\le 0})$.    
\hfill $\Box$
%%%%%%%%%%%%%%%%%%%%%%%%%% end proof %%%%%%%%%%%%%%%%%%%%%%%%%%%%%%%%%%%%
%%%%%%%%%%%%%%%%%%%%%%%%%% lem:G0 %%%%%%%%%%%%%%%%%%%%%%%%%%%%%%%%%%%%%%%
\begin{lemma} \label{lem:G0}
There exists a constant $K_0 > 0$ such that   
%%%%%%%%%%%%%%%%%%%%%%%%%% eqn:G0-asymp %%%%%%%%%%%%%%%%%%%%%%%%%%%%%%%%%
\begin{equation} \label{eqn:G0-asymp} 
|G_0(l; n)| \le  K_0 \cdot \delta_0(n)^{-1} \cdot (1+ l)^{|\rRe \, \gamma_0(n)|} \cdot 
l^{\rho_0 \, l} \cdot \Phi_0^{l}, \qquad 
{}^{\forall} l \in \Z_{\ge 0}, \,\, {}^{\forall} n \in \N,           
\end{equation}
%%%%%%%%%%%%%%%%%%%%%%%%%%%%%%%%%%%%%%%%%%%%%%%%%%%%%%%%%%%%%%%%%%%%%%%%
where for $l = 0$ the convention $l^{\rho_0 \, l} = 1$ is employed.     
\end{lemma}
%%%%%%%%%%%%%%%%%%%%%%%%%%%%%%%%%%%%%%%%%%%%%%%%%%%%%%%%%%%%%%%%%%%%%%%%
%%%%%%%%%%%%%%%%%%%%%%%%%% begin proof %%%%%%%%%%%%%%%%%%%%%%%%%%%%%%%%%%%
{\it Proof}. 
Note that $G_0(l; n)$ in \eqref{eqn:G0p} takes a finite value for every  
$l \ge \kappa := \max_{i \in I_0} (R+1)/\sigma_i$ and $n \in \N$,  
since \eqref{eqn:bdd} implies that $\sigma_i l + \rRe \, \bar{\alpha}_i(n) 
\ge \sigma_i l + \rRe \, \alpha_i(n) \ge \sigma_i l - R \ge 1$ for $i \in I_0$.   
By Stirling's formula \eqref{eqn:stirling} we have  
$G_0(l; n) = u_0(n) \cdot l^{\gamma_0(n) + \rho_0 \, l} 
\cdot \Phi_0^{l} \cdot \{1+ O(1/l) \}$ as $l \to +\infty$  
uniformly with respect to $n \in \N$. 
Thus there exists a constant $M_0 > 0$ such that 
%%%%%%%%
\[ 
|G_0(l; n)| \le  M_0 \cdot 
(1+ l)^{| \rRe \, \gamma_0(n) |} \cdot l^{\rho_0 \, l} \cdot \Phi_0^{l},  
\qquad {}^{\forall} l \ge \kappa, \,\, {}^{\forall} n \in \N.       
\]
%%%%%%%
\par
%%%%%%%
Take the smallest integer $m \ge \max_{i \in I_0} \{ 1 + \sigma_i (\kappa +1) + R \}$ 
and put 
%%%%%%%
\[
\bar{G}_0(l; n) := 
\dfrac{\prod_{i \in I_0} \vG(\sigma_i l + \bar{\alpha}_i(n) + m)}{ \prod_{j \in J_0} 
\vG(\tau_j l + \bar{\beta}_j(n)) }.  
\] 
%%%%%%%
Since $1+ |\sigma_i l + \rRe \, \bar{\alpha}_i(n)| \le 1 + \sigma_i l + 
|\rRe \, \bar{\alpha}_i(n)| \le 1 + \sigma_i l + 
|\rRe \, \alpha_i(n)| + \sigma_i \le m$ for any $0 \le l < \kappa$, 
$n \in \N$ and $i \in I_0$, Lemma \ref{lem:dist} implies that for any 
$0 \le l < \kappa$ and $n \in \N$,  
%%%%%%%
\[
|G_0(l; n)| \le \dfrac{ 2^{|I_0|} \cdot |\bar{G}_0(l; n)|}{ \prod_{i \in I_0} 
\mathrm{dist}(\sigma_i l + \bar{\alpha}_i(n), \, \Z_{\le 0})} 
\le \dfrac{ 2^{|I_0|} \cdot |\bar{G}_0(l; n)|}{ \prod_{i \in I_0} 
\mathrm{dist}(\bar{\alpha}_i(n), \, \Z_{\le 0} + |\sigma_i| \Z_{\le 0})} 
\le \dfrac{ 2^{|I_0|} \cdot |\bar{G}_0(l; n)| }{ \delta_0(n) }.     
\]
%%%%%%%
In view of condition \eqref{eqn:punc} there exists a 
constant $M_0' > 0$ such that 
%%%%%%%
\[ 
2^{|I_0|} \cdot | \bar{G}_0(l; n) | \le M_0' \cdot (1+ l)^{|\rRe \, \gamma_0(n)|} 
\cdot l^{\rho_0 l} \cdot \Phi_0^l,   
\qquad 0 \le {}^{\forall} l < \kappa, \,\,  {}^{\forall} n \in \N. 
\]
%%%%%%%%
Then by $1 \le \delta_0(n)^{-1}$ the estimate \eqref{eqn:G0-asymp} holds 
with the constant $K_0 := \max\{ M_0, \, M_0' \}$. \hfill $\Box$ \par\medskip 
%%%%%%%%%%%%%%%%%%%%%%%%%% end proof %%%%%%%%%%%%%%%%%%%%%%%%%%%%%%%%%%%%%
We proceed to the investigation into $G_0^+(l; n)$ by writing    
%%%%%%%%%%%%%%%%%%%%%%%%%% eqn:Hp %%%%%%%%%%%%%%%%%%%%%%%%%%%%%%%%%%%%%%%
\begin{equation} \label{eqn:Hp}
G_0^+(l; n) = H_0^+ \left( l/n ; n \right), \qquad 
H_0^+(x; n) := \dfrac{\prod_{i \in I_0^+}
 \vG( \bar{l}_i(x) \, n + \bar{\alpha}_i(n))}{\prod_{j \in J_0^+} 
\vG( \bar{m}_j(x) \, n + \bar{\beta}_j(n))},     
\end{equation}
%%%%%%%%%%%%%%%%%%%%%%%%%%%%%%%%%%%%%%%%%%%%%%%%%%%%%%%%%%%%%%%%%%%%%%%%
where $\bar{l}_i(x) := \sigma_i x + \bar{\lambda}_i$ and 
$\bar{m}_j(x) := \tau_j x + \bar{\mu}_j$, and then by putting    
%%%%%%%%%%%%%%%%%%%%%%%%%%%
\begin{align*} 
\Phi_0^+(x) &:= e^{- \rho_0^+ x} \prod_{i \in I_0^+} \bar{l}_i(x)^{\bar{l}_i(x)} 
\prod_{j \in J_0^+} \bar{m}_j(x)^{-\bar{m}_j(x)},  \\
u_0^+(x; n) &:= (2 \pi)^{\frac{|I_0^+|-|J_0^+|}{2}} 
\prod_{i \in I_0^+} \bar{l}_i(x)^{\bar{\alpha}_i(n) - \frac{1}{2} } 
\prod_{j \in J_0^+} \bar{m}_j(x)^{ \frac{1}{2} - \bar{\beta}_j(n)},  \\
\gamma_0^+(n) &:= \sum_{i \in I_0^+} \bar{\alpha}_i(n) - 
\sum_{j \in J_0^+} \bar{\beta}_j(n) + \frac{|J_0^+|-|I_0^+|}{2}.   
\end{align*}
%%%%%%%%%%%%%%%%%%%%%%%%%%
Note that $\Phi_0^+(x)$ and $u_0^+(x; n)$ are well-defined continuous 
functions on $[0, \, \ve]$ with $\Phi_0^+(x)$ being positive while 
$u_0^+(x; n)$ non-vanishing and uniformly bounded in $n \in \N$.  
%%%%%%%%%%%%%%%%%%%%%%%%%% lem:Hp %%%%%%%%%%%%%%%%%%%%%%%%%%%%%%%%%%%%%%%
\begin{lemma} \label{lem:Hp} 
For any $0 < \ve < \ve_0$ there exist an integer $N_0(\ve) \in \N$ and 
a constant $K_0^+(\ve) > 0$ such that for any 
$n \ge N_0(\ve)$ and $0 \le x \le \ve$,      
%%%%%%%%%%%%%%%%%%%%%%%% eqn:Hp-asymp %%%%%%%%%%%%%%%%%%%%%%%%%%%%%%%%%%%
\begin{equation} \label{eqn:Hp-asymp} 
|H_0^+(x; n)| \le K_0^+(\ve) \cdot n^{ \rRe \, \gamma_0^+(n) + \rho_0^+ x \, n} 
\cdot \Psi_0^+(\ve)^n \qquad \mbox{with} \quad   
\Psi_0^+(\ve) := \max_{0 \le x \le \ve} \Phi_0^+(x).   
\end{equation} 
%%%%%%%%%%%%%%%%%%%%%%%%%%%%%%%%%%%%%%%%%%%%%%%%%%%%%%%%%%%%%%%%%%%%%%%%
%%%%%%%%%%%%%%%%%%%%%%%%%%%%%%%%%%%%%%%%%%%%%%%%%%%%%%%%%%%%%%%%%%%%%%%%     
\end{lemma}
%%%%%%%%%%%%%%%%%%%%%%%%%%%%%%%%%%%%%%%%%%%%%%%%%%%%%%%%%%%%%%%%%%%%%%%%
%%%%%%%%%%%%%%%%%%%%%%% begin proof %%%%%%%%%%%%%%%%%%%%%%%%%%%%%%%%%%%%%%
{\it Proof}. 
From the definitions of $I_0^+$, $J_0^+$, $\bar{l}_i(x)$, $\bar{m}_j(n)$,  
there is a constant $c(\ve) > 0$ such that     
%%%%%%%%%%%%%%
\[
\bar{l}_i(x) > c(\ve) \quad (i \in I_0^+), \qquad  
\bar{m}_j(x) > c(\ve) \quad (j \in J_0^+), \qquad  0 \le {}^{\forall} x \le \ve. 
\]
%%%%%%%%%%%%%%
By condition \eqref{eqn:bdd}, $H_0^+(x; n)$ takes a finite value for any 
$x \in [0, \, \ve]$ and $n \ge N_0(\ve) := (R+1)/c(\ve)$ and Stirling's 
formula \eqref{eqn:stirling} implies that $H_0^+(x; n)$ admits an 
asymptotic formula 
%%%%%%%%%%%%%%%%%%%%%%%%
\[ 
H_0^+(x; n) = u_0^+(x; n) \cdot n^{\gamma_0^+(n) + \rho_0^+ x \, n} \cdot 
\Phi_0^+(x)^n \cdot \{ 1+ O(1/n) \} \qquad \mbox{as} \quad n \to + \infty,     
\]
%%%%%%%%%%%%%%%%%%%%%%%%
uniform in $x \in [0, \, \ve]$, where one also uses the equality 
$\sum_{i \in I_0^+} \bar{l}_i(x) - 
\sum_{j \in J_0^+} \bar{m}_j(x) = \rho_0^+ \, x$,  
which is due to balancedness condition \eqref{eqn:balance-bq} 
and definition \eqref{eqn:rho0}.  
From this estimate and the boundedness of $u_0^+(x; n)$ coming from 
\eqref{eqn:bdd} the assertion \eqref{eqn:Hp-asymp} follows readily.  
\hfill $\Box$ \par\medskip 
%%%%%%%%%%%%%%%%%%%%%%% end proof %%%%%%%%%%%%%%%%%%%%%%%%%%%%%%%%%%%%%%%
Now we are able to give an estimate for the left-end component $g_0(n)$ in terms of   
%%%%%%%%%%%%%%%%%%%%%%% eqn:del-Psi0 %%%%%%%%%%%%%%%%%%%%%%%%%%%%%%%%%%%%
\begin{subequations} \label{eqn:del-Psi0}
\begin{align}
\ve_3 &:= 
\begin{cases}
+ \infty & (\mbox{if $\rho_0 \ge 0$}), \\[1mm]
e^{-1} \Phi_0^{-1/\rho_0} & (\mbox{if $\rho_0 < 0$}),  
\end{cases}
\label{eqn:del0} \\[2mm] 
\Psi_0(\ve) &:= 
\begin{cases}
1 & (\mbox{if $\rho_0 > 0$, or $\rho_0 = 0$ with $\Phi_0 \le 1$}), \\[1mm]
(\ve^{\rho_0} \Phi_0)^{\ve} \,\,\,\, & 
(\mbox{if $\rho_0 < 0$, or $\rho_0 = 0$ with $\Phi_0 > 1$}),  
\end{cases}
\label{eqn:Psi0} \\[2mm] 
\vD_0(\ve) &:= \Psi_0(\ve) \cdot \Psi_0^+(\ve). \label{eqn:vD0}
\end{align} 
\end{subequations}
%%%%%%%%%%%%%%%%%%%%%%%%%%%%%%%%%%%%%%%%%%%%%%%%%%%%%%%%%%%%%%%%%%%%%%%%
%%%%%%%%%%%%%%%%%%%%%%% lem:g0-bound %%%%%%%%%%%%%%%%%%%%%%%%%%%%%%%%%%%%
\begin{lemma} \label{lem:g0-bound}  
For any $0 < \ve < \ve_4 := \min\{ \ve_0, \, \ve_2, \, \ve_3 \}$ with $\ve_0$, 
$\ve_2$ and $\ve_3$ defined in \eqref{eqn:ve}, Lemma $\ref{lem:ha}$  
and \eqref{eqn:del0} respectively, there exist $N_0(\ve) \in \N$ and 
$K_0(\ve) > 0$ such that 
%%%%%%%%%%%%%%%%%%%%%%% eqn:g0-bound %%%%%%%%%%%%%%%%%%%%%%%%%%%%%%%%%%%%
\begin{equation*} \label{eqn:g0-bound}
|g_0(n)| \le K_0(\ve) \cdot \delta_0(n)^{-1} \cdot 
n^{ |\rRe \, \gamma_0(n)| + \rRe \, \gamma_0^+(n) +1 } \cdot 
\vD_0(\ve)^n, \qquad {}^{\forall} n \ge N_0(\ve).    
\end{equation*}
%%%%%%%%%%%%%%%%%%%%%%%%%%%%%%%%%%%%%%%%%%%%%%%%%%%%%%%%%%%%%%%%%%%%%%%%    
\end{lemma} 
%%%%%%%%%%%%%%%%%%%%%%%%%%%%%%%%%%%%%%%%%%%%%%%%%%%%%%%%%%%%%%%%%%%%%%%%
%%%%%%%%%%%%%%%%%%%%%%%%%% begin proof %%%%%%%%%%%%%%%%%%%%%%%%%%%%%%%%%%
{\it Proof}. 
It follows from formulas \eqref{eqn:Hp}, \eqref{eqn:Hp-asymp} and 
\eqref{eqn:rho0} that   
%%%%%%
\[
|G_0^+(l; n)| = |H_0^+( l/n; n ) | \le K_0^+(\ve) \cdot 
n^{\rRe \, \gamma_0^+(n)} \cdot n^{- \rho_0 l} \cdot \Psi_0^+(\ve)^n.    
\]
%%%%%%   
Multiplying this estimate by inequality \eqref{eqn:G0-asymp}, we have 
from formula \eqref{eqn:G=G0p},   
%%%%%%%%%%%%%%%%%%%%%%%% eqn:GG0 %%%%%%%%%%%%%%%%%%%%%%%%%%%%%%%%%%%%% 
\begin{equation} \label{eqn:GG0}
|G(k; n)| \le K(\ve) \cdot \delta_0(n)^{-1} \cdot n^{\rRe \, \gamma_0^+(n)} 
\cdot \varphi(l; n) \cdot \Psi_0^+(\ve)^n \cdot 
(1+l)^{|\rRe \, \gamma_0(n)|},    
\end{equation}
%%%%%%%%%%%%%%%%%%%%%%%%%%%%%%%%%%%%%%%%%%%%%%%%%%%%%%%%%%%%%%%%%%%%%%
for any $n \ge N_0(\ve)$ and $0 \le l := k - \lceil r_0 n \rceil < \ve n$, 
where $K(\ve) := K_0 \cdot K_0^+(\ve)$ and 
%%%%%%%%%
\[
\varphi( t; n) := (t/n)^{\rho_0 t} \cdot \Phi_0^t \quad (t > 0) 
\quad \mbox{with} \quad \varphi(0; n) = \lim_{t \to +0} 
\varphi(t; n) = 1. 
\]  
%%%%%%%%%
\par
%%%%% 
A bit of differential calculus shows the following: 
\begin{enumerate} 
\item[(i)] If either $\rho_0 > 0$ or $\rho_0 = 0$ with $\Phi_0 \le 1$, then 
$\varphi(t; n)$ is non-increasing in $t \ge 0$ and hence  
$\varphi(t; n) \le \varphi(0; n) = 1 = \Psi_0(\ve)^n$ for any $t \ge 0$. 
\item[(ii)] If either $\rho_0 < 0$ or $\rho_0 = 0$ with $\Phi_0 > 1$, 
then $\frac{d}{d t} \varphi (t; n) \ge 0$ in $0 \le t \le \ve_3 \, n$ with 
equality only when $t = \ve_3 \, n$, so that $\varphi(t; n) \le 
\varphi(\ve n; n) = (\ve^{\rho_0} \Phi_0)^{\ve n} = \Psi_0(\ve)^n$ 
for any $0 \le t \le \ve n$ $(< \ve_3 \, n)$, where $\ve_3$ and 
$\Psi_0(\ve)$ are defined in \eqref{eqn:del0} and \eqref{eqn:Psi0} respectively. 
\end{enumerate} 
In either case $0 < \varphi(t; n) \le \Psi_0(\ve)^n$ for any 
$0 \le t < \ve n$ and thus \eqref{eqn:GG0} and \eqref{eqn:vD0} lead to 
%%%%%%%%%%%%%%%%%%%%%%%%%% eqn:GGG0 %%%%%%%%%%%%%%%%%%%%%%%%%%%%%%%%%%%%
\begin{equation} \label{eqn:GGG0}
|G(k; n)| \le K(\ve) \cdot \delta_0(n)^{-1}\cdot 
n^{\rRe \, \gamma_0^+(n)} \cdot 
\vD_0(\ve)^n \cdot (1+l)^{|\rRe \, \gamma_0(n)|},    
\end{equation}
%%%%%%%%%%%%%%%%%%%%%%%%%%%%%%%%%%%%%%%%%%%%%%%%%%%%%%%%%%%%%%%%%%%%%%%
for any $n \ge N_0(\ve)$ and $0 \le l : = k - \lceil r_0 n \rceil < \ve n$. 
Since 
%%%%%%%
\[
\sum_{0 \le l < \ve n} (1+ l)^{|\rRe \, \gamma_0(n)|} \le 
\int_0^{\ve n+1} (1+t)^{|\rRe \, \gamma_0(n)|} \, d t  
\le \frac{(\ve n + 2)^{|\rRe \, \gamma_0(n)|+1}}{|\rRe \, \gamma_0(n)|+1} \le 
\{(2 + \ve ) \cdot n\}^{|\rRe \, \gamma_0(n)|+1},  
\]
%%%%%%%
summing up \eqref{eqn:GGG0} over the integers $0 \le l \le 
\lceil (r_0+\ve) n \rceil - \lceil r_0 n \rceil -1$ $(< \ve n)$ yields   
%%%%%%%
\[
|g_0(n)| \le K(\ve) \cdot \delta_0(n)^{-1} \cdot 
(2 + \ve)^{|\rRe \, \gamma_0(n)|+1} \cdot 
n^{ |\rRe \, \gamma_0(n)| + \rRe \, \gamma_0^+(n) +1} \cdot \vD_0(\ve)^n,  
\quad {}^{\forall} n \ge N_0(\ve).         
\]
%%%%%%%
Since $\gamma_0(n)$ is bounded by condition \eqref{eqn:bdd} we can 
take a constant $K_0(\ve) \ge K(\ve) \cdot (2 + \ve)^{|\rRe \, \gamma_0(n)|+1}$ 
to establish the lemma. \hfill $\Box$ 
%%%%%%%%%%%%%%%%%%%%%%%%%% end proof %%%%%%%%%%%%%%%%%%%%%%%%%%%%%%%%%%%
%%%%%%%%%%%%%%%%%%%%%%%%%% prop:end %%%%%%%%%%%%%%%%%%%%%%%%%%%%%%%%%%%%
\begin{proposition} \label{prop:end} 
For any $d > \Phi(r_0)$ there exists a positive constant $\ve_5 \le \ve_4$   
such that 
%%%%%%%%%%%%%%%%%%%%%%%%%%%% eqn:l-end %%%%%%%%%%%%%%%%%%%%%%%%%%%%%%%
\begin{equation} \label{eqn:l-end}
|g_0(n)| \le M_0(d, \ve) \cdot \delta_0(n)^{-1} \cdot d^n, \qquad 
{}^{\forall} n \ge N_0(\ve), \,\, 0 < {}^{\forall} \ve \le \ve_5, 
\end{equation}
%%%%%%%%%%%%%%%%%%%%%%%%%%%%%%%%%%%%%%%%%%%%%%%%%%%%%%%%%%%%%%%%%%%%
for some $M_0(d, \ve) > 0$ and $N_0(\ve) \in \N$ independent of $d$, where 
$\Phi(x)$ is defined in \eqref{eqn:Phi} and $\ve_4$ is given in Lemma 
$\ref{lem:g0-bound}$. 
When $r_1 < + \infty$, a similar statement can be made for the right-end 
component $g_1(n)$ in \eqref{eqn:divide}; for any $d > \Phi(r_1)$ there exists 
a sufficiently small $\ve_6 > 0$ such that  
%%%%%%%
\[
|g_1(n)| \le M_1(d, \ve) \cdot \delta_0(n)^{-1} \cdot d^n, \qquad 
{}^{\forall} n \ge N_1(\ve), \,\, 0 < {}^{\forall} \ve \le \ve_6. 
\]
%%%%%%%
\end{proposition} 
%%%%%%%%%%%%%%%%%%%%%%%%%%%% begin proof %%%%%%%%%%%%%%%%%%%%%%%%%%%%%%%
{\it Proof}. We show the assertion for the left-end component $g_0(n)$ only  
as the right-end counterpart follows by reflectional symmetry \eqref{eqn:refl}. 
Observe that $\Psi_0(\ve) \to 1$, $\Psi_0^+(\ve) \to \Phi(r_0)$ and 
so $\vD_0(\ve) \to \Phi(r_0)$ as $\ve \to +0$. 
Thus given $d > \Phi(r_0)$ there is a constant $0 < \ve_5 < \ve_4$ 
such that $d > \vD_0(\ve)$ for any $0 < \ve \le \ve_5$. 
Then Lemma \ref{lem:g0-bound} enables us to take a constant 
$M_0(d, \ve)$ as in \eqref{eqn:l-end}.    
\hfill $\Box$ \par\medskip\noindent
%%%%%%%%%%%%%%%%%%%%%%%%%%%% end proof %%%%%%%%%%%%%%%%%%%%%%%%%%%%%%%%%
%%%%%%%%%%%%%%%%%%%%%%%%%%% begin proof of thm:dspm %%%%%%%%%%%%%%%%%%%%%
{\it Proof of Theorem $\ref{thm:dspm}$}.  
As is mentioned at the end of \S \ref{subsec:results} only the singleton 
case $\Mm = \{ x_0 \}$ is treated for the sake of simplicity. 
We can take a number $d$ so that $\max\{\Phi(r_0), \, \Phi(r_1) \} < d < 
\Phim$, since $\Phi(x)$ attains its maximum only at the interior  
point $x_0 \in (r_0, \, r_1)$.   
For this $d$ take the numbers $\ve_5$ and $\ve_6$ as in  
Proposition \ref{prop:end} and put $\ve := \min \{ \ve_5, \, \ve_6 \}$. 
For this $\ve$ consider the numbers $\Phi_0^{\ve}$ and $\Phi_1^{\ve}$ 
in Lemma \ref{lem:side}, both of which are strictly smaller than $\Phim$. 
Take a number $d_0$ so that $\max\{ d, \, \Phi_0^{\ve}, \, \Phi_1^{\ve} \}  
< d_0 < \Phim$ and put $\lambda := \Phim/d_0 > 1$.  
Then the estimates in Propositions \ref{prop:h(n)} and \ref{prop:end} and 
Lemma \ref{lem:side} are put together into equation \eqref{eqn:divide} to 
yield    
%%%%%%%%%%%%%%%
\begin{equation*}
g(n) = C(n) \cdot  n^{\gamma(n) + \frac{1}{2} } \cdot 
\Phim^{\, n} \cdot \left\{ 1 + \Omega(n) \right\},     
\end{equation*} 
%%%%%%%%%%%%%%%
where $C(n)$ is defined in \eqref{eqn:PhimC} and $\Omega(n)$ admits 
the estimate \eqref{eqn:uniform}.  \hfill $\Box$ \par \medskip
%%%%%%%%%%%%%%%%%%%%%%%%%%%% end proof %%%%%%%%%%%%%%%%%%%%%%%%%%%%%%%%
Even without assuming properness \eqref{eqn:max} we have 
the following convenient proposition. 
%%%%%%%%%%%%%%%%%%%%%%%%%%%%%% prop:irregular %%%%%%%%%%%%%%%%%%%%%%%%%%%% 
\begin{proposition} \label{prop:irregular}
Suppose that the sum $g(n)$ in \eqref{eqn:g(n)} satisfies 
balancedness \eqref{eqn:ap-bq2}, boundedness \eqref{eqn:bdd},  
admissibility \eqref{eqn:ap-bq1} and genericness \eqref{eqn:punc} 
along with continuity at infinity \eqref{eqn:infty} and 
convergence \eqref{eqn:conv} when $r_1 = + \infty$.  
For any $d > \Phim$ there exist $K > 0$ and $N \in \N$ such that 
%%%%%
\[
|g(n)| \le K \cdot d^n \cdot  \{ \delta_0(n)^{-1} + \delta_1(n)^{-1}  \},  
\qquad {}^{\forall} n \ge N.   
\]
%%%%%   
\end{proposition}
%%%%%%%%%%%%%%%%%%%%%%%%%%%%%% begin proof %%%%%%%%%%%%%%%%%%%%%%%%%%%%%
{\it Proof}. 
Divide $g(n)$ into three components; sums over $[r_0, \, r_0+ \ve]$, 
$[r_0+\ve, \, r_1-\ve]$ and $[r_1-\ve, \, r_1]$. 
Take $\ve > 0$ sufficiently small depending on how $d$ is close to $\Phim$. 
Apply Proposition \ref{prop:end} to the left and right components 
and then use Lemma \ref{lem:side} in the middle one. \hfill $\Box$ 
%%%%%%%%%%%%%%%%%%%%%%%%%%%%%% end proof %%%%%%%%%%%%%%%%%%%%%%%%%%%%%%   
%%%%%%%%%%%%%%%%%%%%%%%%%%%%%% subsec:apriori %%%%%%%%%%%%%%%%%%%%%%%%%%%
\subsection{A Priori Estimates} \label{subsec:apriori} 
%%%%%%%%%%%%%%%%%%%%%%%%%%%%%%%%%%%%%%%%%%%%%%%%%%%%%%%%%%%%%%%%%%%%%%%
We present the a priori estimates used in \S \ref{subsec:t-s}. 
In what follows we often use the inequality   
%%%%%%%%%%%%%%%%%%%%%%%%% eqn:e^x %%%%%%%%%%%%%%%%%%%%%%%%%%%%%%%%%%%%%%%
\begin{equation} \label{eqn:e^x}
|e^x -1 | \le |x| \, e^{|x|} \qquad (x \in \R). 
\end{equation}
%%%%%%%%%%%%%%%%%%%%%%%%%%%%%%%%%%%%%%%%%%%%%%%%%%%%%%%%%%%%%%%%%%%%%%%%
\par
%%%%%%
Given a positive constant $a$, we consider the function  
$\varphi(x; a) := e^{-a x^2}$.  
%%%%%%%%%%%%%%%%%%%%%%%% lem:Lipschitz %%%%%%%%%%%%%%%%%%%%%%%%%%%%%%%%%%
\begin{lemma} \label{lem:Lipschitz} 
If $x$, $y \in \R$ and $| y-x | \le 1$, then   
%%%%%%%%%%%%%%%%%%%%%%%% eqn:Lipschitz %%%%%%%%%%%%%%%%%%%%%%%%%%%%%%%%%%
\begin{equation} \label{eqn:Lipschitz}
| \varphi( y; a) - \varphi( x; a ) | \le M_1( a ) \, 
|y-x| \, \varphi( x; a/2 ),     
\end{equation}
%%%%%%%%%%%%%%%%%%%%%%%%%%%%%%%%%%%%%%%%%%%%%%%%%%%%%%%%%%%%%%%%%%%%%%%%
where $M_1(a) := a \, \ds \sup_{x \in \R} (2 |x| +1) 
e^{-\frac{a}{2} (x^2 - 4|x|-2)} < \infty$. 
\end{lemma}
%%%%%%%%%%%%%%%%%%%%%%%%%%%%%%%%%%%%%%%%%%%%%%%%%%%%%%%%%%%%%%%%%%%%%%%
%%%%%%%%%%%%%%%%%%%%%%%%% begin proof %%%%%%%%%%%%%%%%%%%%%%%%%%%%%%%%%%%
{\it Proof}. 
Put $h := y-x$. 
It then follows from inequality \eqref{eqn:e^x} that  
%%%%
\[
\begin{split}
\left| e^{a x^2 - a (x+h)^2} -1 \right| 
&= \left| e^{-a h (2 x + h)} -1 \right| \le a \, 
|h||2 x + h| \, e^{a |h||2 x + h|} \\
&\le a \, |h| \left( 2 |x| + |h| \right) \, 
e^{\alpha |h|(2|x| + |h|)} \le a \, |h| 
\left( 2 |x| + 1 \right) \, e^{a (2|x| + 1) },   
\end{split}
\]
%%%%
whenever $|h| \le 1$.  
Dividing both sides by $e^{a x^2}$ we have 
%%%%
\[
\begin{split}
\left| e^{- a (x+h)^2} - e^{- a x^2} \right|
&\le a \, |h| \left( 2 |x| + 1 \right) \, e^{-a (x^2 -2|x| - 1) } \\ 
&= a \, \left( 2 |x| + 1 \right) \, e^{-\frac{a}{2} (x^2 -4|x| - 2) } 
\cdot e^{- \frac{a}{2} x^2} \, |h| \le M_1(a) \, 
e^{- \frac{a}{2} x^2} \, |h|,      
\end{split}
\]
%%%%
which proves the lemma. \hfill $\Box$ \par\medskip
%%%%%%%%%%%%%%%%%%%%%%%% end proof %%%%%%%%%%%%%%%%%%%%%%%%%%%%%%%%%%%%%%
Let $b >0$, $m \ge 1$, $0 < \ve \le \frac{a}{4 b}$, and 
suppose that a function $\delta (x)$ admits an estimate  
%%%%%%%%%%%%%%%%%%%%%%%% eqn:pert-a %%%%%%%%%%%%%%%%%%%%%%%%%%%%%%%%%%%%%
\begin{equation} \label{eqn:pert-a}
| \delta( x ) | \le \frac{b}{m} \, |x|^3 \qquad 
(|{}^{\forall}x \, | \le \ve m).  
\end{equation}
%%%%%%%%%%%%%%%%%%%%%%%%%%%%%%%%%%%%%%%%%%%%%%%%%%%%%%%%%%%%%%%%%%%%%%%%
%%%%%%%%%%%%%%%%%%%%%%%% lem:pert %%%%%%%%%%%%%%%%%%%%%%%%%%%%%%%%%%%%%%% 
\begin{lemma} \label{lem:pert}
Under condition \eqref{eqn:pert-a}, the function $\psi(x; a) 
:= e^{- a x^2 + \delta(x)}$ satisfies  
%%%%%%%%%%%%%%%%%%%%%%%% eqn:pert-r %%%%%%%%%%%%%%%%%%%%%%%%%%%%%%%%%%%%%
\begin{equation} \label{eqn:pert-r}
| \psi(x; a) - \varphi(x; a) | \le 
\frac{b \, M_2(a)}{m} \, \varphi(x; a/2 ) 
\qquad (|{}^{\forall} x \,| \le \ve m),  
\end{equation}
%%%%%%%%%%%%%%%%%%%%%%%%%%%%%%%%%%%%%%%%%%%%%%%%%%%%%%%%%%%%%%%%%%%%%%%%% 
where $M_2(a) := \ds \sup_{x \in \R} |x|^3 \, 
e^{- \frac{a}{4} \, x^2} < \infty$. 
\end{lemma}
%%%%%%%%%%%%%%%%%%%%%%%%%%%%%%%%%%%%%%%%%%%%%%%%%%%%%%%%%%%%%%%%%%%%%%%%%
%%%%%%%%%%%%%%%%%%%%%%%% begin proof %%%%%%%%%%%%%%%%%%%%%%%%%%%%%%%%%%%%
{\it Proof}. 
For $|x| \le \ve m$, we have 
%%%%%
\begin{alignat*}{2}
\left| e^{- a x^2 + \delta(x) } - e^{-a x^2} \right| 
& = e^{- a x^2} \left| e^{\delta(x)} -1 \right| 
\le | \delta(x) | \, e^{- a x^2 + | \delta(x) |} 
\qquad & & \mbox{by \eqref{eqn:e^x},} \\
& \le \frac{b}{m} |x|^3 \, e^{-a x^2 + \frac{b}{m} |x|^3 }  
= \frac{b}{m} |x|^3 \, e^{-a x^2 
\left( 1-\frac{b}{a m} |x| \right) } 
\qquad & & \mbox{by \eqref{eqn:pert-a},} \\ 
& \le \frac{b}{m} |x|^3 \, 
e^{-a x^2 \left( 1-\frac{b \ve}{a} \right) } 
\qquad & & \mbox{by $|x| \le \ve m$,} \\
& \le \frac{b}{m} |x|^3 \, e^{-\frac{3 a}{4} x^2 } 
= \frac{b}{m} |x|^3 \, e^{-\frac{a}{4} x^2 } \cdot 
e^{-\frac{a}{2} x^2 } 
\qquad & & \mbox{by $0 < \ve \le \frac{ a }{4 b }$,} \\
& \le \frac{b M_2(a)}{m} \, e^{-\frac{a}{2} x^2 },  
\end{alignat*}
%%%%%
where the last inequality is by the definition of 
$M_2(a)$. \hfill $\Box$ 
%%%%%%%%%%%%%%%%%%%%%%%% end proof %%%%%%%%%%%%%%%%%%%%%%%%%%%%%%%%%%%%%%
%%%%%%%%%%%%%%%%%%%%%%%% lem:approxi %%%%%%%%%%%%%%%%%%%%%%%%%%%%%%%%%%%%
\begin{lemma} \label{lem:approxi} 
Under condition \eqref{eqn:pert-a}, if $|x| \le \ve m$, $|y| \le \ve m$ 
and $|y-x| \le 1/m$, then    
%%%%%%%%%%%%%%%%%%%%%%%% eqn:approxi %%%%%%%%%%%%%%%%%%%%%%%%%%%%%%%%%%%%
\begin{equation} \label{eqn:approxi}
| \psi(y; a) - \varphi(x; a) | \le \frac{M_3(a, b)}{m} \, \varphi(x; a/4 ).    
\end{equation}
%%%%%%%%%%%%%%%%%%%%%%%%%%%%%%%%%%%%%%%%%%%%%%%%%%%%%%%%%%%%%%%%%%%%%%%% 
where $M_3(a, b) := M_1(a) + b M_2(a) + b M_1( a/2 ) M_2( a )$.
\end{lemma}
%%%%%%%%%%%%%%%%%%%%%%%%%%%%%%%%%%%%%%%%%%%%%%%%%%%%%%%%%%%%%%%%%%%%%%%%%
%%%%%%%%%%%%%%%%%%%%%%%% begin proof %%%%%%%%%%%%%%%%%%%%%%%%%%%%%%%%%%%%
{\it Proof}.  Putting $y = x + h$ with $|h| \le 1/m$, we have   
%%%%%%%%%%%%%%%%%%%
\begin{alignat*}{2}
| \psi & ( y; a ) - \varphi( x; a ) | & 
&  \\
& \le | \psi( y; a ) - \varphi( y; a ) | + 
| \varphi( y; a ) - \varphi( x; a ) | 
\qquad & & \mbox{by t.i.,} \\[1mm]
& \le \ts \frac{b M_2(a)}{m} \, \varphi( y; \frac{a}{2} ) 
+ M_1(a) \, |h| \, \varphi( x; \frac{a}{2} ) 
\qquad & & \mbox{by \eqref{eqn:pert-r} and \eqref{eqn:Lipschitz},} \\[1mm]
& \le \ts \frac{b M_2(a)}{m}  
\left\{ |\varphi( y; \frac{a}{2} ) - \varphi( x; \frac{a}{2} ) | + 
\varphi( x; \frac{a}{2} ) \right\} 
+ \ts \frac{M_1(a)}{m} \, \varphi( x; \frac{a}{2} ) 
\qquad & & \mbox{by t.i. and $|h| \le \frac{1}{m}$,} \\[1mm]
& \le \ts \frac{b M_2(a)}{m} \left\{ M_1( \frac{a}{2} ) \, |h| \,   
\varphi( x; \frac{a}{4} ) + \varphi( x; \frac{a}{2} ) \right\} 
+ \ts \frac{M_1(a)}{m} \, \varphi( x; \frac{a}{2} )  
\qquad & & \mbox{by \eqref{eqn:Lipschitz},} \\[1mm]
& \le \ts \frac{b M_2(a)}{m} \left\{ M_1( \frac{a}{2} ) \,    
\varphi( x; \frac{a}{4} ) + \varphi( x; \frac{a}{2} ) \right\} 
+ \ts \frac{M_1(a)}{m} \, \varphi( x; \frac{a}{2} )  
\qquad & & \mbox{by $|h| \le \frac{1}{m} \le 1$,} \\[1mm]
& \le \ts \frac{M_3(a, b )}{m} \, \varphi( x; \frac{a}{4} ) 
\qquad & & \mbox{by $\varphi(x; \frac{a}{2}) \le \varphi(x; \frac{a}{4})$,} 
\end{alignat*}
%%%%%%%%%%%%%%
where t.i. refers to trigonometric inequality. \hfill $\Box$ 
%%%%%%%%%%%%%%%%%%%%%%%% end proof %%%%%%%%%%%%%%%%%%%%%%%%%%%%%%%%%%%%%%
%%%%%%%%%%%%%%%%%%%%%%%% lem:Lipschitz2 %%%%%%%%%%%%%%%%%%%%%%%%%%%%%%%%%
\begin{lemma} \label{lem:Lipschitz2} 
If $x$, $y \in \R$ and $| y-x | \le 1$, then 
$\varphi_1(x; a) := |x| \, e^{-a x^2}$ satisfies    
%%%%%%%%%%%%%%%%%%%%%%%% eqn:Lipschitz2 %%%%%%%%%%%%%%%%%%%%%%%%%%%%%%%%%
\begin{equation} \label{eqn:Lipschitz2}
| \varphi_1( y; a) - \varphi_1( x; a ) | \le M_4( a ) \, 
|y-x| \, \varphi( x; a/4 ),     
\end{equation}
%%%%%%%%%%%%%%%%%%%%%%%%%%%%%%%%%%%%%%%%%%%%%%%%%%%%%%%%%%%%%%%%%%%%%%%%%
where $M_4(a) := 1+ M_1(a) \cdot \ds \max_{x \in \R} (|x| +1) 
e^{-\frac{a}{4} x^2} < \infty$. 
\end{lemma}
%%%%%%%%%%%%%%%%%%%%%%%%%%%%%%%%%%%%%%%%%%%%%%%%%%%%%%%%%%%%%%%%%%%%%%%%%
%%%%%%%%%%%%%%%%%%%%%%%%% begin proof %%%%%%%%%%%%%%%%%%%%%%%%%%%%%%%%%%%
{\it Proof}. Putting $y = x + h$ with $|h| < 1$, one has  
%%%%%%%%%%%%%
\begin{alignat*}{2}
| \varphi_1 
& ( x+h; a ) - \varphi_1( x; a ) | = ||x+h| \, \varphi(x+h; a) - |x| \, \varphi(x; a)| &  &  \\
& \le |x+h| |\varphi(x+h; a) - \varphi(x; a) | + ||x+h|-|x|| \, \varphi(x; a) \quad & & \mbox{by t.i.},  \\
& \le (|x| + 1) \, M_1(a) \, |h| \, \varphi(x; a/2) + |h| \, \varphi(x; a)  \quad & & \mbox{by $|h| < 1$, 
\eqref{eqn:Lipschitz} and t.i.}, \\
& = M_1(a) \cdot (|x|+1) e^{-\frac{a}{4} x^2} \, |h| \, \varphi(x; a/4) + |h| \, \varphi(x; a) & & \\
& \le \{ 1 + M_1(a) \cdot (|x|+1) e^{-\frac{a}{4} x^2} \} |h| \, \varphi(x; a/4) & & 
\mbox{by $\varphi(x; a) \le \varphi(x; a/4)$}, \\
& \le M_4(a) \, |h| \, \varphi(x; a/4).  & & 
\end{alignat*} 
%%%%%%%%%%%%
Thus estimate \eqref{eqn:Lipschitz2} has been proved. \hfill $\Box$
%%%%%%%%%%%%%%%%%%%%%%%%% end proof %%%%%%%%%%%%%%%%%%%%%%%%%%%%%%%%%%%
%%%%%%%%%%%%%%%%%%%%%%%%%%%%%% sec:ds %%%%%%%%%%%%%%%%%%%%%%%%%%%%%%%%%
\section{Dominant Sequences} \label{sec:ds}
%%%%%%%%%%%%%%%%%%%%%%%%%%%%%%%%%%%%%%%%%%%%%%%%%%%%%%%%%%%%%%%%%%%%%
Recall that the hypergeometric series ${}_3g_2(\ba)$ is defined in \eqref{eqn:3g2} 
and the subset $\cS(\Z) \subset \Z^5$ is defined in \eqref{eqn:cS(Z)}.  
In what follows we fix any positive numbers $R$, $\sigma > 0$ and let 
%%%%%%
\[
\A(R, \sigma) := \{\, \ba = (a_0, a_1, a_2; b_1, b_2) \in \C^5 \,:\, 
|\!| \ba |\!| \le R, \,\, \rRe \, s(\ba) > \sigma \, \},  
\] 
%%%%%
where $|\!| \cdot |\!|$ is the standard norm on $\C^5$. 
As an application of \S \ref{sec:dspm} we shall show the following.      
%%%%%%%%%%%%%%%%%%%%%%%%% thm:ds %%%%%%%%%%%%%%%%%%%%%%%%%%%%%%%%%%%%%
\begin{theorem} \label{thm:ds} 
If $\bp = (p_0, p_1, p_2; q_1, q_2) \in \cS(\Z)$ is any vector satisfying either 
%%%%%%%%%%%%%%%%%%%%%%%%% eqn:cond-key %%%%%%%%%%%%%%%%%%%%%%%%%%%%%%%
\begin{equation} \label{eqn:cond-key}
(\mathrm{a}) \quad \vD(\bp) \le 0 \qquad \mbox{or} \qquad 
(\mathrm{b}) \quad  2 q_1^2 - 2 (p_1+p_2) q_1 + p_1 p_2 \ge 0,  
\end{equation}
%%%%%%%%%%%%%%%%%%%%%%%%%%%%%%%%%%%%%%%%%%%%%%%%%%%%%%%%%%%%%%%%%%%%%
where $\vD(\bp)$ is the polynomial in \eqref{eqn:v(p)}, 
then $|D(\bp)| > 1$ and there exists an asymptotic formula  
%%%%%%%%%%%%%%%%%%%%%%%%
\begin{equation*} 
t(\ba) \cdot {}_3g_2(\ba + n \bp) = B(\ba; \bp)   
\cdot D(\bp)^{n} \cdot n^{-s(\sba) -\frac{1}{2} } 
\, \left\{ 1+ O(n^{-\frac{1}{2} }) \right\} \qquad 
\mbox{as} \quad n \to +\infty,   
\end{equation*} 
%%%%%%%%%%%%%%%%%%%%%%%%
uniformly valid with respect to $\ba \in \A(R, \sigma)$, where $D(\bp)$ is 
defined in \eqref{eqn:A(p)} and 
%%%%%%%%%%%%%%%%%%%%%%% eqn:B(a;p) %%%%%%%%%%%%%%%%%%%%%%%%%%%%
\begin{equation} \label{eqn:B(a;p)}
t(\ba) := \sin \pi(b_1-a_0) \cdot \sin \pi(b_2 - a_0),  \quad 
B(\ba; \bp) := \dfrac{\pi^{ \frac{1}{2} } \cdot  
p_0^{a_0- \frac{1}{2} } p_1^{a_1- \frac{1}{2} } p_2^{a_2- \frac{1}{2} } \cdot  
s_2(\bp)^{s(\sba) - 1}}{2^{ \frac{3}{2} } 
\prod_{i=0}^2 \prod_{j=1}^2 (q_j-p_i)^{b_j-a_i- \frac{1}{2} }},    
\end{equation} 
%%%%%%%%%%%%%%%%%%%%%%%%%%%%%%%%%%%%%%%%%%%%%%%%%%%%%%%%%%%%%
with $s_2(\bp) := p_0 p_1 + p_1 p_2 + p_2 p_0 - q_1 q_2$ as in 
Theorem $\ref{thm:rs}$.   
\end{theorem} 
%%%%%%%%%%%%%%%%%%%%%%%%%%%%%%%%%%%%%%%%%%%%%%%%%%%%%%%%%%%%%
%%%%%%%%%%%%%%%%%%%%%%%%% rem:ds %%%%%%%%%%%%%%%%%%%%%%%%%%%%%
\begin{remark} \label{rem:ds} 
Conditions \eqref{eqn:bp} and \eqref{eqn:cond-key} are 
invariant under multiplication of $\bp$ by any positive scalar. 
This homogeneity allows one to restrict $\cS(\R)$ to 
$\cS_1(\R) := \cS(\R) \cap \{ q_1 = 1\}$, which is a $3$-dimensional 
solid tetrahedron.  
A numerical integration shows that the domain in $\cS_1(\R)$ 
bounded by inequalities \eqref{eqn:cond-key} occupies some $43$ \% 
of the whole $\cS_1(\R)$ in volume basis.  
Thus we may say that about $43$ \% of the vectors in $\cS(\Z)$ 
satisfy condition \eqref{eqn:cond-key}. 
\end{remark}
%%%%%%%%%%%%%%%%%%%%%%%%%%%%%%%%%%%%%%%%%%%%%%%%%%%%%%%%%%%%%
\par
%%%%%%
By the definition of ${}_3g_2(\ba)$ one can write $g(n) := 
{}_3g_2(\ba + n \bp) = \sum_{k=0}^{\infty} \varphi(k; n)$ with   
%%%%%%%%%%%%%
\[
\varphi(k; n) :=   
\frac{\vG(k + p_0 n + a_0) \vG(k-(q_1-p_0) n+a_0-b_1+1) 
\vG(k-(q_2-p_0) n+a_0-b_2+1)}{\vG(k+1) \vG(k+(p_0-p_1) n+a_0-a_1+1) 
\vG(k+(p_0-p_2) n+ a_0-a_2+1)}.  
\]
%%%%%%%%%%%% 
We remark that the current $g(n)$ corresponds to the sequence 
$g_0(n)$ in \eqref{eqn:gi(n)}, not to $g(n)$ in \eqref{eqn:gg(n)}. 
In general a gamma factor $\vG(\sigma k + \lambda n + \alpha)$ is said 
to be {\sl positive} resp. {\sl negative} on an interval of $k$, if 
$\sigma k + \lambda n$ is positive resp. negative whenever 
$k$ lies in that interval.   
Since $\bp \in \cS(\Z)$, all lower and an upper gamma factors of 
$\varphi(k; n)$ are positive in $k > 0$, while the remaining two upper 
factors changes their signs when $k$ goes across $(q_1-p_0) n$ or 
$(q_2-p_0) n$.    
Thus it is natural to make a decomposition $g(n) = g_1(n) + g_2(n) + g_2(n)$ 
with  
%%%%%%%%%%%%
\[
g_1(n) := \sum_{k=0}^{(q_1-p_0) n-1} \varphi(k; n), 
\quad 
g_2(n) :=\sum_{k=(q_1-p_0) n}^{(q_2-p_0) n-1} \varphi(k; n), 
\quad 
g_3(n) :=\sum_{k=(q_2-p_0) n}^{\infty} \varphi(k; n),    
\]
%%%%%%%%%%%
where if $q_1 = q_2$ then $g_2(n)$ should be null so 
we always assume $q_1 < q_2$ when discussing $g_2(n)$. 
It turns out that the first component $g_1(n)$ is 
the most dominant among the three, yielding the leading 
asymptotics for $g(n)$.   
The proof of Theorem \ref{thm:ds} is completed at the end 
of \S\ref{subsec:3rd}.   
%%%%%%%%%%%%%%%%%%%%%% subsec:1st %%%%%%%%%%%%%%%%%%%%%%%%%%% 
\subsection{First Component} \label{subsec:1st}
%%%%%%%%%%%%%%%%%%%%%%%%%%%%%%%%%%%%%%%%%%%%%%%%%%%%%%%%%%%%
For the first component $g_1(n)$, applying Euler's reflection formula 
for the gamma function to the two negative gamma factors in the 
numerator of $\varphi(k; n)$, we have  
%%%%%%%%%%%%
\[
t(\ba) \cdot g_1(n) = \pi^2 \cdot (-1)^{(q_1 + q_2) n} \cdot G_1(n)  
\quad \mbox{with} \quad 
G_1(n) := \sum_{k=0}^{L_1 n -1}  
\frac{\vG(\sigma_1 k + \lambda_1 n + \alpha_1)}{\prod_{j=1}^5 \vG(\tau_j k + \mu_j n + \beta_j)}, 
\]
%%%%%%%%%%%%
where $L_1 = q_1-p_0$, $\sigma_1 = 1$, $\lambda_1 = p_0$, $\alpha_1 = a_0$ and 
%%%%%%%%%%%%
\begin{alignat*}{5}
\tau_1 &= 1, \quad & \tau_2 &= 1, \quad & \tau_3 &= 1, \quad & \tau_4 &= -1, \quad & \tau_5 &= -1, \\
\mu_1 &= 0, \quad & \mu_2 &= p_0-p_1, \quad & \mu_3 &= p_0-p_2, \quad & 
\mu_4 &= q_1-p_0 = L_1, \quad & \mu_5 &= q_2-p_0, \\
\beta_1 &= 1, \quad & \beta_2 &= a_0-a_1+1, \quad & \beta_3 &= a_0-a_2+1, \quad & 
\beta_4 &= b_1-a_0, \quad & \beta_5 &= b_2-a_0.   
\end{alignat*}
%%%%%%%%%%%%
\par
%%%%%%%%%%%%  
Under the assumption of Theorem \ref{thm:ds} the sum $G_1(n)$ 
satisfies all conditions in Theorem \ref{thm:dspm}.   
Indeed, balancedness \eqref{eqn:ap-bq2} follows from $s(\bp) = 0$; 
boundedness \eqref{eqn:bdd} is trivial because $\alpha_1$ and $\beta_j$ 
are independent of $n$; admissibility \eqref{eqn:ap-bq1} is fulfilled   
with $r_0 = 0$ and $r_1 = L_1$ due to condition \eqref{eqn:bp};  
genericness \eqref{eqn:punc} is trivial since $I_0 \cup I_1 = \emptyset$ 
with $J_0 = \{1\}$ and $J_1 = \{4\}$ by inequalities in \eqref{eqn:bp}.  
To verify properness \eqref{eqn:max}, notice that the characteristic 
equation \eqref{eqn:max-eq} now reads    
%%%%%%%%%%%
\[
\chi_1(x) = \dfrac{x(x+p_0-p_1)(x+p_0-p_2)}{(-x+q_1-p_0)(-x+q_2-p_0)} - (x+p_0) = 0.   
\]
%%%%%%%%%%%
Thanks to $s(\bp) = 0$ this equation reduces to a linear equation in $x$ 
having the unique root   
%%%%%%%%%%%
\[
x_0 = \frac{p_0(q_1-p_0)(q_2-p_0)}{p_1 p_2-(q_1-p_0)(q_2-p_0)}  
= \frac{p_0(q_1-p_0)(q_2-p_0)}{s_2(\bp)},  
\]
%%%%%%%%%%%
where $s(\bp) = 0$ again leads to $s_2(\bp) = 
p_1 p_2-(q_1-p_0)(q_2-p_0)$, which together with \eqref{eqn:bp} 
yields $s(\bp) -p_0(q_2-p_0) = (q_1-p_1)(q_1-p_2) > 0$ and hence   
$s_2(\bp) > p_0(q_2-p_0) > 0$, that is,  
%%%%%%%%%%%
\[ 
0 < x_0 < L_1 = q_1 - p_0. 
\]
%%%%%%%%%%% 
If $\phi_1(x)$ is the additive phase function for $G_1(n)$ then 
it follows from \eqref{eqn:dash2} and \eqref{eqn:bp} that     
%%%%%%%%%%%
\begin{align*}
\phi_1''(x_0) &= 
\frac{1}{x_0} + \frac{1}{x_0+p_0-p_1} + \frac{1}{x+p_0-p_2} + 
\frac{1}{q_1-p_0-x_0} +  \frac{1}{q_2-p_0-x_0} - \frac{1}{x_0+p_0} \\[2mm]
&= \frac{s_2(\bp)^4}{p_0 p_1 p_2 \prod_{i=0}^2 \prod_{j=1}^2 (q_j-p_i)} > 0. 
\end{align*}
%%%%%%%%%%
Thus in the interval $0 < x < L_1$ the function $\phi_1(x)$ has only one 
local and hence global minimum at $x = x_0$, which is non-degenerate.  
Therefore properness \eqref{eqn:max} is satisfied with $\Mm = \{x_0 \}$ 
and hence Theorem \ref{thm:dspm} applies to the sum $G_1(n)$. 
%%%%%%%%%%%%%%%%%%%%%%%% lem:1st %%%%%%%%%%%%%%%%%%%%%%%%%%%%%
\begin{lemma} \label{lem:1st} 
For any $\bp \in \cS(\Z)$ we have $|D(\bp)| > 1$ and an 
asymptotic representation 
%%%%%
\[
t(\ba) \cdot g_1(n) = B(\ba; \bp)   
\cdot D(\bp)^{n} \cdot n^{-s(\sba) -\frac{1}{2} } 
\, \left\{ 1+ O(n^{-\frac{1}{2} }) \right\} \qquad \mbox{as} \quad 
n \to +\infty,     
\]
%%%%%
uniform with respect to $\ba \in \A(R, \sigma)$, where $D(\bp)$, $t(\ba)$ 
and $B(\ba; \bp)$ are as in \eqref{eqn:A(p)} and \eqref{eqn:B(a;p)}.    
\end{lemma}
%%%%%%%%%%%%%%%%%%%%%%%%%%%%%%%%%%%%%%%%%%%%%%%%%%%%%%%%%%%%
%%%%%%%%%%%%%%%%%%%%%%% begin proof %%%%%%%%%%%%%%%%%%%%%%%%%%
{\it Proof}.   
Substituting $x = x_0$ into formulas \eqref{eqn:G-asymp2} and 
using $s(\bp) = 0$ repeatedly, one has  
%%%%%%%%%%
\[
(\Phi_1)_{\mathrm{max}} = \Phi_1(x_0) 
= \dfrac{p_0^{p_0} p_1^{p_1} p_2^{p_2}}{\prod_{i=0}^2 \prod_{j=1}^2 
(q_j - p_i)^{q_j - p_i}},  \quad 
u_1(x_0) 
= \dfrac{p_0^{a_0-1} p_1^{a_1-1} p_2^{a_2-1} 
s_2(\bp)^{s(\sba) + 1}}{(2 \pi)^2 \prod_{i=0}^2 \prod_{j=1}^2 (q_j-p_i)^{b_j-a_i}}, 
\]
%%%%%%%%%%%
while $\gamma_1 := \gamma(n)$ in definition \eqref{eqn:ga(n)} now reads 
$\gamma_1 = - s(\ba) - 1$.  
Since $\delta_0(n) = \delta_1(n) = 1$ in \eqref{eqn:uniform} by 
$I_0 \cup I_1 = \emptyset$, formula \eqref{eqn:dspm} in 
Theorem \ref{thm:dspm} implies that  
%%%%%%%%%%%
\[
G_1(n) = C_1 \cdot (\Phi_1)_{\mathrm{max}}^n \cdot 
n^{\gamma_1 + \frac{1}{2} } \cdot \left\{ 1 + O(n^{- \frac{1}{2} }) \right\} 
\qquad \mbox{as} \quad n \to + \infty, 
\]
%%%%%%%%%%%
where formula \eqref{eqn:PhimC} allows one to calculate the constant 
$C_1 := C(n)$ as  
%%%%%%%%%%%
\[ 
C_1  
= \sqrt{2 \pi} \, \frac{u_1(x_0)}{\sqrt{\phi_1''(x_0)}} 
= \dfrac{p_0^{a_0- \frac{1}{2} } p_1^{a_1- \frac{1}{2} } p_2^{a_2- \frac{1}{2} } 
s_2(\bp)^{s(\sba) - 1}}{(2 \pi)^{ \frac{3}{2} } 
\prod_{i=0}^2 \prod_{j=1}^2 (q_j-p_i)^{b_j-a_i- \frac{1}{2} }}. 
\] 
%%%%%%%%%% 
In view of the relation between $G_1(n)$ and $g_1(n)$ the above 
asymptotic formula for $G_1(n)$ gives the one for 
$g_1(n)$.  
Finally $|D(\bp)| > 1$ follows from 
Lemma \ref{lem:Phi1(0)} below. \hfill $\Box$ 
%%%%%%%%%%%%%%%%%%%%%%% end proof %%%%%%%%%%%%%%%%%%%%%%%%%%%%
%%%%%%%%%%%%%%%%%%%%%%% lem:Phi1(0) %%%%%%%%%%%%%%%%%%%%%%%%%%%
\begin{lemma} \label{lem:Phi1(0)} 
Under condition \eqref{eqn:bp} one has  
$|D(\bp) | = (\Phi_1)_{\mathrm{max}} > \Phi_1(0) > 1$. 
\end{lemma}
%%%%%%%%%%%%%%%%%%%%%%% begin proof %%%%%%%%%%%%%%%%%%%%%%%%%%%
{\it Proof}. 
First, $|D(\bp)| = (\Phi_1)_{\mathrm{max}}$ is obvious from the definition 
\eqref{eqn:A(p)} of $D(\bp)$ and the expression for 
$(\Phi_1)_{\mathrm{max}}$, while $(\Phi_1)_{\mathrm{max}} > \Phi_1(0)$ 
is also clear from $\Mm = \{ x_0 \}$.   
Regarding $\bp = (p_0,p_1,p_2; q_1, q_2)$ as real variables 
subject to the linear relation $s(\bp) = 0$ and ranging over the 
closure of the domain \eqref{eqn:bp}, we shall find the minimum of  
%%%%%
\[
\Phi_1(0) = \frac{p_0^{p_0}}{(p_0-p_1)^{p_0-p_1} (p_0-p_2)^{p_0-p_2}
(q_1-p_0)^{q_1-p_0}(q_2-p_0)^{q_2-p_0}}.  
\]
%%%%%
For any fixed $(p_0, p_1, p_2)$, due to the constraint $s(\bp) = 0$, 
one can thought of $\Phi_1(0)$ as a function of single variable $q_1$ 
in the interval $p_0 \le q_1 \le p_1 + p_2$. 
Differentiation with respect to $q_1$ shows that $\Phi_1(0)$ attains  
its minimum (only) at the endpoints $q_1 = p_0$, $p_1+p_2$, 
whose value is   
%%%%%%
\[ 
\Psi(p_0, p_1, p_2) :=  
\frac{p_0^{p_0}}{(p_0-p_1)^{p_0-p_1} (p_0-p_2)^{p_0-p_2}
(p_1+p_2-p_0)^{p_1+p_2-p_0}}.   
\]
%%%%%%%%%%%%%%%%%%%%%%%%%%%%%%%%%%%%%%%%%%%%%%%%%%%%%%%%%%%%%
so $\Phi_1(0) > \Psi(p_0, p_1, p_2)$ for any $p_0 < q_1 < p_1+p_2$. 
With a fixed $p_0 > 0$ we think of $\Psi(p_0, p_1, p_2)$ as a function of 
$(p_1, p_2)$ in the closed simplex $p_0 \le p_1 + p_2$, $p_1 \le p_0$, 
$p_2 \le p_0$. 
It has a unique critical value $\Psi(p_0, 2 p_0/3, 2 p_0/3) = 3^{p_0} > 1$ 
in the interior of the simplex, while on its boundary one has 
$\Psi(p_0, \alpha, p_0)  = 
\Psi(p_0, p_0, \alpha)=  \Psi(p_0, \alpha, p_0-\alpha) = p_0^{p_0} 
\alpha^{-\alpha} (p_0-\alpha)^{\alpha-p_0} \ge 1$ for any 
$0 \le \alpha \le p_0$.  
Therefore we have $\Phi_1(0) > \Psi(p_0, p_1, p_2) \ge 1$ 
under condition \eqref{eqn:bp}.  \hfill $\Box$ 
%%%%%%%%%%%%%%%%%%%%% end proof %%%%%%%%%%%%%%%%%%%%%%%%%%%%%%%   
%%%%%%%%%%%%%%%%%%%%%%%% subsec:2nd %%%%%%%%%%%%%%%%%%%%%%%%%%
\subsection{Second Component} \label{subsec:2nd}
%%%%%%%%%%%%%%%%%%%%%%%%%%%%%%%%%%%%%%%%%%%%%%%%%%%%%%%%%%%%%
Taking the shift $k \mapsto k + (q_1-p_0) n$ in $\varphi(k; n)$ 
(see \eqref{eqn:shift}) and applying the reflection formula to the 
unique negative gamma factor in the numerator of 
$\varphi(k+(q_1-p_0) n; n)$, one has     
%%%%%%%%%%%%
\[
g_2(n) = \frac{\pi \cdot (-1)^{(q_2 - q_1) n} \, G_2(n) }{\sin \pi(b_2-a_0)} 
\quad \mbox{with} \quad 
G_2(n) := \sum_{k=0}^{L_2 n -1}  (-1)^k
\frac{ \prod_{i=1}^2 \vG(\sigma_i k + \lambda_i n + \alpha_i)}{\prod_{j=1}^4 
\vG(\tau_j k + \mu_j n + \beta_j)}, 
\]
%%%%%%%%%%%%
where $L_2 = q_2-q_1 > 0$, $\sigma_1 = \sigma_2 = 1$, $\lambda_1 = q_1$, 
$\lambda_2 = 0$, $\alpha_1 = a_0$, $\alpha_2 = a_0-b_1+1$ and 
%%%%%%%%%%%%
\begin{alignat*}{4}
\tau_1 &= 1, \quad & \tau_2 &= 1, \quad & \tau_3 &= 1, \quad & \tau_4 &= -1, \\
\mu_1 &= q_1-p_0, \quad & \mu_2 &= q_1-p_1, \quad & \mu_3 &= q_1-p_2, \quad & 
\mu_4 &= q_2-q_1 = L_2, \\
\beta_1 &= 1, \quad & \beta_2 &= a_0-a_1+1, \quad & \beta_3 &= a_0-a_2+1, \quad & 
\beta_4 &= b_2-a_0.    
\end{alignat*} 
%%%%%%%%%%%%%%%
\par
%%%%%%%%%%%%%%%
Rewriting $k \mapsto 2 k$ or $k \mapsto 2 k+1$ according as $k$ is even or odd,  
we have a decomposition $G_2(n) = G_{20}(n) - G_{21}(n) + H_2(n)$, 
where $G_{2\nu}(n)$ is given by   
%%%%%%%%%%%%%%%
\[
G_{2\nu }(n) := \sum_{k=0}^{\lceil \frac{L_2}{2} n \rceil -1}   
\frac{ \prod_{i=1}^2 \vG(2 \sigma_i k + \lambda_i n + \alpha_i + 
\nu \, \sigma_i) }{\prod_{j=1}^4 \vG(2 \tau_j k + \mu_j n + \beta_j + \nu \, \tau_j )}, 
\qquad \nu = 0, 1, 
\]
%%%%%%%%%%%%%%%
while if $L_2$ or $n$ is even then $H_2(n) := 0$; otherwise, i.e.,  
if both of $L_2$ and $n$ are odd then  
%%%%%%%%%%%%%%%
\[
H_2(n) := \dfrac{\prod_{i=1}^2 \vG((\sigma_i L_2 + \lambda_i) n + 
\alpha_i) }{ \prod_{j=1}^4 \vG( (\tau_j L_2 + \mu_j) n + \beta_j ) }.  
\]
%%%%%%%%%%%%%%%
\par
%%%%%%%%%%%%%%%
Obviously, $G_{20}(n)$ and $G_{21}(n)$ have the same multiplicative 
phase function, which we denote by $\Phi_2(x)$. 
Let $\phi_2(x) := - \log \Phi_2(x)$ be the associated additive phase function. 
In order to make the second component $g_2(n)$ weaker than the 
first one $g_1(n)$, we want to make $\phi_2'(x) \ge 0$ or equivalently 
$\chi_2(x) \ge 0$ for every $0 \le x \le L_2/2$, where $\chi_2(x)$ is 
the common characteristic function \eqref{eqn:max-eq} for the 
sums $G_{20}(n)$ and $G_{21}(n)$, which is given by 
%%%%%%
\[
\chi_2(x) = 
\frac{(2 x + \mu_1)^2(2 x + \mu_2)^2 (2 x + \mu_3)^2}{(-2 x + L_2)^2} -
(2 x + \lambda_1)^2 (2 x +\lambda_2)^2.   
\]
%%%%%%
The non-negativity of $\chi_2(x)$ in the interval $0 \le x \le L_2/2$ 
is equivalent to   
%%%%%%%%%%%%%%%%%%%%%%%%%% eqn:chi2 %%%%%%%%%%%%%%%%%%%%%%%%%%%%%%%
\begin{equation} \label{eqn:chi2} 
\begin{split}
\chi(x; \bp)  
&:= (x + \mu_1)(x + \mu_2) (x + \mu_3) + 
(x + \lambda_1) (x +\lambda_2) (x -L_2)  \\
&= (x+q_1-p_0)(x+q_1-p_1)(x+q_1-p_2)+ x (x+q_1) (x+q_1-q_2) \\ 
&\ge 0 \qquad \mbox{for any} \qquad 0 \le x \le L_2 = q_2 - q_1.  
\end{split}
\end{equation}
%%%%%%%%%%%%%%%%%%%%%%%%%%%%%%%%%%%%%%%%%%%%%%%%%%%%%%%%%%%%%%%%%
\par
%%%%%%
It is easy to see that $G_{20}(n)$ and $G_{21}(n)$ satisfy balancedness 
\eqref{eqn:ap-bq2}, boundedness \eqref{eqn:bdd} and admissibility 
\eqref{eqn:ap-bq1} conditions, where $r_0 = 0$, $r_1 = L_2/2$ and 
$I_0 = \{2\}$, $I_1 = J_0 = \emptyset$, $J_1 = \{4\}$, while 
genericness \eqref{eqn:punc} for $G_{2\nu}(n)$ becomes  
$b_1-a_0 \not\in \Z_{\ge \nu+1}$ for $\nu = 0, 1$.  
%%%%%%%%%%%%%%%%%%%%%%%%%% lem:2nd %%%%%%%%%%%%%%%%%%%%%%%%%%%%%%%
\begin{lemma} \label{lem:2nd} 
Under the assumption of Lemma $\ref{lem:1st}$, if $\bp$ satisfies the 
additional condition \eqref{eqn:chi2} then there exist positive constants  
$0 < d_2 < |D(\bp)|$, $C_2 > 0$ and $N_2 \in \N$ such that 
%%%%%
\[
|t(\ba) \cdot g_2(n)| \le C_2 \cdot d_2^n,  
\qquad {}^{\forall} n \ge N_2, \, \, {}^{\forall}\ba \in \A(R, \sigma).   
\]
%%%%%
\end{lemma} 
%%%%%%%%%%%%%%%%%%%%%%%%%% begin proof %%%%%%%%%%%%%%%%%%%%%%%%%%%%
{\it Proof}. 
Condition \eqref{eqn:chi2} implies that $\Phi_2(x)$ is decreasing 
everywhere in $0 \le x \le L_2/2$ and is strictly so near 
$x = 0$ since $\chi(0; \bp) = (q_1-p_0)(q_1-p_1)(q_1-p_2) > 0$ by 
condition \eqref{eqn:bp}.  
Hence $\Phi_2(x)$ attains its maximum (only) at the left end 
$x = 0$ of the interval, having the value  
%%%%%%
\[
(\Phi_2)_{\mathrm{max}} 
= \Phi_2(0) = \frac{q_1^{q_1}}{(q_1-p_0)^{q_1-p_0} (q_1-p_1)^{q_1-p_1}
(q_1-p_2)^{q_1-p_2}(q_2-q_1)^{q_2-q_1}}  
= \Phi_1(L_1).  
\]
%%%%% 
whereas $(\Phi_2)_{\mathrm{max}} = \Phi_1(L_1) < (\Phi_1)_{\mathrm{max}} 
= |D(\bp)|$ follows from Lemma \ref{lem:Phi1(0)}. 
Thus if $d_2$ is any number such that $(\Phi_2)_{\mathrm{max}} 
< d_2 < |D(\bp)|$, then Proposition \ref{prop:irregular} shows that 
%%%%%%
\[
|G_{2\nu}(n)| \le \dfrac{K_2 \cdot d_2^n}{ \min \{1, \, 
\mathrm{dist}(b_1-a_0, \, \Z_{\ge \nu+1}) \}} 
\le \dfrac{K_2 \cdot d_2^n}{\delta(b_1-a_0)}, 
\quad {}^{\forall} n \ge N_2, \,\, \nu = 0, 1,  
\]
%%%%%% 
for some $K_2 > 0$ and $N_2 \in \N$, where $\delta(z) := \min \{ 1, \, 
\mathrm{dist}(z, \, \N)\}$ for $z \in \C$.  
%%%%%%
\par
%%%%%%
We have to take care of $H_2(n)$ when $L_2$ and $n$ are both odd. 
Stirling's formula \eqref{eqn:stirling} yields 
%%%%%%
\[
H_2(n) = \frac{1}{2 \pi} \cdot \frac{\prod_{i=1}^2 
(\sigma_i L_2 + \lambda_i)^{\alpha_i- \frac{1}{2}} }{ \prod_{j=1}^4 
(\tau_j L_2 + \mu_j)^{\beta_j - \frac{1}{2} } } 
\cdot \Phi_2(L_2/2)^n \cdot n^{\gamma_2} \cdot \{1 + O(1/n) \} 
\quad \mbox{as} \quad n \to +\infty, 
\]
%%%%%%
where $\gamma_2 := \alpha_1+\alpha_2-\beta_1-\beta_2-\beta_3-\beta_4+1$. 
Since $\Phi_2(L_2/2) < (\Phi_2)_{\mathrm{max}} < d_2$, upon retaking $K_2 > 0$ 
suitably, one has $|H_2(n)| \le K_2 \cdot d_2^n \le K_2 \cdot 
d_2^n/ \delta(b_1-a_0)$ for any $n \ge N_2$.  
%%%%%
\par
%%%%%
Then from the relation between $g_2(n)$ and 
$G_2(n) = G_{20}(n)-G_{21}(n) + H_2(n)$ one has   
%%%%%%
\[
|t(\ba) \cdot g_2(n) | \le 3 \pi K_2 \cdot M_2(\ba) \cdot d_2^n  
\quad \mbox{with} \quad  
M_2(\ba) := \dfrac{|\sin \pi(b_1-a_0)|}{\delta(b_1-a_0)}.  
\]
%%%%%
Since $M_2(\ba)$ is bounded for $\ba \in \A(R, \sigma)$ the lemma 
follows (here $\sigma$ is irrelevant). \hfill $\Box$ \par\medskip
%%%%%%%%%%%%%%%%%%%%%%% end proof %%%%%%%%%%%%%%%%%%%%%%%%%%%%%%%%
Lemma \ref{lem:2nd} tempts us to ask when condition \eqref{eqn:chi2} 
is satisfied.   
%%%%%%%%%%%%%%%%%%%%%%% lem:key %%%%%%%%%%%%%%%%%%%%%%%%%%%%%%%%%%
\begin{lemma}  \label{lem:key} 
For any $\bp \in \cS(\R)$ condition \eqref{eqn:cond-key} implies 
condition \eqref{eqn:chi2}.  
\end{lemma} 
%%%%%%%%%%%%%%%%%%%%%%% begin proof %%%%%%%%%%%%%%%%%%%%%%%%%%%%%%
{\it Proof}. 
We use the following general fact. 
Let $\chi(x)$ be a real cubic polynomial with positive leading coefficient 
and $\vD$ be its discriminant. 
If $\vD < 0$ then $\chi(x)$ has only one real root so that once 
$\chi(c_0) > 0$ for some $c_0 \in \R$ then $\chi(x) > 0$ for every 
$x \ge c_0$. 
Even if $\vD = 0$, once $\chi(c_0) > 0$ then $\chi(x) \ge 0$ for every  
$x \ge c_0$ with possible equality $\chi(c_1) = 0$, $c_1 > c_0$,  
only if $\chi(x)$ attains a local minimum at $x = c_1$. 
Currently, $\chi(x; \bp)$ has discriminant $\vD(\bp)$ in formula 
\eqref{eqn:v(p)} and $\chi(0; \bp) = (q_1-p_0)(q_1-p_1)(q_1-p_2) > 0$ by 
condition \eqref{eqn:bp}.   
Thus if $\vD(\bp) \le 0$ then $\chi(x; \bp) \ge 0$ for every $x \ge 0$; 
this is just the case (a) in condition \eqref{eqn:cond-key}.  
%%%%%%
\par
%%%%%
We proceed to the case (b) in \eqref{eqn:cond-key}. 
The derivative of $\chi(x; \bp)$ in $x$ is given by 
%%%%%%
\[
\chi'(x; \bp) = 6 x^2 +4 (2 q_1-q_2) x
+ (3 q_1-p_1-p_2)(p_1+p_2-q_2)  + 2 q_1^2 -2 (p_1+p_2) q_1 + p_1 p_2.  
\]
%%%%%%
Note that $2 q_1 - q_2$, $3 q_1-p_1-p_2$, $p_1+p_2-q_2 > 0$ 
by condition \eqref{eqn:bp}. 
Having axis of symmetry $x = -(2 q_1-q_2)/3 < 0$, the 
quadratic function $\chi'(x; \bp)$ is increasing in $x \ge 0$ and hence    
%%%%%%
\begin{align*}
\chi'(x; \bp) \ge \chi'(0; \bp) 
&= (3 q_1-p_1-p_2)(p_1+p_2-q_2)  + 2 q_1^2 -2 (p_1+p_2) q_1 + p_1 p_2 \\[1mm]
&> 2 q_1^2 -2 (p_1+p_2) q_1 + p_1 p_2 \ge 0 
\qquad \mbox{for any} \quad x \ge 0,  
\end{align*}
%%%%%%
where the last inequality stems from (b) in condition \eqref{eqn:cond-key}.  
Thus $\chi(x; \bp) \ge \chi(0; \bp) > 0$ for any $x \ge 0$, so condition 
\eqref{eqn:chi2} is satisfied. \hfill $\Box$ \par\medskip
%%%%%%%%%%%%%%%%%%%%%%% end proof %%%%%%%%%%%%%%%%%%%%%%%%%%%%%%%%
The converse to the implication in Lemma \ref{lem:key} is also true, 
accordingly conditions \eqref{eqn:cond-key} and \eqref{eqn:chi2} are 
equivalent for any $\bp \in \cS(\R)$, but the proof of this fact is omitted 
as it is not needed in this article.  
In the situation of Lemma \ref{lem:2nd} we proceed to the third component. 
%%%%%%%%%%%%%%%%%%%%%%% subsec:3rd %%%%%%%%%%%%%%%%%%%%%%%%%%%%%%%
\subsection{Third Component} \label{subsec:3rd}
%%%%%%%%%%%%%%%%%%%%%%%%%%%%%%%%%%%%%%%%%%%%%%%%%%%%%%%%%%%%%%%%
For the third component $g_3(n)$, taking the shift 
$k \mapsto k+(q_2-p_0)n$ in $\varphi(k; n)$, one has      
%%%%%%%%%%%%%%%
\[
g_3(n) = \sum_{k=0}^{\infty} \frac{\prod_{i=1}^3 
\vG(\sigma_i k+ \lambda_i n+ \alpha_i)}{\prod_{j=1}^3
\vG(\tau_j k+ \mu_j n +\beta_j)}, 
\]
%%%%%%%%%%%%%%%
where $\sigma_1 = \sigma_2 = \sigma_3 = \tau_1 = \tau_2 = \tau_3 = 1$ and 
%%%%%%%%%%%%%%%
\begin{alignat*}{6}
\lambda_1 &= q_2, \quad & \lambda_2 &= q_2-q_1, \quad & \lambda_3 &= 0, \quad &
\alpha_1 &= a_0, \quad & \alpha_2 &= a_0-b_1+1, \quad & \alpha_3 &= a_0-b_2+1, 
\\ 
\mu_1 &= q_2-p_0, \quad & \mu_2 &= q_2-p_1, \quad & \mu_3 &= q_2-p_2, \quad &
\beta_1 &= 1, \quad  &  \beta_2 &= a_0-a_1+1, \quad &  \beta_3 &= a_0-a_2+1.  
\end{alignat*}
%%%%%%%%%%%%%%%
\par
%%%%%%
It is easy to see that $g_3(n)$ satisfies balancedness 
\eqref{eqn:ap-bq2}, boundedness \eqref{eqn:bdd}, admissibility 
\eqref{eqn:ap-bq1} with $r_0 = 0$ and $r_1 = +\infty$.  
Notice that $I_0 = \{3 \}$ if $q_1 < q_2$ and $I_0 = \{2, 3\}$ if 
$q_1 = q_2$, while $I_1 = J_0 = J_1 = \emptyset$. 
Genericness \eqref{eqn:punc} becomes $b_2-a_0 \not \in \N$ if $q_1 < q_2$,  
and $b_1-a_0$, $b_2-a_0 \not \in \N$ if $q_1 = q_2$.  
Continuity at infinity \eqref{eqn:infty} is satisfied with   
$\bsigma^{\sbsigma} = \btau^{\sbtau} = 1$; convergence 
condition \eqref{eqn:conv} is equivalent to $\rRe \, s(\ba) \ge \sigma$.  
Under the assumption of Lemma $\ref{lem:2nd}$ we have the following. 
%%%%%%%%%%%%%%%%%%%%%%%%%% lem:3rd %%%%%%%%%%%%%%%%%%%%%%%%%%%%%%
\begin{lemma} \label{lem:3rd} 
There exist positive constants $0 < d_3 < |D(\bp)|$, $C_3 > 0$ and 
$N_3 \in \N$ such that 
%%%%%
\[
|t(\ba) \cdot g_3(n)| \le C_3 \cdot d_3^n, 
\qquad {}^{\forall} n \ge N_3, \, \, {}^{\forall}\ba \in \A(R, \sigma).   
\]
%%%%%
\end{lemma}
%%%%%%%%%%%%%%%%%%%%%%%%% begin proof %%%%%%%%%%%%%%%%%%%%%%%%%%%%%
{\it Proof}.  
In view of $s(\bp) = 0$ the characteristic function 
\eqref{eqn:max-eq} for $g_3(n)$ is given by 
%%%%%%%%%%%%%%%
\begin{align*}
\chi_3(x) 
&= (x+q_2-p_0)(x+q_2-p_1)(x+q_2-p_2) - (x+q_2)(x+q_2-q_1)x \\
&= s(\bp) \, x^2 + \{ 2 s(\bp) \, q_2 + s_2(\bp) \} x 
+ (q_2-p_0)(q_2-p_1)(q_2-p_2) \\
&= s_2(\bp) \, x + (q_2-p_0)(q_2-p_1)(q_2-p_2).  
\end{align*}
%%%%%%%%%%%%%%% 
Since $s_2(\bp) > 0$ and $(q_2-p_0)(q_2-p_1)(q_2-p_2) > 0$ from 
condition \eqref{eqn:bp}, one has $\chi_3(x) > 0$ and hence the additive phase 
function $\phi_3(x)$ satisfies $\phi_3'(x) > 0$ for any $x \ge 0$.     
Thus $\Phi_3(x) = e^{-\phi_3(x)}$ is strictly decreasing in $x \ge 0$ and 
attains its maximum (only) at $x = 0$, having the value 
%%%%%%%%%%%%%%%
\[
(\Phi_3)_{\mathrm{max}} = \Phi_3(0) = 
\frac{q_2^{q_2} (q_2-q_1)^{q_2-q_1}}{(q_2-p_0)^{q_2-p_0} 
(q_2-p_1)^{q_2-p_1} (q_2-p_2)^{q_2-p_2}} = \Phi_2(L_2/2),        
\]
%%%%%%%%%%%%%%
whereas $(\Phi_3)_{\mathrm{max}} = \Phi_2(L_2/2) < (\Phi_2)_{\mathrm{max}} 
= \Phi_2(0) = \Phi_1(L_1) < (\Phi_1)_{\mathrm{max}} = |D(\bp)|$. 
Thus if $d_3$ is any number with $(\Phi_3)_{\mathrm{max}} < d_3 < |D(\bp)|$ 
then Proposition \ref{prop:irregular} implies that for any $n \ge N_3$, 
%%%%%%
\[
|g_3(n)| \le \frac{K_3 \cdot d_3^n}{\delta(b_2-a_0)} \quad \mbox{if} \quad q_1 < q_2; 
\qquad 
|g_3(n)| \le \frac{K_3 \cdot d_3^n}{\delta(b_1-a_0) \cdot \delta(b_2-a_0)} 
\quad \mbox{if} \quad q_1 = q_2,   
\]
%%%%%
where the function $\delta(z)$ is defined in the proof of Lemma \ref{lem:2nd}. 
Since $\sin \pi(b_j-a_0) / \delta(b_j-a_0)$, $j = 1, 2$, are bounded 
for $\ba \in \A(R, \sigma)$, the lemma follows immediately.  
\hfill $\Box$ \par\medskip 
%%%%%%%%%%%%%%%%%%%%%%%%%%%%%% end proof %%%%%%%%%%%%%%%%%%%%%%%%%%%%%%%%
Theorem \ref{thm:ds} is now an immediate consequence of Lemmas \ref{lem:1st}, 
\ref{lem:2nd}, \ref{lem:key} and \ref{lem:3rd}. 
Theorems \ref{thm:rs} and \ref{thm:ds} then imply that  if  the shift vector 
$\bp \in \cS(\Z)$ satisfies condition \eqref{eqn:cond-key} then $f(n)$ in 
\eqref{eqn:f(n)} and $g(n)$ in \eqref{eqn:gg(n)} are recessive and dominant 
solutions to the recurrence relation \eqref{eqn:3trr-x} whose coefficients 
$q(n)$ and $r(n)$ are given by \eqref{eqn:q(n)} and \eqref{eqn:r(n)}.  
Now it is almost ready to apply the general error estimate \eqref{eqn:error} 
to $X(n) = f(n)$ and $Y(n) = g(n)$, where a precise asymptotic formula for the 
ratio $R(n) = f(n+2)/g(n+2)$ is available from Theorems \ref{thm:rs} and \ref{thm:ds}. 
%%%%%%%%%%%%%%%%%%%%%%%%%%%%%% sec:casorati %%%%%%%%%%%%%%%%%%%%%%%%%%%%%%
\section{Casoratian and Error Estimates} \label{sec:casorati}
%%%%%%%%%%%%%%%%%%%%%%%%%%%%%%%%%%%%%%%%%%%%%%%%%%%%%%%%%%%%%%%%%%%%%%%
All that remain are to evaluate the initial term $\omega(0)$ of the Casoratian 
determinant 
%%%%%%
\[
\omega(n) := f(n) \cdot g(n+1) - f(n+1) \cdot g(n), 
\]       
%%%%%
and to incorporate the ensuing formula with the asymptotic representation 
for $R(n)$ to complete the proofs of Theorems \ref{thm:cf-straight} and 
\ref{thm:cf-cyclic}. 
The first task is done in \S\ref{subsec:casorati}, while the second in 
\S\ref{subsec:error}. 
%%%%%%%%%%%%%%%%%%%%%%%%%%%%%% subsec:casorati %%%%%%%%%%%%%%%%%%%%%%%%%%%
\subsection{Casoratian} \label{subsec:casorati} 
%%%%%%%%%%%%%%%%%%%%%%%%%%%%%%%%%%%%%%%%%%%%%%%%%%%%%%%%%%%%%%%%%%%%%%% 
In order to evaluate $\omega(0)$, following \cite[formulas (7), (8) and (10)]{EI} 
we define      
%%%%%%%%%%%%%%%%%%%%%%%%%%%%%% tab:invol %%%%%%%%%%%%%%%%%%%%%%%%%%%%%%%%%
\begin{table}
\begin{align*}
\sigma_0^{(0)}(\ba) &:= \ba = (a_0, a_1, a_2; b_1, b_2), \\ 
\sigma_1^{(0)}(\ba) &:= (a_0+1-b_1, a_1+1-b_1, a_2+1-b_1; 2-b_1, b_2+1-b_1), \\
\sigma_2^{(0)}(\ba) &:= (a_0+1-b_2, a_1+1-b_2, a_2+1-b_2; b_1+1-b_2, 2-b_2), \\ 
\sigma_0^{(\infty)}(\ba) &:= (a_0, a_0+1-b_1, a_0+1-b_2; a_0+1-a_1, a_0+1-a_2), \\
\ba^* &:= (1-a_0, 1-a_1, 1-a_2; 2-b_1, 2-b_2). 
\end{align*}
\caption{Five parameter involutions (including identity).}
\label{tab:invol}
\end{table}
%%%%%%%%%%%%%%%%%%%%%%%%%%%%%%%%%%%%%%%%%%%%%%%%%%%%%%%%%%%%%%%%%%%%%%%
\begin{align*}
y_i^{(0)}(\ba; z) &:= z^{1-b_i} {}_3f_2(\sigma_i^{(0)}(\ba); z), \quad i = 0, 1, 2, \quad b_0 := 1,  \\
y_0^{(\infty)}(\ba; z) &:= e^{\mathrm{i} \pi s(\sba)} z^{-a_0}  {}_3f_2(\sigma_0^{(\infty)}(\ba); 1/z), 
\end{align*}
%%%%%%%%%%% 
where $\sigma_i^{(\nu)}$ are  involutions on the parameters $\ba$ as in 
Table \ref{tab:invol}, and put $y_i^{(\nu)}(\ba) := y_i^{(\nu)}(\ba; 1)$.   
Note that $y_0^{(0)}(\ba) = {}_3f_2(\ba)$ and $y_0^{(\infty)}(\ba) = 
e^{\mathrm{i} \pi s(\sba)} {}_3g_2(\ba)$. 
Moreover let $\1 := (1,1,1;1,1)$.    
%%%%%%%%%%%%%%%%%%%%%%%%%%%%%% lem:wronskian %%%%%%%%%%%%%%%%%%%%%%%%%%%%
\begin{lemma} \label{lem:wronskian} 
For any $\ba \in \C^5$ with $\rRe \, s(\ba) > 1$ one has   
%%%%%%%%%%%%%%%%%%%%%%%%%%%%% eqn:wronskian %%%%%%%%%%%%%%%%%%%%%%%%%%%%%
\begin{equation} \label{eqn:wronskian}
\begin{split}
W(\ba) 
&:= y_0^{(0)}(\ba) \cdot y_0^{(\infty)}(\ba+\1) - y_0^{(0)}(\ba+\1) \cdot y_0^{(\infty)}(\ba) 
\\[1mm]
&= - e^{\mathrm{i} \pi s(\sba)} 
\dfrac{\vG(a_0) \vG(a_1) \vG(a_2) \vG(a_0-b_1+1) \vG(a_0-b_2+1) 
\vG(s(\ba)-1)}{\vG(b_1-a_1) \vG(b_1-a_2) \vG(b_2-a_1) \vG(b_2-a_2)}.  
\end{split}
\end{equation}
%%%%%%%%%%%%%%%%%%%%%%%%%%%%%%%%%%%%%%%%%%%%%%%%%%%%%%%%%%%%%%%%%%%%%%%  
\end{lemma}
%%%%%%%%%%%%%%%%%%%%%%%%%%%%%%%%%%%%%%%%%%%%%%%%%%%%%%%%%%%%%%%%%%%%%%%
%%%%%%%%%%%%%%%%%%%%%%%%%%%%%% begin proof %%%%%%%%%%%%%%%%%%%%%%%%%%%%%%
{\it Proof}. 
A careful inspection of Bailey \cite[\S 10.3, formulas (3) and (5)]{Bailey1} shows that 
%%%%%%%%%%%
\begin{align*}
w(\ba; z) 
&:= y_0^{(0)}(\ba; z) \cdot y_1^{(0)}(\ba+\1; z)  
- y_0^{(0)}(\ba+\1; z) \cdot y_1^{(0)}(\ba; z) \\[1mm]
&= \dfrac{\vG(a_0) \vG(a_1) \vG(a_2) \vG(a_0-b_1+1) \vG(a_1-b_1+1) 
\vG(a_2-b_1+1)}{\vG(b_1) \vG(1-b_1) \vG(b_2-a_0)\vG(b_2-a_1) \vG(b_2-a_2)} \\
&\phantom{=} \times z^{1-b_1-b_2}(1-z)^{s(\sba)-1} \cdot y_2^{(0)}(\ba^*; z),  
\qquad |z| < 1,   
\end{align*}
%%%%%%%%%%
where $\ba^*$ is defined in Table \ref{tab:invol}, 
while Okubo, Takano and Yoshida \cite[Lemma 2]{OTY} shows that 
%%%%%%%%%%
\[
\lim_{z \uparrow 1} (1-z)^{s(\sba)-1} \cdot y_2^{(0)}(\ba^*; z) = \vG(s(\ba)-1), 
\qquad \rRe \, s(\ba) > 1. 
\]
%%%%%%%%%%% 
It follows from these facts that $w(\ba) := w(\ba; 1)$ admits a representation 
%%%%%%%%%%%
\[
w(\ba) = 
\dfrac{\vG(a_0) \vG(a_1) \vG(a_2) \vG(a_0-b_1+1) \vG(a_1-b_1+1) \vG(a_2-b_1+1) 
\vG(s(\ba)-1)}{\vG(b_1) \vG(1-b_1) \vG(b_2-a_0)\vG(b_2-a_1) \vG(b_2-a_2)}.      
\]
%%%%%%%%%%
By the connection formula 
$y_0^{(\infty)}(\ba) = C_0(\ba) \, y_0^{(0)}(\ba) + C_1(\ba) \, y_1^{(0)}(\ba)$ in 
\cite[formula (16)]{EI}, where    
%%%%%%%%%%
\[
C_0(\ba) = \dfrac{e^{\mathrm{i} \pi s(\sba)} \cdot 
\sin \pi a_1 \cdot \sin \pi a_2}{\sin \pi b_1 \cdot \sin \pi(b_2-a_0)}, 
\qquad
C_1(\ba) = 
-\dfrac{e^{\mathrm{i} \pi s(\sba)} \cdot 
\sin \pi (b_1-a_1) \cdot \sin \pi (b_1-a_2)}{\sin \pi b_1 \cdot \sin \pi(b_2-a_0)},     
\]
%%%%%%%%%% 
and the periodicity $C_i(\ba+\1) = C_i(\ba)$, $i = 0, 1$, we have  
$W(\ba) = C_1(\ba) \, w(\ba)$. 
This together with  the reflection formula for the gamma function yields 
formula \eqref{eqn:wronskian}. \hfill $\Box$ 
%%%%%%%%%%%%%%%%%%%%%%%%%%%%%% end proof %%%%%%%%%%%%%%%%%%%%%%%%%%%%%%%
%%%%%%%%%%%%%%%%%%%%%%%%%%%%%% thm:casorati %%%%%%%%%%%%%%%%%%%%%%%%%%%%%
\begin{theorem} \label{thm:casorati} 
The initial value of the Casoratian $\omega(n)$ is given by 
%%%%%%%%%%%%%%%%%%%%%%%%%%%%%% eqn:casorati %%%%%%%%%%%%%%%%%%%%%%%%%%%%%
\begin{equation} \label{eqn:casorati} 
\omega(0) = \dfrac{\pi^2 \cdot \rho(\ba; \bk) \cdot \vG(a_0) \vG(a_1) \vG(a_2) 
\vG( s(\ba) ) }{ t(\ba) \prod_{i=0}^2 
\prod_{j=1}^2 \vG(b_j-a_i +(l_j-k_i)_+)}, 
\end{equation}
%%%%%%%%%%%%%%%%%%%%%%%%%%%%%%%%%%%%%%%%%%%%%%%%%%%%%%%%%%%%%%%%%%%%%%
where $\rho(\ba; \bk) \in \Q[\ba]$ is the polynomial in \eqref{eqn:rho} 
and $t(\ba) := \sin \pi(b_1-a_0) \cdot \sin \pi(b_2-a_0)$.  
\end{theorem}
%%%%%%%%%%%%%%%%%%%%%%%%%%%%%%%%%%%%%%%%%%%%%%%%%%%%%%%%%%%%%%%%%%%%%%%  
%%%%%%%%%%%%%%%%%%%%%%%%%%%%%% begin proof %%%%%%%%%%%%%%%%%%%%%%%%%%%%%%
{\it Proof}.  
From definitions \eqref{eqn:f(n)} and \eqref{eqn:gg(n)} we find that  
%%%%%%%%%%%%
\begin{align*}
\omega(0) 
&= f_0(0) \cdot g_1(0) - f_1(0) \cdot g_0(0) 
= {}_3f_2(\ba) \cdot {}_3g_2(\ba + \bk) -  {}_3f_2(\ba + \bk) \cdot {}_3g_2(\ba) \\[1mm]
&= y_0^{(0)}(\ba) \, e^{-\mathrm{i} \pi s(\sba+\sbk)} y_0^{(\infty)}(\ba+\bk) 
- y_0^{(0)}(\ba+\bk) \, e^{-\mathrm{i} \pi s(\sba)} y_0^{(\infty)}(\ba) \\
&=  e^{-\mathrm{i} \pi s(\sba)} \{ y_0^{(0)}(\ba) \, y_0^{(\infty)}(\ba+\bk) 
- y_0^{(0)}(\ba+\bk) \, y_0^{(\infty)}(\ba) \} \\[1mm]
&=  e^{-\mathrm{i} \pi s(\sba)} r(\ba; \bk) 
\{ y_0^{(0)}(\ba) \, y_0^{(\infty)}(\ba+\1) - y_0^{(0)}(\ba+\1) \, y_0^{(\infty)}(\ba) \} 
= e^{-\mathrm{i} \pi s(\sba)} r(\ba; \bk) \, W(\ba), 
\end{align*}
%%%%%%%%%%%%
where the fourth equality follows from $s(\bk) = 0$ and the fifth from the 
three-term relation 
%%%%%
\[
y_0^{(\nu)}(\ba+\bk) = r_1(\ba; \bk) \, y_0^{(\nu)}(\ba) + 
r(\ba; \bk) \, y_0^{(\nu)}(\ba+\1), \quad \nu = 0, \infty, 
\]
%%%%%
where $r_1(\ba; \bk)$ and $r(\ba; \bk)$ are the $(1, 1)$ and $(1, 2)$ entries 
of  the connection matrix $A(\ba; \bk)$ as in \cite[formulas (33) and (34)]{EI}.  
Using formula \eqref{eqn:wronskian} one has  
%%%%%%%%%%%%
\begin{align*} 
\omega(0) 
&= - r(\ba; \bk) \, \dfrac{\vG(a_0) \vG(a_1) \vG(a_2) \vG(a_0-b_1+1) 
\vG(a_0-b_2+1) \vG(s(\ba)-1)}{\vG(b_1-a_1) \vG(b_1-a_2) 
\vG(b_2-a_1) \vG(b_2-a_2) } \\[2mm]
&= - \dfrac{ \pi^2 \cdot r(\ba; \bk) \cdot \vG(a_0) \vG(a_1) \vG(a_2) 
\vG(s(\ba)-1) }{ t(\ba) \prod_{i=0}^2 \prod_{j=1}^2 \vG(b_j-a_i) }  \\[2mm]
&= \dfrac{ \pi^2 \cdot \rho(\ba; \bk) \cdot \vG(a_0) \vG(a_1) \vG(a_2) \cdot 
\{ s(\ba)-1 \} \vG(s(\ba)-1) }{ t(\ba)  
 \prod_{i=0}^2 \prod_{j=1}^2 (b_j - a_i; \, (l_j-k_i)_+) \cdot 
\prod_{i=0}^2 \prod_{j=1}^2 \vG(b_j-a_i) }  = \mbox{RHS of \eqref{eqn:casorati}},  
\end{align*}
%%%%%%%%%%%%
where the second equality follows from the reflection formula for the gamma 
function, the third from \eqref{eqn:rho} and the final one from the recursion 
formula for the gamma function. \hfill $\Box$ 
%%%%%%%%%%%%%%%%%%%%%%%%%%%%%% end proof %%%%%%%%%%%%%%%%%%%%%%%%%%%%%%
%%%%%%%%%%%%%%%%%%%%%%%%%%%%%% subsec:error %%%%%%%%%%%%%%%%%%%%%%%%%%%
\subsection{Error Estimates} \label{subsec:error}
%%%%%%%%%%%%%%%%%%%%%%%%%%%%%%%%%%%%%%%%%%%%%%%%%%%%%%%%%%%%%%%%%%%%%
We are now in a position to establish our main results in \S \ref{subsec:mrcf} by 
means of the general estimate \eqref{eqn:error} upon putting $X(n) = f(n)$ and 
$Y(n) = g(n)$. 
In this subsection, unless otherwise mentioned explicitly, Landau's symbols 
$O(\, \cdot \,)$ are uniform in any compact subset of 
%%%%%%%%%
\[
\A := \{ \ba \in \C^5 \,:\, \rRe \, s(\ba) > 0\}. 
\]
%%%%%%%%%%%%%%%%%%%% begin proof of thm:cf-straight %%%%%%%%%%%%%%%%%%%%%%
{\it Proof of Theorem $\ref{thm:cf-straight}$}. 
In the straight case in Definition \ref{def:admissible} the sequences in 
\eqref{eqn:f(n)} and \eqref{eqn:gg(n)} are given by 
$f(n) = {}_3f_2(\ba + n \bk)$ and $g(n) = {}_3g_2(\ba + n \bk)$ respectively.  
Under the assumption of Theorem \ref{thm:cf-straight} we can use   
Theorems \ref{thm:rs} and \ref{thm:ds} with $\bp$ replaced by $\bk$ 
to get    
%%%%%%%%%%%%%%
\begin{align*}
f(n) &= {}_3f_2(\ba + n \bk) = \vG(s(\ba)) \cdot s_2(\bk)^{-s(\sba)} \cdot 
n^{-2 s(\sba)} \cdot \{1 + O(1/n) \}, \\[1mm]
g(n) &= {}_3g_2(\ba + n \bk) = \dfrac{B(\ba; \bk)}{t(\ba)} \cdot D(\bk)^n \cdot 
n^{-s(\sba)-\frac{1}{2} } \cdot \left\{ 1 + O(n^{-\frac{1}{2} }) \right\},  
\end{align*}
%%%%%%%%%%%%%
where $D(\bk)$, $t(\ba)$ and $B(\ba; \bk)$ are given by \eqref{eqn:A(p)} and 
\eqref{eqn:B(a;p)} with $\bp$ replaced by $\bk$, and hence 
%%%%%%%%%%%%%
\[
R(n) = \frac{f(n+2)}{g(n+2)} = 
\dfrac{ t(\ba) \cdot \vG(s(\ba)) \cdot s_2(\bk)^{-s(\sba)}}{ B(\ba; \bk) \cdot D(\bk)^2 } 
\cdot D(\bk)^{-n} \cdot n^{-s(\sba)+ \frac{1}{2} } \cdot 
\left\{1 + O(n^{-\frac{1}{2}}) \right\}.  
\]
%%%%%%%%%%%% 
\par
%%%%%%%%%%%%
Combining this formula with \eqref{eqn:casorati} in Theorem 
\ref{thm:casorati}, we have  
%%%%%%%%%%%%
\[
\omega(0) \, R(n) = \rho(\ba; \bk) \cdot e_{\rs}(\ba; \bk) \cdot 
\gamma(\ba; \bk) \cdot D(\bk)^{-n} \cdot n^{-s(\sba) + \frac{1}{2} } 
\cdot \left\{1 + O(n^{-\frac{1}{2} }) \right\}, 
\]
%%%%%%%%%%%%
where $e_{\rs}(\ba; \bk)$ and $\gamma(\ba; \bk)$ are defined in  
\eqref{eqn:es(a;k)} and \eqref{eqn:gamma(a)}. 
To cope with the error term in formula \eqref{eqn:error}, we also need to 
care about how $R(n) \cdot Y(0)/ X(0)$ depends on $\ba \in \A$.  
Observe that  
%%%%%%%%%%%
\begin{gather*}
\frac{R(n) \cdot Y(0)}{X(0)} = \frac{R(n) \cdot {}_3g_2(\ba)}{{}_3f_2(\ba)} = 
\psi_1(\ba) \cdot \psi_2(\ba) 
\cdot D(\bk)^{-n} \cdot n^{-s(\sba) + \frac{1}{2}} \cdot 
\left\{ 1 + O(n^{-\frac{1}{2} }) \right\}, \\ 
\mbox{with} \quad 
\psi_1(\ba) := \frac{\vG(s(\ba)) \cdot s_2(\bk)^{-s(\sba)}}{B(\ba; \bk) 
\cdot D(\bk)^2},  \quad  
\psi_2(\ba) := \frac{t(\ba) \cdot {}_3g_2(\ba)}{{}_3f_2(\ba)}.     
\end{gather*}
%%%%%%%%%%% 
It is obvious that $\psi_1(\ba)$ is holomorphic in $\A$. 
It is also easy to see that $\psi_2(\ba)$ is holomorphic in 
$\A_0 := \{\ba \in \A : {}_3f_2(\ba) \neq 0 \}$. 
Indeed ${}_3g_2(\ba)$ has a pole when $a_0-b_1+1 \in \Z_{\le 0}$ or 
$a_0-b_2+1 \in \Z_{\le 0}$ but the pole is canceled by a zero of $t(\ba) = 
\sin \pi(b_1-a_0) \cdot \sin \pi(b_2-a_0)$; similarly ${}_3g_2(\ba)$ has a 
pole when $a_0 \in \Z_{\le 0}$ but it is canceled by a pole of 
${}_3f_2(\ba)$.         
Now estimate \eqref{eqn:error} leads to asymptotic formula 
\eqref{eqn:cf-straight}, in which Landau's symbol is uniform in 
any compact subset of $\A_0$. \hfill $\Box$ 
%%%%%%%%%%%%%%%%%%%% end proof %%%%%%%%%%%%%%%%%%%%%%%%%%%%%%%%%%%%%
\par\medskip\noindent
%%%%%%%%%%%%%%%%%%%% begin proof of thm:cf-cyclic %%%%%%%%%%%%%%%%%%%%%%
{\it Proof of Theorem $\ref{thm:cf-cyclic}$}. 
In the twisted case in Definition \ref{def:admissible}, 
if $n = 3 m + i$, $m \in \Z_{\ge0}$, $i = 0, 1, 2$, then the sequences $f(n)$ 
in \eqref{eqn:f(n)} and $g(n)$ in \eqref{eqn:gg(n)} are given by  
%%%%%%%%%%%
\[
f(3 m + i) = {}_3f_2(\ba + \bj_i + m \bp), \qquad 
g(3 m + i) = {}_3g_2(\ba + \bj_i + m \bp),  
\] 
%%%%%%%%%%%
where $\bj_0 = \0$, $\bj_1 = \bk$, $\bj_2 = \bl$ and $\bp$ is the 
shift vector in formula \eqref{eqn:adm-kp}.  
%%%%%%
\par
%%%%%%
Observe that $\bp$ belongs to $\cS(\Z)$ and satisfies condition 
\eqref{eqn:cond-key},  if and only if the seed vector $\bk \in \Z^5$ 
fulfills condition \eqref{eqn:cond-cyclic}.  
Indeed, since $p_0 = p_1 = p_2 = l_1+l_2$, $q_1 = 3 l_1$ and $q_2 = 3 l_2$, the 
inequalities in \eqref{eqn:bp} becomes 
$l_1 + l_2 < 3 l_1 \le 3 l_2 < 2(l_1 + l_2)$, which is equivalent to $l_1 \le l_2 < 2 l_1$. 
Case (a) in condition \eqref{eqn:cond-key} now reads 
%%%%%%%
$\vD(\bp) = -27(l_2^2-4 l_1 l_2 + l_1^2) (l_2^2+2 l_1 l_2 -2l_1^2) 
(2 l_2^2 -2l_1 l_2-l_1^2) \le 0$,  
%%%%%%% 
which together with $l_1 \le l_2 < 2 l_1$ yields $l_1 \le l_2 \le \tau l_1$ in condition 
\eqref{eqn:cond-cyclic}, where $\tau = (1+\sqrt{3})/2$.    
On the other hand, case (b) in condition \eqref{eqn:cond-key} becomes   
$l_2^2 - 10 l_1 l_2 + 7 l_1^2 \ge 0$, that is, $l_2 \le (5-3 \sqrt{2}) l_1$ or 
$l_2 \ge (5 + 3\sqrt{2}) l_1$, but neither of which is possible when 
$l_1 \le l_2 < 2 l_1$.  
%%%%%%%
\par
%%%%%%%
Thus under the assumption of Theorem \ref{thm:cf-cyclic} one can apply  
Theorems \ref{thm:rs} and \ref{thm:ds} to the shift vector $\bp$ in 
\eqref{eqn:adm-kp} with $\ba$ replaced by $\ba+\bj_i$ to obtain    
%%%%%%%%%%% 
\begin{align*}
f(3 m + i) &= {}_3f_2(\ba + \bj_i + m \bp) = \vG(s(\ba)) \cdot s_2(\bp)^{-s(\sba)} 
\cdot m^{-2 s(\sba)} \cdot \{ 1 + O(1/m)\},  \\[1mm]  
g(3 m + i) &=  {}_3g_2(\ba + \bj_i + m \bp) = 
\dfrac{B(\ba + \bj_i ; \bp)}{t(\ba + \bj_i)} \cdot D(\bp)^m 
\cdot m^{-s(\sba)-\frac{1}{2}} 
\cdot  \left\{ 1 + O(m^{-\frac{1}{2} }) \right\},   
\end{align*}
%%%%%%%%%%% 
where $s(\bj_i) = 0$ is also used. 
Substituting the settings \eqref{eqn:adm-kp} and \eqref{eqn:cyclic} into 
definitions \eqref{eqn:A(p)} and \eqref{eqn:B(a;p)} and taking $s(\bk) = 0$ 
into account, one has $D(\bp) = E(l_1, l_2)^3$ and  
%%%%%%%%%%%
\[
t(\ba + \bj_i) = (-1)^{i(l_1+l_2)} \cdot t(\ba), \quad 
B(\ba+ \bj_i; \bp) = (-1)^{i (l_1+l_2)} \cdot 
B(\ba; \bp) \cdot E(l_1, l_2)^i, \quad i = 0, 1, 2,   
\]
%%%%%%%%%%
where $E(l_1, l_2)$ is defined in \eqref{eqn:A(l1,l2)}.  
These formulas and $s_2(\bp) = 3 (l_1^2 -l_1l_2+l_2^2)$ lead to 
%%%%%%%%%%
\begin{align*}
f(n) 
&=  3^{s(\sba)} \cdot \vG(s(\ba)) \cdot (l_1^2 -l_1 l_2 + l_2^2)^{-s(\sba)}\cdot  
n^{-2 s(\sba)} \cdot \{1 + O(1/n)\},  \\[1mm]  
g(n) &= 3^{s(\sba)+ \frac{1}{2} } \cdot \frac{B(\ba; \bp)}{t(\ba)} 
\cdot E(l_1, l_2)^n \cdot n^{-s(\sba)- \frac{1}{2} } \cdot 
\left\{ 1 + O(n^{-\frac{1}{2}}) \right\},   
\end{align*}
%%%%%%%%%%%
so the ratio $R(n) = f(n+2)/g(n+2)$ is estimated as 
%%%%%%%%%%
\[
R(n) = \dfrac{ t(\ba) \cdot \vG(s(\ba)) \cdot 
(l_1^2-l_1l_2+l_2^2)^{-s(\sba)} }{ 3^{\frac{1}{2} } \cdot 
B(\ba; \bp) \cdot E(l_1, l_2)^2} \cdot 
E(l_1, l_2)^{-n} \cdot n^{-s(\sba)+ \frac{1}{2} } \cdot 
\left\{1 + O(n^{-\frac{1}{2} }) \right\}.  
\]
%%%%%%%%%%%
\par
%%%%%%%%%%%
Substituting $\bp = (l_1+l_2, l_1+l_2, l_1+l_2; 3l_1, 3l_2)$ into definition 
\eqref{eqn:B(a;p)} yields 
%%%%%%%%%%%
\[
B(\ba; \bp) = 
\dfrac{\pi^{ \frac{1}{2} } (l_1+l_2)^{a_0+a_1+a_2- \frac{3}{2} } \cdot 3^{s(\sba)-1} 
\cdot (l_1^2-l_1l_2+l_2^2)^{s(\sba)-1}}{2^{ \frac{3}{2} } \cdot 
(2l_1-l_2)^{2 b_1-b_2+s(\sba)- \frac{3}{2} } 
(2l_2-l_1)^{2 b_2-b_1+s(\sba)-\frac{3}{2} }},  
\]
%%%%%%%%%%
which is put together with formula \eqref{eqn:casorati} in Theorem \ref{thm:casorati} 
to give   
%%%%%%%%%%%%
\[
\omega(0) \, R(n) = \rho(\ba; \bk) \cdot e_{\rt}(\ba; \bk) \cdot \gamma(\ba; \bk) 
\cdot E(l_1, l_2)^{-n} \cdot n^{-s(\sba) + \frac{1}{2} } \cdot 
\left\{1 + O(n^{-\frac{1}{2} }) \right\}, 
\]
%%%%%%%%%%%%
where $e_{\rt}(\ba; \bk)$ and $\gamma(\ba; \bk)$ are given in 
\eqref{eqn:ec(a;k)} and \eqref{eqn:gamma(a)}.  
The treatment of $R(n) \cdot Y(0)/X(0)$ is similar to the one in the straight case 
and the estimate \eqref{eqn:error} leads to  asymptotic formula 
\eqref{eqn:cf-cyclic}, in which Landau's symbol is uniform in any compact 
subset of $\A_0$. \hfill $\Box$ 
%%%%%%%%%%%%%%%%%%%%%%%% end proof %%%%%%%%%%%%%%%%%%%%%%%%%%%%%%%%%%%%%%
%%%%%%%%%%%%%%%%%%%%%%%% sec:rn-s %%%%%%%%%%%%%%%%%%%%%%%%%%%%%%%%%%%%%%%
\section{Back to Original Series and Specializations} \label{sec:bos-s}
%%%%%%%%%%%%%%%%%%%%%%%%%%%%%%%%%%%%%%%%%%%%%%%%%%%%%%%%%%%%%%%%%%%%%%%
Theorems \ref{thm:cf-straight} and \ref{thm:cf-cyclic} are stated in terms of 
the renormalized series ${}_3f_2(\ba)$. 
It is interesting  
to reformulate them in terms of the original series ${}_3F_2(1)$. 
Multiplying equations \eqref{eqn:cf-straight} and \eqref{eqn:cf-cyclic} by   
%%%%%%%%%
\[
\frac{\vG(b_1+l_1) \vG(b_2+l_2)}{\vG(a_0+k_0) \vG(a_1+k_1) \vG(a_2+k_2)} 
\cdot \frac{\vG(a_0) \vG(a_1) \vG(a_2)}{\vG(b_1) \vG(b_2)} = 
\frac{(b_1; \, l_1)(b_2; \, l_2)}{(a_0; \, k_0) (a_1; \, k_1) (a_2; \, k_2)}
\]
%%%%%%%%
and using relation \eqref{eqn:Fvsf} between ${}_3f_2(\ba)$ and ${}_3F_2(\ba)$, 
we find that 
%%%%%%%%%%%%%%%%%%%%%%%%%%%%%% eqn:cf-o %%%%%%%%%%%%%%%%%%%%%%%%%%%%%%%
\begin{subequations} \label{eqn:cf-o}
\begin{align} 
\frac{{}_3F_2(\ba+\bk)}{{}_3F_2(\ba)} - 
\overset{n}{\underset{j=0}{\cfL}} \,\, \frac{r^*(j)}{q^*(j)} 
&= \dfrac{ c_{\rs}^*(\ba; \bk) }{{}_3F_2(\ba)^2}  
\cdot D(\bk)^{-n} \cdot n^{-s(\sba) + \frac{1}{2} } \cdot 
\left\{ 1 + O (n^{-\frac{1}{2} } ) \right\},  \label{eqn:cf-straight-o} \\[2mm] 
\frac{{}_3F_2(\ba+\bk)}{{}_3F_2(\ba)} - 
\overset{n}{\underset{j=0}{\cfL}} \,\, \frac{r^*(j)}{q^*(j)} 
&= \dfrac{ c_{\rt}^*(\ba; \bk) }{{}_3F_2(\ba)^2}  
\cdot E(l_1, l_2)^{-n} \cdot n^{-s(\sba) + \frac{1}{2} } \cdot 
\left\{ 1 + O (n^{-\frac{1}{2} } ) \right\},  \label{eqn:cf-cyclic-o}
\end{align}
\end{subequations}
%%%%%%%%%%%%%%%%%%%%%%%%%%%%%%%%%%%%%%%%%%%%%%%%%%%%%%%%%%%%%%%%%%%%%
as $n \to + \infty$, where the quantities marked with an asterisk are defined by 
%%%%%%%%%%%%%%%%
\begin{alignat*}{3}
q^*(0)  &:= u(\ba) \ts \prod_{i=0}^2 (a_i; \, k_i), \qquad & 
q^*(n)  &:= q(n), \quad n \ge 1, \qquad &  & \\[2mm]
r^*(0)  &:= (b_1; \, l_1)(b_2; \, l_2), \qquad & 
r^*(1) &:= v(\ba) \ts \prod_{i=0}^2 (a_i; \, k_i), \qquad & 
r^*(n) &:= r(n), \quad n \ge 2, 
\end{alignat*}
%%%%%%%%%%%%%%%%%%%%%%%%%%% eqn:c* %%%%%%%%%%%%%%%%%%%%%%%%%%%%%%%%%
\begin{equation} \label{eqn:c*}
c_{\iota}^*(\ba; \bk)  := \rho(\ba; \bk) \cdot e_{\iota}(\ba; \bk) \cdot 
\gamma^*(\ba; \bk), \qquad  \iota = \rs, \rt,  
\end{equation}
%%%%%%%%%%%%%%%%%%%%%%%%%%%%%%%%%%%%%%%%%%%%%%%%%%%%%%%%%%%%%%%%%%%
with $\rho(\ba; \bk)$ and $e_{\iota}(\ba; \bk)$ unaltered while  
%%%%%%
\[
\gamma^*(\ba; \bk) := \frac{\vG(b_1+l_1) \vG(b_2+l_2) \vG(b_1) \vG(b_2) 
\vG^2(s(\ba)) }{\vG(a_0+k_0) \vG(a_1+k_1) \vG(a_2+k_2) 
\prod_{i=0}^2 \prod_{j=1}^2 \vG(b_j-a_i+(l_j-k_i)_+)}. 
\]
%%%%%% 
It follows from \eqref{eqn:Fvsf} that $\A_0^* 
:= \{ \ba \in \A \,:\, b_1, \, b_2 \not \in \Z_{\le 0}, \,\, 
{}_3 F_2 (\ba) \neq 0 \} \subset \A_0$, where $\A$ and $\A_0$ are 
defined in \S \ref{subsec:error}, so Landau's symbols in \eqref{eqn:cf-o} 
are uniform in any compact subset of $\A_0^*$.  
%%%%%%
\par
%%%%%% 
Take an index $\lambda \in \{0, 1, 2\}$ such that $k_{\lambda} > 0$ 
and put $\{\lambda, \mu, \nu\} = \{0, 1, 2\}$.  
For any nonzero vector $\bk \in \Z_{\ge0}^5$ with $s(\bk) =0$ such an 
index $\lambda$ always exists since $k_0+k_1+k_2 = l_1+l_2 > 0$. 
In formulas \eqref{eqn:cf-o} take the limit $a_{\lambda} \to 0$ and 
make the substitutions $a_i \mapsto a_i-k_i$, $b_j \mapsto b_j-l_j$ 
for $i = \mu, \nu$ and $j = 1, 2$.  
This procedure is referred to as the $\lambda$-th {\sl specialization}. 
If this is well defined then ${}_3F_2(\ba) \to 1$ as 
$a_{\lambda} \to 0$ so formulas \eqref{eqn:cf-straight-o} and 
\eqref{eqn:cf-cyclic-o} lead to     
%%%%%%%%%%%%%%%%%%%%%%%%%%%%%% eqn:cf-s %%%%%%%%%%%%%%%%%%%%%%%%%%%%%
\begin{subequations} \label{eqn:cf-s}
\begin{align} 
{}_3F_2\! 
\begin{pmatrix} 
k_{\lambda}, & a_{\mu}, & a_{\nu} \\ 
     & b_1, & b_2    
\end{pmatrix}
- \overset{n}{\underset{j=0}{\cfL}} \,\, \frac{ \hat{r}(j) }{ \hat{q}(j) } 
&= \hat{c}_{\rs}(\ba; \bk)    
\cdot D(\bk)^{-n} \cdot n^{- \hat{s} + \frac{1}{2} } \cdot 
\left\{ 1 + O (n^{-\frac{1}{2} } ) \right\},  \label{eqn:cf-straight-s} \\[2mm] 
{}_3F_2\! 
\begin{pmatrix} 
k_{\lambda}, & a_{\mu}, & a_{\nu} \\ 
     & b_1, & b_2    
\end{pmatrix} 
 - \overset{n}{\underset{j=0}{\cfL}} \,\, \frac{ \hat{r}(j)}{ \hat{q}(j) } 
&= \hat{c}_{\rt}(\ba; \bk) \cdot E(l_1, l_2)^{-n} \cdot n^{- \hat{s} + \frac{1}{2} } 
\cdot \left\{ 1 + O (n^{-\frac{1}{2} } ) \right\},  \label{eqn:cf-cyclic-s}
\end{align}
\end{subequations} 
%%%%%%%%%%%%%%%%%%%%%%%%%%%%%%%%%%%%%%%%%%%%%%%%%%%%%%%%%%%%%%%%%%%%%
where $\hat{q}(n)$ and $\hat{r}(n)$ are derived from $q^*(n)$ and $r^*(n)$,  
while $\hat{c}_{\iota}(\ba; \bk) := \hat{\rho}(\ba; \bk) \cdot 
\hat{e}_{\iota}(\ba; \bk) \cdot \hat{\gamma}(\ba; \bk)$, $\iota = \rs, \rt$,    
are obtained from \eqref{eqn:c*} through the specialization; in particular one has    
%%%%%%
\[
\hat{\gamma}(\ba; \bk) :=  
\frac{\vG(b_1) \vG(b_2) \vG^2(\hat{s})}{\vG(k_{\lambda}) 
\vG(a_{\mu}) \vG(a_{\nu}) 
\ds \prod_{j=1,2} (b_j-l_j; \, (l_j-k_{\lambda})_+) \cdot 
\prod_{i = \mu, \nu} \prod_{j=1,2} \vG(b_j-a_i + (k_i-l_j)_+) }, 
\]
%%%%%% 
with $\hat{s} := b_1+b_2-a_{\mu}-a_{\nu}-k_{\lambda}$. 
Landau's symbols in \eqref{eqn:cf-s} are uniform in compact subsets of 
%%%%%%
\[
\hat{\A} := \{\, (a_{\mu}, a_{\nu}; b_1, b_2) \in \C^4 \,:\, \rRe \, \hat{s} > 0, \, \, 
b_1, b_2 \not \in \Z_{\le 0} \, \}. 
\]
%%%%%%
\par
%%%%%%
The specialization is indeed well defined. 
It follows from \cite[Proposition 4.9]{EI} and Lemma \ref{lem:factorial} that 
for any $\bc \in \Q^5$ the restriction 
$\rho(\ba+\bc; \bk) \big|_{a_0 = a_1 = a_2 = 0}$ is a nonzero polynomial in 
$\Q[b_1, b_2]$ and hence $\rho(\ba+\bc; \bk)|_{a_{\lambda} = 0}$ is a 
nonzero polynomial in $\Q[a_{\mu}, a_{\nu}, b_1, b_2]$. 
This is also the case for $\sigma(\bk)$ and $\bl$ in place of $\bk$ in 
\S \ref{subsec:types}. 
Thus formula \eqref{eqn:uv2u} implies that the specialization for $q^*(0)$, 
%%%%%%
\[
\hat{q}(0) := \lim_{a_{\lambda} \to 0} u(\ba) \prod_{i=0}^2 (a_i; \, k_i) 
\quad \mbox{followed by} \quad a_i \mapsto a_i - k_i, \,\, 
b_j \mapsto b_j - l_j, \quad i = \mu, \nu, \,\, j = 1, 2,    
\] 
%%%%%
is well defined and the ensuing $\hat{q}(0)$ is a nontrivial rational function 
in $\Q(a_{\mu}, a_{\nu}, b_1, b_2)$. 
In a similar manner formula \eqref{eqn:uv2v} tells us that the specialization 
for $r^*(1)$, that is, $\hat{r}(1) \in \Q(a_{\mu}, a_{\nu}, b_1, b_2)$ is  
well defined and nontrivial. 
The specialization for $q^*(n)$ with $n \ge 1$ is also well defined,  
since $q^*(n) = q(n)$ is of the form ${}^{\sigma^i} \! u(\ba + \bc)$, where 
$i \in \{0,1,2\}$ and $\bc$ is a vector in $\Z_{\ge 0}^5$ whose 
$\lambda$-th upper component, say $c_{\lambda}$, is positive, in which case 
one can take $\lim_{a_{\lambda} \to 0} {}^{\sigma^i} \! u(\ba + \bc)$ 
without trouble, because the critical factorial 
$(a_{\lambda}; \, k_{\lambda}) \to 0$ in the denominator of 
\eqref{eqn:uv2u} is now replaced by a safe one 
$(a_{\lambda}+ c_{\lambda}; \, k_{\lambda}) \to 
(c_{\lambda}; \, k_{\lambda}) \neq 0$. 
The resulting $\hat{q}(n)$ is nontrivial in $\Q(a_{\mu}, a_{\nu}, b_1, b_2)$.  
A similar argument can be made for $\hat{r}(n)$ with $n \ge 2$. 
Thus the procedure of specialization is well defined over the rational function 
field $\Q(a_{\mu}, a_{\nu}, b_1, b_2)$.        
%%%%%%%%%%%%%%%%%%%%%%%% sec:example %%%%%%%%%%%%%%%%%%%%%%%%%%%%%%%%%%
\section{Some Examples} \label{sec:example} 
%%%%%%%%%%%%%%%%%%%%%%%%%%%%%%%%%%%%%%%%%%%%%%%%%%%%%%%%%%%%%%%%%%%%%
To illustrate Theorems \ref{thm:cf-straight} and \ref{thm:cf-cyclic} we 
present a couple of the simplest examples.          
%%%%%%%%%%%%%%%%%%%%%%%%%%%%% ex:1 %%%%%%%%%%%%%%%%%%%%%%%%%%%%%%%%%%%
\begin{example} \label{ex:1} 
The simplest example of twisted type is given by  
%%%%%%
\[
\bk = \begin{pmatrix} 1, & 1, & 0 \\ & 1, & 1 \end{pmatrix}, \qquad 
\bl = \begin{pmatrix} 1, & 2, & 1 \\ & 2, & 2 \end{pmatrix}, \qquad 
\bp = \begin{pmatrix} 2, & 2, & 2 \\ & 3, & 3 \end{pmatrix},   
\]
%%%%%%
together with $\sigma(a_0, a_1, a_2; b_1, b_2) = (a_2, a_0, a_1; b_1, b_2)$.  
The recipe described in \S \ref{subsec:cr} readily yields   
$\rho(\ba; \bl) = b_1 b_2 - a_0 a_2$ and $\rho(\ba; \bk) = 
\rho(\ba+\bk; \sigma(\bk))  = 1$, so formula \eqref{eqn:uv2} yields   
%%%%%%%%
\[
u(\ba) = \dfrac{ b_1 b_2 - a_0 a_2}{a_0 a_1}, \qquad 
v(\ba) = -\dfrac{(b_1-a_0)(b_2-a_0)}{a_0 a_1}.      
\] 
%%%%%
\par
%%%%%
Thus the partial denominators and numerators of the continued fraction in 
Theorem \ref{thm:cf-cyclic} are given as in Table \ref{tab:pdn-ex1}, 
%%%%%%%%%%%%%%%%%%%%%%%%%%%%%%%% tab:pdn-ex1 %%%%%%%%%%%%%%%%%%%%%%%%%%%%%
\begin{table}[t]
%%%%%%%%%%%%%
\begin{alignat*}{2} 
q_0(n) &:= \dfrac{(3 n+b_1)(3 n+b_2)-(2 n + a_0)(2 n+ a_2)}{(2 n + a_0)(2 n + a_1)} 
\qquad & (n &\ge 0),  \\[2mm]
q_1(n) &:= \dfrac{(3 n + b_1 +1)(3 n + b_2+1)-(2 n+a_0+1)(2 n+ a_1+1)}{(2 n + a_1+1)(2 n + a_2)} 
\qquad & (n &\ge 0),  \\[2mm]
q_2(n) &:= 
\dfrac{(3 n+b_1+2)(3 n+ b_2+2)-(2 n+ a_1+2)(2 n+ a_2+1)}{(2 n+a_0+1)(2 n + a_2 + 1)} 
\qquad & (n &\ge 0),  \\[2mm] 
r_0(n) &:= - \dfrac{(n+b_1-a_2)(n+b_2-a_2)}{(2 n+a_0-1)(2 n+a_2-1)} 
\qquad & (n &\ge 1), \\[2mm]
r_1(n) &:= - \dfrac{(n+b_1-a_0)(n+b_2-a_0)}{(2 n+a_0)(2 n+a_1)} 
\qquad & (n &\ge 0),  \\[2mm]
r_2(n) &:= - \dfrac{(n+b_1-a_1)(n+b_2-a_1)}{(2 n+a_1+1)(2 n+a_2)} 
\qquad & (n &\ge 0).   
\end{alignat*}
\caption{Partial denominators and numerators in Example \ref{ex:1} with $r_0(0) := 1$.} 
\label{tab:pdn-ex1}
\end{table} 
%%%%%%%%%%%%%%%%%%%%%%%%%%%%%%%%%%%%%%%%%%%%%%%%%%%%%%%%%%%%%%%%%%%%%%%
so the continued fraction for ${}_3 f_2(\ba+\bk)/{}_3f_2(\ba)$ is well defined when 
%%%%%%%%%%%%
\[
a_0, \, a_1, \, a_2 \not \in \Z_{\le 0}; \qquad 
b_j -a_i \not\in \Z_{\le 0}, \quad b_j-a_2 \not\in \Z_{\le -1}, \quad  
i = 0, 1, \,\, j = 1, 2.    
\]
%%%%%%%%%%% 
In the error estimate \eqref{eqn:cf-cyclic} we have $E(l_1, l_2) = E(1,1) = 4$ and 
%%%%%%%%%%%
\[
c_{\rt}(\ba; \bk) = \dfrac{\pi^{\frac{3}{2}} \cdot 
\vG(a_0) \vG(a_1) \vG(a_2) \vG^2(s(\ba))}{2^{a_0+a_1+a_2+1} \cdot 3^{s(\sba) - \frac{1}{2}} 
\prod_{j=1}^2 \prod_{i=0}^2 \vG(b_j-a_i + \delta_{i 2})}.    
\]
%%%%%%%%%%% 
\par
%%%%%%
Passing to the continued fraction for ${}_3 F_2(\ba+\bk)/{}_3F_2(\ba)$, 
we have $q_0^*(0) = b_1 b_2-a_0 a_2$, $r^*_0(0) = b_1 b_2$,  
$r^*_1(0) = -(b_1-a_0)(b_2-a_0)$ and $q^*_i(n) = q_i(n)$, $r^*_i(n) = r_i(n)$ 
for all other $(i, n)$ in formula \eqref{eqn:cf-cyclic-o}.   
The $0$-th specialization of \eqref{eqn:cf-cyclic-o} then leads to 
$\hat{q}_0(0) = \hat{r}_0(0) = (b_1-1)(b_2-1)$, $\hat{r}_1(0) = -(b_1-1)(b_2-1)$ 
in formula \eqref{eqn:cf-cyclic-s} while all other $\hat{q}_i(n)$ and 
$\hat{r}_i(n)$ are given as in Table \ref{tab:pdn-cf}, where circumflex 
``$\hat{\phantom{q}}$'' is deleted for the sake of simplicity. 
Clearly we can make $\hat{q}_0(0) = \hat{r}_0(0) = 1$, $\hat{r}_1(0) = -1$ 
up to equivalence of continued fractions and we have established      
Theorem \ref{thm:ete-i}. 
\end{example}
%%%%%%%%%%%%%%%%%%%%%%%%%%%%%%%%%%%%%%%%%%%%%%%%%%%%%%%%%%%%%%%%%%%%%%
%%%%%%%%%%%%%%%%%%%%%%%% ex:2 %%%%%%%%%%%%%%%%%%%%%%%%%%%%%%%%%%%%%%%%%
\begin{example} \label{ex:2} 
The next simplest example of twisted type is given by  
%%%%%%
\[
\bk = \begin{pmatrix} 2, & 0, & 0 \\ & 1, & 1 \end{pmatrix}, \qquad 
\bl = \begin{pmatrix} 2, & 2, & 0 \\ & 2, & 2 \end{pmatrix}, \qquad 
\bp = \begin{pmatrix} 2, & 2, & 2 \\ & 3, & 3 \end{pmatrix},    
\]
%%%%%% 
where $\bp$ and $\sigma$ are the same as in Example \ref{ex:1}. 
In this case the recipe in \S\ref{subsec:cr} gives 
%%%%%%
\begin{align*}
\rho(\ba; \bl) &= b_1 b_2 + a_0 a_1 - (a_2 - 1) (a_0 + a_1 + 1), \\
\rho(\ba; \bk) &= b_1 b_2 -a_1 a_2 - (a_0 + 1)(b_1+b_2 -a_1-a_2), \\ 
\rho(\ba+\bk; \sigma(\bk)) &= (b_1+1)(b_2+1)-(a_0+2)a_2-(a_1+1) (b_1+b_2-a_0-a_2).   
\end{align*}
%%%%%%
With these data the formula \eqref{eqn:uv2} yields   
%%%%%%
\[ 
u(\ba) = \frac{(b_1 - a_1) (b_2 - a_1) \cdot \rho(\ba; \bl) }{a_0 (a_0 + 1) \cdot 
\rho(\ba+\bk; \sigma(\bk)) }, \qquad 
v(\ba) = - \frac{(b_1 - a_2 + 1) (b_2 - a_2 + 1) \cdot 
\rho(\ba; \bk) }{a_0 (a_0 + 1) \cdot \rho(\ba+\bk; \sigma(\bk)) }. 
\]
%%%%%% 
In the error estimate \eqref{eqn:cf-cyclic} in Theorem \ref{thm:cf-cyclic} we have 
$E(l_1, l_2) = E(1,1) = 4$ and 
%%%%%%
\[
c_{\rt}(\ba; \bk) = 
\dfrac{\pi^{\frac{3}{2}} \cdot \rho(\ba; \bk) \cdot \vG(a_0) \vG(a_1) \vG(a_2) 
\vG^2(s(\ba))}{2^{a_0+a_1+a_2+1} \cdot 3^{s(\sba) - \frac{1}{2}} \cdot 
\prod_{j=1}^2\vG(b_j-a_0) \cdot \prod_{i=1}^2 \prod_{j=1}^2 \vG(b_j-a_i+1)}. 
\]
%%%%%%
The $0$-th specialization leads to a continued fraction expansion 
\eqref{eqn:cf-cyclic-s} for ${}_3F_2(2, a_1,a_2; b_1, b_2)$. 
\end{example}
%%%%%%%%%%%%%%%%%%%%%%%%%%%%%%%%%%%%%%%%%%%%%%%%%%%%%%%%%%%%%%%%%%%%%%
%%%%%%%%%%%%%%%%%%%%%%%% ex:3 %%%%%%%%%%%%%%%%%%%%%%%%%%%%%%%%%%%%%%%%%
\begin{example} \label{ex:3}
The simplest example of straight type is given by $\bk = (2,2,2; 3,3)$, 
$\bl = 2 \bk$ and $\bp = 3 \bk$. 
The recipe in \S \ref{subsec:cr} shows that 
$\rho(\ba; \bk) = a_0 a_1 a_2 (b_1+b_2+1) + b_1 b_2 \{s(\ba) - s_2(\ba) \}$ and   
$\rho(\ba; 2 \bk)$ is a polynomial of degree $10$ (explicit formula is omitted).  
Formula \eqref{eqn:uv3} yields 
%%%%%%
\[
u(\ba) = \dfrac{\rho(\ba; 2 \bk)}{\rho(\ba+\bk; \bk) \prod_{i=0}^2 a_i(a_i+1)}, 
\qquad 
v(\ba) = - \dfrac{\rho(\ba; \bk) 
\prod_{i=0}^2 \prod_{j=1}^2 (b_j-a_i+1)}{\rho(\ba+\bk; \bk) \prod_{i=0}^2 a_i(a_i+1)}. 
\]
%%%%%%
In the error estimate \eqref{eqn:cf-straight} in Theorem \ref{thm:cf-straight} 
we have $D(\bk) = 2^6 = 64$ and 
%%%%%%
\[
c_{\rs}(\ba; \bk) = 
\dfrac{\pi^{\frac{3}{2}} \cdot \rho(\ba; \bk) \cdot \vG(a_0) \vG(a_1) \vG(a_2) 
\vG^2(s(\ba))}{2^{a_0+a_1+a_2+9} \cdot 3^{2 s(\sba) - 1} \cdot 
\prod_{i=0}^2 \prod_{j=1}^2 \vG(b_j-a_i+1)}. 
\]
%%%%%%
\end{example} 
%%%%%%%%%%%%%%%%%%%%%%%% References %%%%%%%%%%%%%%%%%%%%%%%%%%%%%%%%%%%
\bibliographystyle{plain}
\bibliography{conti-frac-ref}
%%%%%%%%%%%%%%%%%%%%%%%%%%%%%%%%%%%%%%%%%%%%%%%%%%%%%%%%%%%%%%%%%%%%%%
\end{document}